\documentstyle []{article}
\def\picture#1by#2(#3){
\vbox to #2 {
   \hrule width #1 height 0pt depth 0pt \vfill \special{picture #3}}
}

\def\scaledpicture#1by#2(#3scaled#4){{
\dimen0=#1  \dimen1=#2
\divide\dimen0 by 1000 \multiply\dimen0 by #4
\divide\dimen1 by 1000 \multiply\dimen1 by #4
\picture \dimen0 by \dimen1 (#3 scaled #4)}}
\def\dfigure#1by#2(#3scaled#4offset#5:#6)
   {\medskip
    \vglue 2mm minus 2mm
    $$
      \hbox{
        \hglue#5
        {\scaledpicture #1 by #2 (#3 scaled #4)}
      }
    $$
    \par\goodbreak
    \vglue 2mm minus 2mm
    \medskip}

\def\qmod#1#2{{\hbox{}^{\displaystyle{#1}}}\!\big/\!\hbox{}_{
\displaystyle{#2}}}

\def\inw#1{\mathop{#1}\limits^\circ}

\def\qqmod#1#2{{\hbox{}^{\displaystyle{#1}}}\!\big/\hskip -3pt \big/\!\hbox{}_{
\displaystyle{#2}}}

\def\resto#1#2{{
#1\hskip 0.4ex\vline_{\hskip 0.4ex\raisebox{-1ex}
{{${\scriptstyle #2}$}}}}}

\newfam\msbfam

\font\tenmsb=msbm10
\font\eightmsb=msbm10 at 8pt
\font\sevenmsb=msbm10 at 7pt
\font\fivemsb=msbm10 at 5pt

\textfont\msbfam=\tenmsb
\scriptfont\msbfam=\sevenmsb
\scriptscriptfont\msbfam=\fivemsb

\def\Bbb{\fam\msbfam\tenmsb}
\font\curly=rsfs10

\def\A{{\Bbb A}}

\def\C{{\Bbb C}}

\def\E{{\Bbb E}}
\def\F{{\Bbb F}}
\def\G{{\Bbb G}}
\def\H{{\Bbb H}}

\def\N{{\Bbb N}}

\def\P{{\Bbb P}}

\def\R{{\Bbb R}}
\def\Z{{\Bbb Z}}

\def\cringle{\mathaccent23}
\def\union{\mathop{\bigcup}}
\def\qed {\hfill\vrule height6pt width6pt depth0pt \bigskip}

\def\map{\longrightarrow}
\def\textmap#1{\mathop{\vbox{\ialign{
                                 ##\crcr
     ${\scriptstyle\hfil\;\;#1\;\;\hfil}$\crcr
     \noalign{\kern-1pt\nointerlineskip}
     \rightarrowfill\crcr}}\;}}

\def\textlmap#1{\mathop{\vbox{\ialign{
                                 ##\crcr
     ${\scriptstyle\hfil\;\;#1\;\;\hfil}$\crcr
     \noalign{\kern-1pt\nointerlineskip}
     \leftarrowfill\crcr}}\;}}

\newfam\meuffam

\font\tenmeuf=eufm10
\font\sevenmeuf=eufm7
\font\fivemeuf=eufm5

\textfont\meuffam=\tenmeuf
\scriptfont\meuffam=\sevenmeuf
\scriptscriptfont\meuffam=\fivemeuf

\def\germ{\fam\meuffam\tenmeuf}

\def\ag{{\germ a}}
\def\bg{{\germ b}}
\def\cg{{\germ c}}
\def\dg{{\germ d}}
\def\eg{{\germ e}}

\def\g{{\germ g}}
\def\hg{{\germ h}}

\def\kg{{\germ k}}
\def\lg{{\germ l}}

\def\sg{{\germ s}}
\def\tg{{\germ t}}
\def\ug{{\germ u}}

\def\zg{{\germ z}}

\def\Rg{{\germ R}}

\newtheorem{sz}{Satz}[section]
\newtheorem{thry}[sz]{Theorem}
\newtheorem{pr}[sz]{Proposition}
\newtheorem{re}[sz]{Remark}
\newtheorem{co}[sz]{Corollary}
\newtheorem{dt}[sz]{Definition}
\newtheorem{lm}[sz]{Lemma}

\begin{document}
\def\Pr{{\rm Pr}}
\def\tr{{\rm Tr}}
\def\End{{\rm End}}
\def\Aut{{\rm Aut}}
\def\Spin{{\rm Spin}}
\def\U{{\rm U}}
\def\SU{{\rm SU}}
\def\SO{{\rm SO}}
\def\PU{{\rm PU}}
\def\GL{{\rm GL}}
\def\spin{{\rm spin}}
\def\u{{\rm u}}
\def\su{{\rm su}}
\def\so{{\rm so}}
\def\ub{\underbar}
\def\pu{{\rm pu}}
\def\Pic{{\rm Pic}}
\def\Iso{{\rm Iso}}
\def\NS{{\rm NS}}
\def\deg{{\rm deg}}
\def\Hom{{\rm Hom}}
\def\Aut{{\rm Aut}}
\def\h{{\germ h}}
\def\Herm{{\rm Herm}}
\def\Vol{{\rm Vol}}
\def\pf{{\bf Proof: }}
\def\id{{\rm id}}
\def\i{{\germ i}}
\def\im{{\rm im}}
\def\rk{{\rm rk}}
\def\ad{{\rm ad}}
\def\h{{\bf H}}
\def\coker{{\rm coker}}
\def\dbar{\bar{\partial}}
\def\Lo{{\Lambda_g}}
\def\niq{=\kern-.18cm /\kern.08cm}
\def\Ad{{\rm Ad}}
\def\RSU{\R SU}
\def\ad{{\rm ad}}
\def\dva{\bar\partial_A}
\def\da{\partial_A}
\def\p{\partial\bar\partial}
\def\sp{\Sigma^{+}}
\def\sm{\Sigma^{-}}
\def\spm{\Sigma^{\pm}}
\def\smp{\Sigma^{\mp}}
\def\oo{{\scriptstyle{\cal O}}}
\def\ooo{{\scriptscriptstyle{\cal O}}}
\def\sw{Seiberg-Witten }
\def\pa{\partial_A\bar\partial_A}
\def\Dr{{\raisebox{0.15ex}{$\not$}}{\hskip -1pt {D}}}
\def\gr{{\scriptscriptstyle|}\hskip -4pt{\g}}
\def\subsetint{{\  {\subset}\hskip -2.45mm{\raisebox{.28ex}
{$\scriptscriptstyle\subset$}}\ }}
\def\nr{\parallel}
\def\ra{\rightarrow}

\title{The universal  Kobayashi-Hitchin correspondence on Hermitian manifolds}
\author{ M. L\"ubke,   A. Teleman}
\maketitle

\tableofcontents

\begin{abstract} We prove a very general  Kobayashi-Hitchin correspondence
on arbitrary compact Hermitian manifolds, and we discuss differential
geometric properties of the corresponding moduli spaces.
This correspondence  refers to moduli spaces
of ``universal holomorphic  oriented pairs". Most of the classical moduli
problems in complex geometry (e. g. holomorphic bundles with reductive
structure groups, holomorphic pairs, holomorphic Higgs pairs,  Witten
triples, arbitrary quiver moduli problems)  are special cases of this
universal classification problem.  Our Kobayashi-Hitchin correspondence
relates the complex geometric concept ``polystable oriented holomorphic
pair" to the existence of a
 reduction solving a generalized Hermitian-Einstein equation.  The proof is
based on the Uhlenbeck-Yau continuity method. Using idea from  Donaldson
theory, we further introduce and investigate   canonical Hermitian  metrics
on such moduli spaces. We discuss in detail remarkable classes of moduli
spaces in the non-K\"ahlerian framework: Oriented holomorphic
structures, Quot-spaces,  oriented holomorphic pairs and oriented vortices,
non-abelian Seiberg-Witten  monopoles.
\end{abstract}


\setcounter{section}{-1}

\section{Introduction}

In order to understand the aim and the motivation of this article let us consider the
following  {\it classical  complex geometric moduli problems}:
\paragraph{ Holomorphic structures with fixed determinant.}

Let $X$ be a compact complex $n$-dimensional manifold and $E$ a differentiable rank
$r$ vector bundle on
$X$. We fix a holomorphic structure ${\cal L}$ on the determinant line  bundle $\det E$
of $X$.  The problem is to classify all holomorphic structures ${\cal E}$ on $E$ which
induce the fixed holomorphic structure ${\cal L}$ on $\det E$  modulo the
group $\Gamma(X,SL(E))$ of automorphisms   of determinant 1.

In order to get a Hausdorff moduli spaces with good properties one considers only
{\it stable} (or more general {\it semistable}) holomorphic structures.  The stability condition depends on the choice of a Gauduchon metric on $X$,
or more generally a Hermitian metric $g$ (see \cite{LT}) which plays the same role as
the choice of   a polarization of the base manifold in algebraic geometry. ${\cal E}$ is
called $g$-stable if
\begin{equation}\label{BundStab}
\mu_g({\cal F})<\mu_g({\cal E})
\end{equation}
for any nontrivial subsheaf ${\cal F}\subset {\cal E}$ with torsion free quotient. Such  a
subsheaf must be reflexive. Note  that the  slope map $\mu_g$ is \ub{not} a topological
invariant in the general non-K\"ahlerian framework, but it is always a holomorphic invariant.

The {\it classical Kobayashi-Hitchin correspondence} (\cite{Ko1}, \cite{Do1}, \cite{Do2},
\cite{UY1},
\cite {UY2}, \cite{Bu},
\cite{LY},  \cite{LT})  states that a  holomorphic structure is stable  if and only if it
is simple   (i. e.   it admits  no non-trivial trace free
infinitesimal automorphisms) and admits a   Hermitian-Einstein metric, i. e. a
metric $h$ solving the Hermitian-Einstein equation:
$$i\Lambda_g F_h=\frac{2\pi}{(n-1)! Vol_g(X)}\mu_g({\cal E})\id_E \ .
$$
General solutions of the Hermitian-Einstein equation
correspond to {\it polystable} holomorphic structures, i. e. to bundles which are direct
sums of stable bundles of the same slope.

{\it From  a historical point of view}, around 1980 the
correspondence was independently suggested by Hitchin
 and Kobayashi (according to Donaldson \cite{Do1} who first
linked the two names together\footnote{We thank the referee for
pointing out this fact and reference \cite{KO}.}). First
indications that a  connection between these two concepts might
exist had been e.g. Kobayashi's proof that a Hermitian-Einstein
bundle is stable in the sense of Bogomolov \cite{Ko4}, and the
first author's proof of the Chern class inequality for
Hermitian-Einstein bundles \cite{Lu1}.  The ``simple" implication
in this case, i.e. the polystability of Hermitian-Einstein
bundles, was first written down  in Japanese by Kobayashi in
\cite{Ko3} and announced in \cite{Ko1}; the first author   then
found a simpler proof \cite{Lu2}. The ``difficult" implication,
i.e. the existence of a Hermitian-Einstein connection in a stable
bundle, was subsequently proved by Donaldson for Riemann surfaces
\cite{Do3}, algebraic surfaces \cite{Do1} and manifolds
\cite{Do2}, by Uhlenbeck and Yau for K\"ahler manifolds
\cite{UY1}, \cite{UY2}, by Buchdahl for arbitrary Hermitian
surfaces \cite{Bu}, and finally by Li and Yau for arbitrary
Hermitian manifolds \cite{LY}. In the last two articles, as in the
present work, an important technical tool is the use of a  {\it
Gauduchon} metric in the base manifold. An extensive discussion of
the classical correspondence, including complete proofs and
several applications, can be found in [LT].

Let us fix a Hermitian metric $h$ on $E$. The Kobayashi-Hitchin   correspondence yields
an isomorphism of moduli spaces between the moduli space of stable holomorphic
structures with fixed determinant and the moduli space of irreducible integrable
Hermitian-Einstein  connections $A$ with fixed determinant on $(E,h)$. This isomorphism plays a
fundamental role in Donaldson theory, since -- on complex surfaces -- the latter moduli space can
be further identified with a moduli space of  projectively ASD connections
(see
\cite{DK},
\cite{LT} and section \ref{IsoBdls} in this article for details).

\paragraph{Higgs pairs.}

We fix  a differentiable vector   bundle $E$ on $X$, a holomorphic structure
${\cal L}$ on $\det E$ and a holomorphic vector bundle ${\cal F}_0$ on $X$.
In the classical theory of Higgs fields one takes ${\cal F}_0=\Omega^1_X$
(see \cite{Hi}, \cite{Si}).

The problem here is to classify modulo $\Gamma(X,SL(E))$ all pairs $({\cal E},\Phi)$
where
${\cal E}$ is again a holomorphic structure on $E$ which induces ${\cal L}$ on $\det E$,
and
$\Phi\in H^0(\End({\cal E})\otimes{\cal F}_0)$ is a holomorphic ${\cal F}_0$-twisted
endomorphism.

For such objects one also has a stability condition: this condition asks that (\ref{BundStab})
holds for all nontrivial $\Phi$-invariant  subsheaves ${\cal F}$.     This complex geometric
condition can be again characterized in a differential geometric way (see \cite{Hi},
\cite{Si}): a Higgs holomorphic pair
$({\cal E},\Phi)$ is stable if and only if is simple   and
${\cal E}$   admits a metric
$h$ satisfying the Higgs equation
$$i\Lambda F_h +\frac{1}{2}[\Phi,\Phi^{*_h}]=\frac{2\pi}{(n-1)!
Vol_g(X)}\mu_g({\cal E})\id_E
$$
Again, this result    yields an isomorphism of moduli spaces: the complex geometric
moduli space of stable holomorphic Higgs pairs and the differentiable geometric moduli
space of irreducible integrable solutions $(A,\Phi)$ of the Higgs equation. Arbitrary
solutions of this equation correspond to polystable holomorphic Higgs pairs (see
\cite{Hi}, \cite{Si}).

\paragraph{Holomorphic pairs.}  (see \cite{Bra}, \cite{Th}, \cite{HL})  Let
$E$ be a rank $r$ bundle on $X$. One wants to classify -- modulo
$\Gamma(X,GL(E))$ -- the pairs
$({\cal E},\varphi)$, where ${\cal E}$ is a holomorphic structure on $E$ and $\varphi$ is a
holomorphic section in
${\cal E}$.

The stability condition for this problem depends on a real parameter $\tau$. A pair
$({\cal E},\varphi)$ is $\tau$-stable if
$$\max(\mu_g({\cal E}),\sup\limits_{{\cal F}\in{\cal R}(E)}\mu_g({\cal
F}))<\tau<\inf\limits_{\matrix{\scriptstyle{\cal F}\in{\cal
R}(E)\cr \scriptstyle
\varphi\in H^0({\cal F})}}\mu_g({\cal E}/{\cal F})\ ,
$$
where ${\cal R}({\cal E})$ denotes the set of subsheaves of ${\cal E}$ with torsion
free quotients.

There is again a differential geometric characterization of the stable pairs. $({\cal
E},\varphi)$ is $\tau$-stable if and only if it is simple (i. e. it admits no non-trivial
infinitesimal automorphisms)  and ${\cal E}$ admits a  metric
$h$ satisfying the {\it vortex equation}
$$i\Lambda F_h+\frac{1}{2}\varphi\otimes\varphi^{*_h}=\frac{2\pi}{(n-1)!
Vol_g(X)}\ \tau\id_E\ .
$$

\paragraph{Witten triples.} (see \cite{W}, \cite{Du}) Consider a
line  bundle  $L$  on  a complex surface
$X$.   We want to  classify modulo
$\Gamma(X,GL(L))={\cal C}^{\infty}(X,\C^*)$ the triples $({\cal L},\varphi,\alpha)$
consisting of a holomorphic structure ${\cal L}$ on $L$, a holomorphic section
$\varphi\in H^0({\cal L})$ and a holomorphic morphism    $\alpha:{\cal L}\ra {\cal
K}_X$.

The stability condition depends again on a real parameter $\tau$. A triple $({\cal
L},\varphi,\alpha)$ is $\tau$-stable if   either $\deg({\cal L})<\tau$ and $\varphi\ne
0$, or $\deg({\cal L})>\tau$ and $\alpha\ne 0$.

A    triple $({\cal L},\varphi,\alpha)$ is stable if and only if it is
simple (i. e. $(\varphi,\alpha)\ne 0$) and   ${\cal L}$ admits a metric $h$ satisfying
the {\it mixed vortex equation}
$$i\Lambda_g F_h+\frac{1}{2}(|\varphi|^2-|\alpha|^2)=\frac{2\pi}{Vol_g(X)}\ \tau\ .
$$
This result yields  again an isomorphism of moduli spaces.
\\ \\
{\bf Remark:} In the vortex equation and the mixed vortex equation one can replace
the constant $\frac{2\pi}{(n-1)!
Vol_g(X)}\ \tau$ on the right by a real function $t$. The corresponding stability
condition depends only on $\int_X t\ vol_g$ when $g$ is Gauduchon.  The case when
$t=\frac{s_g}{2}$ (where $s_g$ stands for the scalar curvature of $g$) plays a crucial
role in the Seiberg-Witten theory on   complex surfaces (see \cite{W}, \cite{Bi}, \cite{Du}
\cite{OT1},
\cite{OT2}, \cite{Te2}).

\paragraph{Oriented holomorphic pairs.} We fix again a differentiable vector bundle $E$
of rank $r$ on $X$ and a holomorphic structure ${\cal L}$ on $\det E$. This time we
want to classify modulo $\Gamma(X, SL(E))$ the pairs $({\cal E},\varphi)$, where
${\cal E}$ is a holomorphic structure on $E$ which induces ${\cal L}$ on $\det E$ and
$\varphi\in H^0({\cal E})$ is a holomorphic section.

The stability condition in this case is (see \cite{OST}, \cite{OT2},
\cite{OT8},
\cite{Te2})
$$\max(\mu_g({\cal E}),\sup\limits_{{\cal F}\in{\cal R}(E)}\mu_g({\cal
F}))< \inf\limits_{\matrix{\scriptstyle{\cal F}\in{\cal
R}(E)\cr \scriptstyle
\varphi\in H^0({\cal F})}}\mu_g({\cal E}/{\cal F})\ ,
$$
whereas the corresponding Hermitian-Einstein type equation is
$$i\Lambda F_h+\frac{1}{2}(\varphi\otimes\varphi^{*_h})_0=0
$$

In the case $n=r=2$, one gets a complex geometric interpretation of the moduli spaces
of $PU(2)$-monopoles on complex surfaces (see \cite{Te2} for the
K\"ahlerian case  and section \ref{NASW} in this article for the Gauduchon
case).
\vskip 0.3cm

In order to illustrate the generality  of the   results  in this article, consider now  the
following {\it artificial moduli problem}:

Fix a system of differentiable vector
bundles
$(E_i)_{0\leq i\leq n}$  of ranks $r_i$, a system of holomorphic vector bundles $({\cal
F}_j)_{1\leq j\leq m}$ of ranks $\rho_j$  on
$X$, and a holomorphic structure ${\cal L}_0$ on $\det E_0$.

Classify systems $({\cal E}_i,  \varphi_{ij},\psi_i)$ where  ${\cal E}_i$ is a holomorphic
structure on $E_i$  such that ${\cal E}_0$ induces ${\cal L}_0$ on $\det E_0$,
$\varphi_{ij}\in H^0({\cal E}^\vee_i\otimes {\cal F}_j) $ are holomorphic
homomorphisms, and $\psi_i\in \Gamma(X,\P(\End({\cal E}_i))$ are holomorphic
sections in the projectivizations of the endomorphism bundles of ${\cal E}_i$.

The question is: which is the correct stability (polystability) condition in this case, and
which is the corresponding Hermitian-Einstein equation?

All the ``classical" moduli problems listed above, as well as this artificial moduli
problem are just special case of the following very general complex  geometric
classification  problem:

\paragraph{A universal moduli problem.} We note that in all these examples one has
a system of   {\it fixed  holomorphic structures},  a system of {\it variable  holomorphic
structures} on fixed ${\cal C}^\infty$-bundles, and a system of  {\it variable  sections}
in associated bundles, which are holomorphic with respect to   both   variable
and fixed holomorphic structures.

The best way to formulate a universal generalization of all moduli problems of this type is
the following (see \cite{OT3},   \cite{OT5}):

Consider an exact sequence
$$\{1\}\map G\map \hat G\map G_0\map \{1\}
$$
of complex reductive groups, and choose a $\hat G$-bundle $\hat Q$ on $X$ and a
holomorphic action $\hat\alpha:\hat G\times F\ra F$ on a K\"ahlerian manifold $F$.   Let
$Q_0:=\hat Q/G$ be the associated $G_0$-bundle and $E$  the associated $F$-bundle
$E:=\hat Q\times_{\hat\alpha} F$.  Fix a holomorphic    structure  ${\cal Q}_0$ on
$Q_0$. Our universal classification problem is:
\vskip 0.3cm

{\it Classify the pairs $(\hat {\cal Q}, \varphi)$, where $\hat {\cal Q}$ is a holomorphic
structure on $\hat Q$ inducing ${\cal Q}_0$ on $Q_0$ and $\varphi\in \Gamma(X,E)$
is holomorphic with respect to the holomorphic structure induced by $\hat Q$ on $E$,
modulo the gauge group ${\cal G}:=\Aut_{Q_0}(\hat Q)=\Gamma(X,\hat
Q\times_{\Ad} G)$.}
\vskip 0.3cm

A pair $(\hat {\cal Q},\varphi)$ of this type will be called an {\it oriented pair of type}
$(\hat Q,\hat \alpha,{\cal Q}_0)$.

This class of moduli problems was first considered in \cite{OT3} in the case when
$X$ is a complex curve; in this case the corresponding moduli spaces were  used to
introduce the  {\it twisted gauge theoretical Gromov-Witten invariants}. The particular case $G_0=\{1\}$ (the case of
{\it non-oriented pairs}) was previously studied by Mundet i Riera  in \cite{Mu1} on K\"ahler
manifolds.
\\ \\
{\bf Example:} For instance, to get our ``artificial moduli problem" above as
a special case of our  universal problem, take
$$ G:=SL(r_0,\C)\times\prod_{i=1}^n GL(r_i,\C)\ ,\ \hat
G:=GL(r_0,\C)\times\prod_{i=1}^n GL(r_i,\C)\times\prod_{j=1}^m GL(\rho_j,\C)\ ,
$$
$$G_0=\C^*\times\prod_{j=1}^m GL(\rho_j,\C)\ ,\ F:=\prod_{\matrix{\scriptstyle
0\leq i\leq n\cr\scriptstyle 1\leq j\leq
m}}\Hom(\C^{r_i},\C^{\rho_j})\times\prod_{i=0}^n\P(\End(\C^{r_i}))\ .
$$
The action $\hat \alpha$ is the natural action of $\hat G$ on $F$, whereas the fixed
holomorphic structure
${\cal Q}_0$ is defined by the system of fixed holomorphic structures $({\cal L},({\cal
F}_j)_{1\leq j\leq m})$.

Similarly  one can see easily that all classical moduli problems  in our
list are also special cases of this universal problem, as well as all ``quiver problems"
considered in \cite{AG}.
\vskip 0.3cm

The  main goals of our article are:\\ \\
I. To give the correct stability and
polystability condition for this universal moduli problem and to
characterize the polystable pairs in terms of the existence of a reduction
of $\hat {\cal Q}$ to a maximal compact subgroup $\hat K\subset\hat G$
which satisfies a certain generalized Hermitian-Einstein equation.  In this
article, we   solve this problem completely (see Theorem \ref{SimImp},
Theorem
\ref{purpose} for the precise statements).

In particular, these results give a Kobayashi-Hitchin correspondence for
holomorphic  principal bundles  with arbitrary reductive groups on arbitrary compact
Hermitian manifolds.
\\ \\
II. To study, in the non-K\"ahlerian framework, metric properties of this
large class of moduli spaces: we show that the (smooth part of the) moduli
space  associated with such a problem  comes always with a natural Hermitian
metric, and   that this Hermitian metric  has interesting  properties (for
instance is K\"ahler as soon as the base manifold is semi-K\"ahler). The
construction of this {\it canonical Hermitian metric} uses ideas from
Donaldson theory.
\\ \\
II. To apply our general results to remarkable classes of moduli
spaces: moduli spaces of oriented connections, Douady Quot-spaces, moduli
spaces of non-abelian monopoles on Gauduchon surfaces.

\paragraph{Previous results on universal moduli problems.}

In this paragraph we acknowledge previous contributions on the subject, from which we benefited a lot.    On the other hand we
explain        carefully the progress made in the present article
compared with   previous work.

The idea to a find a ``universal" generalization of the correspondence for {\it all possible
coupled} moduli problems is of course very natural. Important progress in this direction was
first achieved by Banfield
\cite{Ba}  who  proved a Kobayashi-Hitchin correspondence for {\it non-oriented} pairs whose
second       component is a section in vector bundle associated with a   {\it linear representation}
$\alpha$.  Banfield's stability condition is a natural extension to coupled problems of the stability
condition for  holomorphic bundles with arbitrary reductive structure group as developed by
Ramanathan and Ramanan-Ramanathan
\cite{Ra}, \cite{RR}.   In his very interesting Ph. D. thesis (see
\cite{Mu1}),   Mundet i Riera   succeeded to   generalize Banfield's result to the case of
an arbitrary (not necessary linear)    action  on a K\"ahlerian  fibre.    In the present
article we use essentially his important contribution on the subject.

We also mention that  a very general Kobayashi-Hitchin correspondence on
K\"ahler manifolds is stated in
\cite{BGM}, where the same type of holomorphic pairs is considered,
but the symmetry group of the  problem is a more general  ``Lie subgroup"    ${\cal
G}\subset
\Aut(\hat Q)$; but there is a price to pay for this generality: the inequality in the
corresponding stability condition depend on a vector varying in an infinite dimensional Lie
algebra, so this is not a practically checkable condition.

All these contributions concern   the particular case when the
base manifold is K\"ahlerian and they all use similar methods.

Compared with   the existing  literature, the results and the   formalism developed
in the present article has several important advantages :
\begin{enumerate}
\item  Our methods applies  to the
general non-K\"ahlerian  framework.

We mention that in the non-K\"ahlerian framework several important difficulties arise:
\begin{enumerate}
\item  The left term of the
Hermitian-Einstein equation  can no longer be regarded as a moment map,
\item  The
K\"ahler identities do not hold,
\item The degree of holomorphic bundles is no longer a topological invariant.
\end{enumerate}
For these reasons, the methods used in \cite{Mu1}, \cite{BGM} do not appear to apply to the
non-K\"ahlerian framework.
\item We treat not only the correspondence between stable objects and irreducible solutions of
the  (generalized) Hermitian-Einstein equations (see \cite{Mu1}, \cite{BGM}) but the general
correspondence between  {\it polystable} objects  and {\it arbitrary}    solutions of these equations.

This is also an important point, because even in the classical Donaldson theory, the possible
presence of the ``reductions" in the moduli spaces must be always carefully taken into account.
\item  Our stability condition has  {\it holomorphic character}; more precisely, the terms
in our stability inequalities   are intrinsic {\it holomorphic invariants}
associated with the  oriented holomorphic pair $(\hat {\cal
Q},\varphi)$  and  certain {\it meromorphic reductions} of $\hat {\cal
Q}$.  Only the {\it given} parameters (the Hermitian base manifold and
the Hamiltonian action on the fibre) are used in the construction of
these invariants. {\it Holomorphically isomorphic oriented pairs will have
the same set of invariants}.

These invariants {\it   do not require and do not depend on a
background
$\hat K$ -  reduction} of
$\hat {\cal Q}$  (contrary
to
\cite{Mu1} and
\cite{BGM}). This is    important (even in the K\"ahler and algebraic
case!), because    a
$\hat K$-reduction is a differential geometric object; its
construction   requires  a partition of unity (as does the construction
of a Hermitian structure of a vector bundle), so a stability condition
which depends on the choice of a $\hat K$-reduction is not practically
checkable.

The complex geometric character of our stability condition is a consequence of the
complex geometric equivariance properties -- obtained in
\cite{Te3}    -- of the maximal weight function associated with a
finitely dimensional {\it energy complete Hamiltonian triple}  (see
Proposition \ref{FundProp}).  Without the energy completeness
condition, there seems to be no way to prove these complex geometric
equivariance properties.
\end{enumerate}

The main feature of the classical stability condition for bundles   is its
complex geometric nature: even on general Hermitian  manifolds the  slope map on the
set of non-trivial subsheaves of a given holomorphic bundle   has a pure complex
geometric character, because it does not depend on the choice of a Hermitian metric on
this bundle (see
\cite{Ga},
\cite{LT}).

Therefore our universal Kobayashi-Hitchin correspondence respects this fundamental feature
of the classical Kobayashi-Hitchin correspondence, i. e. it deals only with   complex
geometric stability   and polystability conditions.

In fact this difficulty already arises   in the finite dimensional framework (see
section \ref{FinDim}): the generalized maximal  weight function $\lambda$ associated with a
 $K$-moment map $\mu$ for a general holomorphic action  $\alpha:G\times
F\ra F$ on a K\"ahler manifold    is not   equivariant  with respect to the whole
complex reductive group
$G$, but only with respect to the fixed maximal compact subgroup $K$. Therefore, the
inequality
$\lambda(x)>0$ involved in the  ``analytic stability condition" is a priori {\it not} a
complex geometric condition. Moreover, in the general case, it is not possible to give
numerical characterizations of   semistable and polystable points.

In the finite dimensional framework these difficulties have been explained and overcome
in
\cite{Te3}, where it was shown that one has to impose a certain completeness
condition on the pair $(\alpha,\mu)$ in order to have a maximal weight function
$\lambda$ with good complex geometric equivariance properties.  For this class of
actions (which includes all holomorphic Hamiltonian actions on compact K\"ahler
manifolds and all linear actions) one can  give    complex geometric  numerical
characterizations of stable, semistable and polystable orbits.

Using this formalism, we were able to extend the    Kobayashi-Hitchin
correspondence  for ``universal" oriented holomorphic pairs to the case when the base
manifold $X$ is an arbitrary compact  Hermitian manifold. Moreover,  we give a
complex geometric numerical characterization not only of the stable, but also of the
polystable pairs (i. e. of all pairs which  admit a $\hat K$-reduction satisfying the
generalized Hermitian-Einstein  equation). This gives a complete complex geometric
interpretation of the whole moduli space of   solutions of the generalized
Hermitian-Einstein equation.

Extending the Kobayashi-Hitchin correspondence to the non-K\"ahlerian
framework  is important for many reasons: for instance this generalization
furnishes complex geometric
descriptions of the moduli spaces of instantons with arbitrary compact
symmetry groups  and of
the moduli spaces of Seiberg-Witten monopoles (abelian and non-abelian)
on all complex
surfaces.  Furthermore,    important applications of the classical
Kobayashi-Hitchin correspondence to the
classification of non-K\"ahlerian surfaces  can be found in  \cite{Te1},
\cite{Te4}
\cite{LT}, \cite{LYZ}.

Our proof is based on the continuity method developed by   Uhlenbeck and Yau for the
classical Kobayashi-Hitchin correspondence (relating polystable bundles to
Hermitian-Einstein connections). This seems to be the appropriate method in     the
non-K\"ahlerian framework, because in this case the left hand term in the generalized
Hermitian-Einstein equation {\it is no longer a formal moment map}.  On the other
 hand, the original proof of Uhlenbeck-Yau for the classical Kobayashi-Hitchin
 correspondence does not extend mechanically to our class of coupled problems with
 general symmetry  groups, so the classical strategy  developed by
 these authors must be completed with new
 techniques.

Usually  moduli problems for bundles with reductive structure group $G$ are
treated using a faithful representation $G\ra GL(r,\C)$ and reducing the
 problem to the case $G=GL(r,\C)$, i. e. to the more familiar case of vector bundles.   We
will {\it not} follow  this tradition, because we realized that this method does not
really simplify the problem, and  dealing  directly with the general case is more natural
and interesting.   Consequently, we  need an analogue of the well-known theorem of
Uhlenbeck-Yau on weakly holomorphic subbundles,  for general reductive structure
groups $G$.  In other words, we show  how one can construct a meromorphic
reduction of a holomorphic principal $G$-bundle ${\cal Q}$ to a parabolic subgroup
$L\subset G$ using an  $L^2_1$-section $s\in \Gamma(\ad({\cal Q}))$
  which satisfies an algebraic property and a certain weak holomorphy condition.
This correspondence, which follows easily from the results in
\cite{UY1},
\cite{UY2}, is   presented in the Appendix.

\paragraph{The structure of the article.}

We begin with the ``finite dimensional"
Koba\-ya\-shi-Hitchin correspondence. This result gives a numerical characterization  of the
polystable orbits with respect to a  holomorphic Hamiltonian  action
$(\alpha,\mu)$ satisfying our technical completeness condition. This result can be
regarded as a ({\it non-algebraic}!) complex geometric   generalization of  the
well-known Hilbert criterion in classical GIT. A different proof   can be found in
\cite{Mu1}      (for the stable case)     and \cite{Te3} for the general semistable
case. Here we give a proof which is based on the same continuity method which will be
used later in the infinite dimensional gauge     theoretical framework.

The main  ideas of the results presented here are present already in  the  classical  GIT contributions
(\cite{Ki}, \cite{MFK}, \cite{KN}); the goal of this chapter is just to treat  carefully the delicate
{\it non-algebraic non-compact} case and to point out the specific difficulties (see
\cite{Te3}, \cite{BT}).

In the second chapter we state our ``universal" moduli problem and we introduce the
stability, semistability and   polystability conditions for   universal oriented
pairs. Much care is dedicated   here to the notion of {\it degree} of a meromorphic
reduction to a parabolic subgroup  with respect to an $\ad$-invariant linear form; this
concept   is very delicate in the non-K\"ahler framework.

In the third chapter we introduce the   Hermitian-Einstein equation  for universal oriented
pairs and we prove the ``simple" implication of the Kobayashi-Hitchin correspondence:
any pair which admits a Hermitian-Einstein reduction is polystable.

The fourth chapter deals with the ``difficult" implication: any polystable pair admits an
admissible Hermitian-Einstein reduction.  The proof of this, based on the Uhlenbeck-Yau continuity method, uses involved analytical techniques.

In the fifth chapter we give   applications of our results:  The starting
point of the chapter is to notice that, as in Donaldson theory, the obtained
correspondence between solutions of the generalized  Hermitian-Einstein
equation and polystable oriented  pairs yields always an isomorphism of
moduli spaces (the Kobayashi-Hitchin isomorphism). In the first section of
the chapter we treat in detail the case  of the correspondence between
oriented holomorphic structures and oriented connections.  We give a simple
argument to explain why moduli spaces of {\it oriented} connections are
important: the ``cobordism argument" used for relating the    Donaldson
invariants to the Seiberg-Witten invariants generalizes to non-simply
connected  4-manifolds, but one should use   Donaldson invariants
constructed with moduli spaces of oriented, projectively ASD
$U(2)$-connections, instead of moduli spaces of
$PU(2)$-instantons.    In the second section of the chapter we turn
to coupled problems. For both classes of moduli problems (non-coupled
and coupled)  we use the Kobayashi-Hitchin isomorphism to construct a
canonical Hermitian metric on (the smooth part of) the moduli space. We
prove in full generality that this metric   is always  K\"ahler when the
base manifold is semi-K\"ahler, and that its K\"ahler form
$\Omega$ satisfies the equation
$\partial\bar\partial\Omega=0$, when the base manifold is Gauduchon. The
main idea of the proof comes from Donaldson theory: for the construction of
$\Omega$ we use a method similar to the one used in \cite{DK} to construct
de Rham representatives of the
$\mu$-classes.   This  generalizes  the results concerning the
metric properties of the classical moduli spaces of holomorphic vector
bundles, obtained with different methods in \cite{LT}.   As an application
we endow Douady Quot-spaces with canonical Hermitian metrics, by identifying
them with suitable moduli spaces of stable oriented holomorphic pairs. The
more delicate non-K\"ahlerian case is treated carefully.  The third section
of the chapter is dedicated to non-abelian Seiberg-Witten theory for
Gauduchon surfaces. We give a geometric interpretation of the moduli spaces
of non-abelian monopoles in this general framework in terms of oriented rank
2 holomorphic pairs.

The Appendix  contains   technical results (i.e. identities for Chern connections, local
diagonalization results for arbitrary compact Lie groups,  meromorphic reductions to
parabolic subgroups)  which  are all of independent interest.
\paragraph{Notations and conventions.}
We tried to respect the standard conventions and notations used in the literature. However, for
mathematical reasons, we  also had to introduce some new ones:

For instance, in \cite{Mu1} the notations $E$, $E_G$ are used for a principal $K$ bundle and its
$G=K^\C$-extension $E\times_K G$.  In our article,  $\hat P$ stands for a $\hat K$ principal bundle
(standard differential geometric notation) and
$\hat Q$ for a $\hat G$-principal bundle, where $\hat K$ is a maximal compact subgroup in $\hat
G$.  It is very important that we do {\it not} fix  an identification     $\hat Q=\hat
P\times_{\hat K}\hat G$ when stability is introduced, because fixing such an
identification is equivalent to choosing a background $\hat K$-reduction of $\hat Q$,
and, as explained above, {\it stability should be a holomorphic
condition referring to holomorphic objects}.

For the same reason, the first component of our holomorphic pairs  is a bundle-holomorphic
structure   on $\hat Q$ (denoted by $\hat J$ or $\hat {\cal Q}$) and not an integrable
$\hat K$-connection
$\hat A$  (compare with \cite{Mu1}, \cite{BGM}). Passing from holomorphic structures
to connections is possible only after fixing a
$\hat K$-reduction of $\hat Q$.

Third, the real and complex  gauge groups   are denoted by    {\curly  G}$_K$ and  {\curly
G}$_G$ in   \cite{Mu1}. To simplify the notation, we preferred the symbols
${\cal K}$ and
${\cal G}$.

Note also that  our stability condition is slightly simpler than the stability condition in
\cite{Mu1}. Indeed,
 the ``initial parameter" in \cite{Mu1} is a pair $(L,\xi)$ consisting of
a parabolic subgroup of
$G$ and a dominant character of $L$. On the other hand, our initial parameter is just a
Hermitian type vector $\xi$ in the Lie algebra $\g$.  The point is that
we only need a parabolic subgroup of $G$ which is {\it canonically
associated with} $\xi$.  Since the functoriality of this parabolic  group
with respect to both
$G$ and $\xi$ plays an important role, we used the notation $G(\xi)$ for this parabolic
subgroup.

Moreover, the central constant  $c$ in the ``$c$-stability condition"
formulated in \cite{Mu1} can be absorbed in the moment map $\mu$, so without any
loss of generality, one can omit this parameter.
\\ \\
{\bf Acknowledgment:} This work grew out as part of an ample research project on
the universal Kobayashi-Hitchin correspondence and its applications,  which
started in Z\"urich in collaboration with Ch. Okonek.

\section{  The finite dimensional
Kobayashi-Hitchin
  correspondence }\label{FinDim}

In this section we give a brief presentation of the stability theory in the  finitely dimensional
K\"ahlerian {\it non-algebraic}  framework and   explain the analogue
of the Hilbert criterion in this framework. The Hilbert criterion we prove here  gives a
numerical characterization of the polystable orbits, i. e. of the orbits which intersect
the vanishing locus of the moment map.   Our proof --  based on  the continuity
method --  is a very good preparation for understanding the proof of our
universal Kobayashi-Hitchin correspondence on Hermitian manifolds in section
\ref{DiffImpl}, which    will follow
the same strategy but will require involved analytical techniques.

\subsection{Analytic Stability, Symplectic stability}\label{AnStabSympStab}

For complete proofs and more details about the notions and the results
introduced in this section we refer to \cite{Te3}.  The analogous results in the algebraic geometric
framework are well known (\cite{Ki}, \cite{MFK}, \cite{KN}), but the non-algebraic  framework
raises specific difficulties.

Let $G$ be a complex reductive group. We denote by $H(G)$ the  subset
$H(G)\subset\g$ defined by
$$H(G):=\{\xi\in\g|\ \overline{\exp(i\R\xi)}\ \hbox{is
compact}\}=\union\limits_{\matrix{K\subset G\cr\hbox{maximal compact}}} i\kg\ .
$$

$H(G)$  is a locally closed cone in $\g$.  The elements in $H(G)$ will be called
vectors of Hermitian type (see \cite{Te3}).   There is an obvious {\it non-surjective} injective
map
$\Hom(\C^*,G)\to H(G)$ given by
$$\theta\mapsto \left.\frac{d}{dt}\right|_{t=0}(\R\ni t\mapsto\theta(1+t))=d_1\theta(1) \ .
$$
Therefore, compared to classical GIT,  in non-algebraic complex geometry, one has to consider more
general (non-algebraic!) one-parameter subgroups, which do not exponentiate     to morphisms
$\C^*\to G$.

To any $\xi\in H(G)$ one can
associate a parabolic subgroup $G(\xi)$ of $G$ by
$$G(\xi):=\{g\in G|\ \lim_{t\ra\infty}\exp(t\xi)g\exp(-t\xi)\ \hbox{exists in}\ G\}
$$

The subgroup $G(\xi)$ fits   in a short exact sequence
\begin{equation}\label{ExSeq}
\{1\}\map U(\xi)\map G(\xi)\map Z(\xi)\map \{1\}\ ,
\end{equation}
where
$$U(\xi):=\{g\in G|\ \lim_{t\ra\infty}\exp(t\xi)g\exp(-t\xi)=e\}
$$
is the unipotent radical of $G(\xi)$ and
$$Z(\xi):=\{g\in G|\ \ \ad_g(\xi)=\xi\}\ ,
$$
is the reductive centralizer of $\xi$.

$G(\xi)$ is the direct product of its subgroups $U(\xi)$,
$Z(\xi)$. Note that that the normal subgroup $U(\xi)$  is canonically associated
with the parabolic subgroup $G(\xi)$, whereas the  subgroup $Z(\xi)\subset
G(\xi)$ depends on the choice of the vector $\xi$.

The Lie algebras of these groups are
$$\ug(\xi):=\bigoplus_{\lambda<0}{\rm Eig}([\xi,\cdot],\lambda)\ ,\
 \zg(\xi)=z_\g(\xi)=\ker  ([\xi,\cdot])\ ,\ \g(\xi)=\bigoplus_{\lambda\leq 0}{\rm
Eig}([\xi,\cdot],\lambda)\ .
$$

Following \cite{Te3}  we introduce in $H(G)$ the following important equivalence
relation:
\begin{dt}
We say that two vectors
$\xi$, $\zeta\in H(G)$ are equivalent  ($\xi\sim \zeta$) if one of the following
equivalent conditions is satisfied:
\begin{enumerate}
\item $\zeta\in \g(\xi)$ and there exists $g\in G(\xi)$ such that
$\ad_g(\xi)=\zeta$.
\item   $\zeta-\xi\in \ug(\xi)$.
\item $\xi\in \g(\zeta)$ and there exists $g\in G(\zeta)$ such that
$\ad_g(\zeta)=\xi$.
\item  $\xi-\zeta\in \ug(\zeta)$.
\end{enumerate}
\end{dt}

\begin{pr}\label{RepSyst} Let $K\subset G$ be a maximal compact subgroup of $G$.
Then the composition
$$i\kg\hookrightarrow H(G)\map \qmod{H(G)}{\sim}$$
is a homeomorphism.
\end{pr}

In particular, $i\kg$ is a complete system of representatives for the equivalence
relation $\sim$.
\begin{co} If $G$ is a closed reductive subgroup of the reductive group $\hat G$, then the
inclusion
$H(G)\ra H(\hat G)$ induces an injection
$$\qmod{H(G)}{\sim_G}\ra \qmod{H(\hat G)}{\sim_{\hat G}}
$$
which is a homeomorphism on its image.
\end{co}
{\bf Example:} Consider the case $G=GL(r,V)$, where $V$ is a finite dimensional
complex vector space. Then $H(G)\subset gl(V)$ is the subset of endomorphisms which are
diagonalizable and have real eigenvalues. Every such endomorphism $\xi$ defines a
filtration
$$V(\xi)=(V_\lambda)_{\lambda\in{\rm Spec}(\xi)}\ ,\
V_\lambda:=\bigoplus_{\nu\leq\lambda} {\rm Eig}(\xi,\nu)
$$
of $V$.

The parabolic subgroup $G(\xi)$ consists of those elements $g\in GL(V)$ which
leave the filtration $V(\xi)$ invariant.  Two elements $\xi$, $\zeta\in H(G)$ are
equivalent if and only if they have the same eigenvalues and define  the same
filtration.
\\
\\

Let $\alpha:G\times F\ra F$ be  a holomorphic action of a complex reductive group
on a complex manifold $F$.  A {\it Hamiltonian triple} for $\alpha$ is a triple
$\tau=(K,g,
\mu)$ consisting of a maximal compact subgroup $K$ of $G$, a $K$-invariant K\"ahler
metric
$g$ on
$F$,  and a moment map $\mu$ for the induced $K$-action on $(F,\omega_g)$.

The data of a Hamiltonian triple for the holomorphic action $\alpha$ allows one to
construct a {\it K\"ahlerian quotient} of $F$ by $G$ (see \cite{HH}).   As in the
algebraic geometric GIT (where the quotient depends on the  choice of a linearization
of the action), a K\"ahler quotient depends in general
essentially on the choice of the Hamiltonian triple $\tau$.

\begin{dt} Let $\tau=(K,g, \mu)$ be  Hamiltonian triple for the action $\alpha$. A
point $x\in F$ will be called:
\begin{enumerate}
\item (symplectically) $\tau$-semistable if $\overline{Gx}\cap\mu^{-1}(0)\ne
\emptyset$.
\item (symplectically) $\tau$-stable if $\g_x=\{0\}$ and $Gx \cap\mu^{-1}(0)\ne
\emptyset$.
\item  (symplectically) $\tau$-polystable if   $Gx
\cap\mu^{-1}(0)\ne
\emptyset$.
\end{enumerate}
\end{dt}

We will denote by $F^{\rm sst}_\tau$, $F^{\rm st}_\tau$, $F^{\rm
pst}_\tau$ the  subsets of symplectically $\tau$-(semi, poly) stable
points in $F$.

By the results of \cite{HH} one has
\begin{thry}
\begin{enumerate}
\item $F^{\rm sst}_\tau$ is open in $F$.
\item  There exists a good quotient  $\pi:F^{\rm sst}_\tau\ra Q_\tau$,
where $Q_\tau$ is a Hausdorff complex space.
\item The closure of every semistable orbit contains a unique polystable orbit.
\item The inclusions $\mu^{-1}(0)\hookrightarrow  F^{\rm
pst}_\tau\hookrightarrow F^{\rm sst}_\tau$ induce homeomorphisms
$$\qmod{\mu^{-1}(0)}K\textmap{\simeq} \qmod{F^{\rm
pst}_\tau}{G}\textmap{\simeq}Q_\tau
$$
In particular two semistable orbits have the same image in $Q_\tau$ if and only if
the  intersection of their closures contains a polystable orbit.
\end{enumerate}
\end{thry}

Therefore the K\"ahler quotient $Q_\tau$ can  identified  with the symplectic
quotient $\qmod{\mu^{-1}(0)}{K}$, so this results states that, in the K\"ahler
case, this symplectic quotient (which in general can be very singular) inherits a
natural complex space structure.

Let again $\tau=(K,g,\mu)$ be a Hamiltonian triple for $\alpha$.
\begin{re}\label{increasing} For any $s\in i\kg$,  the map $t\mapsto  \langle
\mu(e^{ts}x),-is\rangle$ is increasing.
\end{re}
Indeed the derivative of this map at $t$ is just $\|s^\#_{e^{ts}}\|^2$.

We put
$$\lambda^s(x):=\lim_{t\ra\infty}\langle \mu(e^{ts}x),-is\rangle\ .
$$
Note that
\begin{equation}\label{energy}
\lambda^s(x)=\langle \mu(e^{ts}x),-is\rangle+
E_g(c^s_x)
\end{equation}
 where
$c^s_x:[0,\infty)\ra F$   is the parameterized curve  $t\mapsto e^{ts} x$ and
$E_g$ stands for the energy with respect to the metric $g$.

Using the standard equivariance property of the moment map one gets the
$K$-equivariance property of the map $\lambda$
\begin{equation}\label{K-eq}
\lambda^s(x)=\lambda^{\ad_{k }(s)}(k x)\ \forall s\in
i\kg\ ,\ k\in K
\end{equation}

The following stability criterion is well-known  (see for instance \cite{Mu1},
\cite{Te3}):
\begin{pr} Let $x\in F$. The following conditions are equivalent
\begin{enumerate}
\item $x$ is symplectically $\tau$-stable.
\item $\lambda^s(x)>0$ for all $s\in i\kg\setminus\{0\}$.
\end{enumerate}
\end{pr}

A point $x$ satisfying the second condition is  usually called  {\it analytically
$\tau$-stable}. This criterion  can be considered as the K\"ahlerian analogue of the
Hilbert criterion in GIT.   However, this criterion is not satisfactory   for our purposes
for the following   reasons:
\begin{enumerate}
\item The map $\lambda$ is not a purely complex geometric object; it depends on
the choice of a Hamiltonian triple for the action $\alpha$, which is a differential
geometric object.
\item  As it stands, the analytic stability condition is not a $G$-invariant condition,
because $\lambda$ is only a $K$-equivariant function, not a $G$-equivariant one.
\item   The criterion does not work for \ub{semi}stable points: for general
holomorphic actions $\alpha$ on non-compact manifolds the condition
$\lambda^s(x)\geq 0$ for all $s\in i\kg$ does \ub{not} imply that $x$ is
symplectically semistable.
\end{enumerate}

The simplest way to overcome these difficulties is  to restrict ourselves to {\it
energy complete} Hamiltonian triples (see \cite{Te3}):
\begin{dt} A Hamiltonian triple $\tau=(K,g,\mu)$ is called energy complete if for
every $(s,x)\in  i\kg\times F$ the following implication holds
$$E_g(c^s_x)<\infty\Rightarrow \lim_{t\ra\infty} c^s_x(t)\
 \hbox{exists}\ in\ F\ .
$$
\end{dt}

We refer to \cite{Te3} for the following
\begin{pr}\label{FundProp} The map $\lambda:i\kg\times F\ra\R\cup\{\infty\}$
associated with an energy complete Hamiltonian triple $\tau=(K,g,\mu)$ has a unique
extension
$$\lambda:H(G)\times F\ra\R\cup\{\infty\}$$
which is $G$-equivariant, i. e.
$$\lambda^s(x)=\lambda^{\ad_{g }(s)}(g x)\ \forall  s\in
H(G)\  \forall x\in F  \forall   g\in G\ .$$
This function has the following properties
\begin{enumerate}
\item homogeneity: $\lambda^{ts}(x)=t\lambda^s(x)$ for any $t\in \R_{\geq
0}$, $s\in H(G)$.
\item parabolic invariance:
$\lambda^s(x)=\lambda^s(h x) \ \hbox{for every}\ x\in F,\ h\in
G(s)\ . $

\item   $\sim$ invariance: $\lambda^s(x)=\lambda^\sigma(x)$ if $s\sim\sigma$.
\item  semicontinuity: if $(x_n,s_n)_n\ra (x,s)$ then $\lambda^s(x)\leq
\liminf\limits_{n\ra\infty}\lambda^{s_n}(x_n)$.

\end{enumerate}

\end{pr}

   The $G$-equivariant extension   $\lambda:H(G)\times F\ra\R\cup\{\infty\}$
will be called (for historical reasons, see   section  \ref{ProjAct}) the
\ub{maximal} \ub{weight} \ub{function} associated  with the (energy
complete) Hamiltonian triple
$\alpha$.

For our gauge theoretical results we will need the following  result
\begin{pr}\label{BigGroupEquiv} Let $G$ be a reductive normal subgroup of a
reductive group
$\hat G$, and let $\hat\alpha:\hat G \times F\ra F$ be a holomorphic action.

Let $K\subset \hat K$ be  maximal compact subgroups of   $G$ and $\hat  G$ and
let $\tau=(K,g,\mu)$ be an energy complete Hamiltonian triple for
$\hat\alpha|_{G\times F}$  such that
$\mu$ is
$\hat K$-equivariant.

Then the corresponding maximal
weight function $\lambda:H(G)\times F\ra
\R\cup\{\infty\}$ is  $\hat G$-\ub{equivariant}, i. e. it satisfies one of the following
equivalent conditions:
\begin{enumerate}
\item  $\hat G$-equivariance: $\lambda^s(x)=\lambda^{\ad_{g }(s)}(g x)$ for all $s\in
H(G)$, $g\in \hat G$.
\item  $\hat G$-parabolic invariance:  $\lambda^s(x)=\lambda^s(h x)$ for every
$ x\in F,\ h\in
\hat G(s)$.
\end{enumerate}
\end{pr}

  We denote by $G_x$ the stabilizer of $x\in F$, and
by $\g_x$ its Lie algebra (the {\it infinitesimal stabilizer} of $x$). A complex Lie subalgebra
$\hg\subset\g$ will be  called {\it reductive} if it is the complexification of the Lie algebra of a
compact subgroup of
$G$.
\begin{dt} Let $\alpha:G\times F\ra F$ be a holomorphic a action of a complex
reductive group, and let $\lambda:H(G)\times F\ra\R\cup\{\infty\}$ be the  maximal
weight function associated with an energy complete Hamiltonian triple $\tau$.  A
point
$x\in F$ will be called
\begin{enumerate}\label{AnStab}
\item  analytically $\tau$-semistable, if $\lambda^s(x)\geq 0$ for every $s\in H(G)$.
\item analytically $\tau$-stable if it is   $\lambda$-semistable and strict inequality
$\lambda^s(x)> 0$ holds for all $s\in H(G)\setminus\{0\}$.
\item   analytically  $\tau$-polystable if its infinitesimal stabilizer $\g_x$ is
reductive,  it is  $\lambda$-semistable and $\lambda^s(x)> 0$ for  every $s\in H(G)$
which is not equivalent to an element in $\g_x$.
\end{enumerate}
\end{dt}
\begin{re}\label{Obv}
If $x$ is analytically $\tau$-semistable then
$$\lambda^s(x)=\langle \mu(x), -is\rangle=0\ \forall  s\in i\kg_x\ .
$$
\end{re}
Indeed, for $s\in i\kg_x$, the path $c^s_x$ is constant, hence
$\lambda^s(x)=\langle \mu(x), -is\rangle$. The right hand term defines a linear
form on  $i\kg_x$ which takes nonnegative values, hence it vanishes.
\begin{re}\label{HermCrit} Under the assumptions  of Definition \ref{AnStab}, one
has:
\begin{enumerate}
   \item   $x$ is  analytically $\tau$-semistable iff $\lambda^s(x)\geq 0$
for all $s\in i\kg$.
\item if $\g_x=\kg_x^\C$, then $x$ is analytically  $\tau$-polystable iff
$\lambda^s(x)=0$ for all $s\in i\kg_x$ and $\lambda^s(x)>0$ for all $s\in
i\kg\setminus i\kg_x$.
\end{enumerate}
\end{re}

Note that if $\g_x\ne\kg_x^\C$, one has no  numerical characterization of
the polystability of $x$ in terms of elements $s\in i\kg$;
on the other hand  one can always choose a Hamiltonian triple
$\tau'=(K',g',\mu')$ {\it conjugate}  to
$\tau$ -- i. e. a triple of the form
$$
(K',g',\mu')=(\Ad_\gamma(K), (\gamma^{-1})^*
g,
\ad_{\gamma^{-1}}^t\circ \mu\circ
\gamma^{-1})
$$
such that
$\g_x={\kg'_x}^\C$; two conjugate triples define the same maximal weight
function  on $H(G)\times F$ (see \cite{Te3}).

The comparison results in \cite{Te3} show that
\begin{pr} Let $\tau=(K,g, \mu)$ be an energy complete Hamiltonian triple. A point
$x\in F$ is symplectically $\tau$-(semi, poly)stable if and only if it is analytically
(semi, poly)stable.
\end{pr}

This result should be regarded as a (non-algebraic) {\it complex geometric
Hilbert criterion}. Note that      the proof in the \ub{semi}stable case is
rather difficult.

In the next section we will give a short proof for the polystable case, which is based
on the continuity method and follows the same strategy as the proof of   our
universal Kobayashi-Hitchin correspondence.

\subsection{The Continuity Method in the finite dimensional case}\label{ContMeth}

In this section   we give an easy proof based on the continuity method for
the following result:
\begin{pr}\label{HilbCrit} Let $\tau=(K,g,\mu)$ be an energy complete Hamiltonian
triple and let
$x\in F$. Then the following conditions are equivalent:
\begin{enumerate}
\item  $x$ is symplectically $\tau$-polystable.
\item $x$  is analytically $\tau$-polystable
\end{enumerate}
\end{pr}

We begin with the following simple
\begin{lm}\label{MomStab} Let $\mu$ be a moment map for an action
$\alpha:K\times F\ra F$ of a compact Lie group on a symplectic manifold $F$. Then
$$\mu(x)\in z_\kg(\kg_x)\ .
$$
\end{lm}
\pf   Regarding $\mu$ as a
map
$F\ra\kg$ via an $\ad_k$-invariant inner product,  we get for any $u\in \kg$
\begin{equation}\label{InfEquiv}d\mu(u^\#_{x})=\frac{d}{dt}|_{t=0}\ \mu(e^{tu}
(x))=\frac{d}{dt}|_{t=0}\ \ad_{e^{tu}}(\mu(x))=[u,\mu(x)]\ .
\end{equation}
The left hand term vanishes when $u\in \kg_x$.
\qed
\\ \\
{\bf Proof} (of Proposition  \ref{HilbCrit}):

We start with the simple implication $1. \Rightarrow
2.$ If
$x$ is symplectically
$\tau$-polystable, let
$x_0\in Gx\cap
\mu^{-1}(0)$. It is easy to see that $\g_{x_0}= \kg_{x_0}^\C$, hence $\g_{x_0}$
(so also
$\g_x$) are reductive subalgebras.  Moreover,  one checks easily
that $\lambda^s(x_0)=0$ for $s\in i\kg_{x_0}$ and $\lambda^s(x_0)>0$ for
$s\in i\kg\setminus i\kg_{x_0}$. Therefore $x_0$ (so also $x$) is analytically
$\tau$-polystable (see \cite{Te3} for details).

Suppose now that $x$ is analytically polystable.
\\ \\
\ub{Step  1}.  Reducing the problem to the case
\begin{equation}\label{CentrStab}
\g_x\subset z (\g,G)\ ,
\end{equation}
where
$$z (\g,G):=\{u\in \g|\ \ad_g(u)=u\ \forall g\in G\} \ .
$$
Note that $z(\g,G)=z(\g)$ when $G$ is connected.

Put
$$G':=Z_G({\g_x})=\{g\in G|\ \ad_{g}(u)=u\ \forall u\in \g_x\}\ .$$
Since $\g_x$ is a reductive subalgebra, it is easy to see that $G'$ is a reductive
subgroup of $G$ with Lie algebra $\g'=z_\g(\g_x)$.

Consider the restricted action $\alpha':=\alpha|_{G'\times F}$ of $G'$ on $F$.

The infinitesimal stabilizer $\g'_x$ of $x$ with respect to $\alpha'$ is $\g_x\cap
\g'$. Since this subalgebra of $\g'$ is contained in $\g_x$,   its element  are
invariant under $\ad_{g}$ for all $g\in G'$, by the definition of this group. Therefore
$\g'_x\subset z (\g',G')$.

Since $\g_x$ is a reductive subalgebra, it follows that, replacing $\tau$ by a
conjugate triple if necessary, we may suppose that
$\g_x=\kg_x^\C$.

We have obviously $G'={K'}^\C$,  $\g'={\kg'}^\C$, where
$$K':=Z_K(\kg_x)\ ,\ {\kg'} =z_\kg(\kg_x)\ .
$$
The triple $\tau'=(K',g,\mu':=\mu|_{\kg'})$ is an energy complete    Hamiltonian
triple for
$\alpha'$. It is easy to see that $x$ is analytically $\tau'$-polystable. Suppose we
have proved  the implication $2.\Rightarrow 1.$ in the case when (\ref{CentrStab})
holds. It will follow that
$x$ is symplectically $\tau'$-polystable, hence there exists $g'\in G'$ such that
$\mu'(g'x)=0$.

We state that, in our case,  the relation $\mu'(g'x)=0$ implies the apparently stronger
relation  $\mu(g'x)=0$.

Indeed, note first, that since $g'\in G'=Z_G(\g_x)$ then
$\g_{g'x}=\ad_{g'}(\g_x)=\g_x=\kg_x^\C$.

By Lemma \ref{MomStab} we get  $\mu(g'x)\in z_\kg({\kg_x})=\kg'$.  But the
relation
$\mu'(g'x)=0$ means that the projection of $\mu(g'x)$ on $\kg'$ vanishes.
Therefore
$\mu(g'x)=0$.
\\ \\
\ub{Step 2}. Proving the implication $2.\Rightarrow 1.$ in the case (\ref{CentrStab})
using the continuity method.

Since $\g_x\subset z (\g,G)$, we   have
\begin{equation}\label{NewStabil}\kg_{gx}=\g_{gx}\cap\kg=\ad_g(\g_x)\cap
\kg=\g_x\cap\kg =\kg_x \ \forall g\in G\ .
\end{equation}

We seek an element $g\in G$ such that
\begin{equation}\label{equation}\mu(g x)=0
\end{equation}

In order to solve this equation, we use the continuity method: we choose a suitable
point $x_0\in Gx$,   we solve for
$\varepsilon\in (0,1]$ the perturbed equation
\begin{equation}\label{newequation}i\mu(e^s x_0)+\varepsilon s=0\ ,\ s\in
i\kg_{x_0}^\bot
\end{equation}
and then we make $\varepsilon\ra 0$.
\\ \\ \\
a) We choose an element $x_0\in Gx$ such that the
perturbed equation (\ref{newequation})  has a solution in
$i\kg_{x_0}^\bot$ for
$\varepsilon=1$.
\\

Set $s_1:=-i\mu(x)$.  Put
$x_0:=e^{-s_1} x$. Then
$$i\mu(e^{s_1} x_0)+s_1=i\mu(x)+s_1=0\ .
$$
By Remark \ref{Obv} we have $s_1\in  i\kg_x^\bot$, so by  (\ref{NewStabil}), we
get as required
$s_1\in i\kg_{x_0}^\bot$.\\ \\ \\
b) We prove that the set
$$T:=\{\varepsilon\in (0,1]|\ \exists s\in i\kg_{x_0}^\bot\ \hbox{such that}\
i\mu(e^s x_0)+\varepsilon s=0\}
$$
is open in $(0, 1]$.
\\

We apply the implicit function theorem to the function
$$l :i\kg_{x_0}^\bot\times(0,1]\map i\kg_{x_0}^\bot\
,\ l(\varepsilon,s)=i\mu(e^s x_0)+\varepsilon s\ .
$$
The fact the $l$ takes values in  $i\kg_{x_0}^\bot$ follows again from
Remark \ref{Obv}. Let $s,\ \dot s\in i\kg_{x_0}^\bot$ and decompose
$$\dot s=\lambda+[k,s]\ ,\ \hbox{with}\ \lambda\in iz_\kg(s)\ ,\ k\in
z_\kg(s)^\bot
$$
as in the proof of Proposition \ref{Differential} in the Appendix.
Put
$$\sigma:= \left((d_s\exp)(\dot s)\right)e^{-s}\ ,\
\sigma_h:=p_{i\kg}(\sigma)\ ,\ \sigma_a:=p_{\kg}(\sigma)\ .$$
Using the known formula (\ref{InfEquiv}), we get
$$\frac{\partial l}{\partial s}(\varepsilon,s)(\dot s)=id\mu(\sigma^\#_{e^s
x_0})+\varepsilon \dot s=id\mu((\sigma^\#_h)_{e^s
x_0})+id\mu((\sigma^\#_a)_{e^s
x_0})+\varepsilon \dot s
$$
$$
= id\mu((\sigma^\#_h)_{e^s
x_0})+[\sigma_a, i \mu( e^s
x_0)]+\varepsilon  \dot s
$$
At a pair $(\varepsilon,s)$ where  $l(\varepsilon, s)=0$, this reduces to
$$\frac{\partial l}{\partial s}(\varepsilon,s)(\dot s)=id\mu((\sigma^\#_h)_{e^s
x_0})-\varepsilon[\sigma_a, s]+\varepsilon  \dot s\ .
$$
Therefore
$$\left\langle\frac{\partial l}{\partial s}(\varepsilon,s)(\dot s),
\sigma_h\right\rangle= \langle d\mu((\sigma^\#_h)_{e^s
x_0}), -i\sigma_h\rangle-\varepsilon\langle[\sigma_a,
s],\sigma_h\rangle+\varepsilon\langle \dot s,\sigma_h\rangle=
$$
$$\omega(-i(\sigma_h^\#)_{e^s x_0}, (\sigma^\#_h)_{e^s
x_0})-\varepsilon\langle[\sigma_a,
s],\sigma_h\rangle+\varepsilon\langle \dot s,\sigma_h\rangle\geq \|(\sigma^\#_h)_{e^s
x_0}\|^2+\varepsilon\|\dot s\|^2
$$
by Proposition \ref{Differential} in Appendix. Therefore $ \frac{\partial
l}{\partial s}(\varepsilon,s)$ is an invertible operator.
\\ \\ \\
c) We prove that the set $T$ is closed in $(0,1]$.

Let $(\varepsilon_n)_n$ be a sequence converging to $\varepsilon_0\in (0,1]$.
Taking the inner product of the equation $l(\varepsilon_n,s_{\varepsilon_n})=0$ with
$s_n:=s_{\varepsilon_n}$ we get
$$\langle \mu(e^{s_n}x_0), -is_n\rangle+  \varepsilon_n\|s_n\|^2=0
$$
We know by Remark \ref{increasing} that the function $t\mapsto \langle \mu(e^{t
s}x_0), -is\rangle$ is always increasing, so
$$\|s_n\|^2=-\frac{1}{\varepsilon_n}\langle \mu(e^{s_n}x_0), -is_n\rangle\leq
-\frac{1}{\varepsilon_n}\langle \mu(x_0),-is_n\rangle\leq
\frac{1}{\varepsilon_n}\|\mu(x_0)\|\|s_n\|\ ,
$$
which shows that the sequence $(s_n)_n$  is bounded.  The limit of a convergent
subsequence will give a solution of the equation $l(\varepsilon_0,\cdot)=0$.
\\ \\ \\
d) Let $(\varepsilon_n)_n$    be a sequence in $(0,1]$ with  $\varepsilon_n\searrow
0$. We prove that any sequence
$(s_n)_n$, $s_n\in i\kg_{x_0}^\bot$  with
$l(\varepsilon_n,s_n)=0$ is bounded.

Taking a subsequence if necessary, it suffices to prove that the condition
$\|s_n\|\ra\infty$ leads to a contradiction. Put
$t_n:=\|s_n\|$, $\sigma_n:=\frac{1}{t_n} s_n$. Taking a subsequence again, we may
suppose that $\sigma_n$ converges to an element $\sigma\in i\kg_{x_0}^\bot$
whose norm will be of course 1.    By the polystability condition,  one has
$$\lambda^\sigma(x_0)=\lim_{t\ra\infty}\langle \mu(e^{t\sigma}
x_0),-i\sigma\rangle>0\ , $$
so  for  $t_0$ sufficiently large it holds
$$ \langle \mu(e^{t_0\sigma}
x_0),-i\sigma\rangle>0\ .
$$
Fix such a $t_0$. For \ub{all} $n$ sufficiently large we will still have
$$ \langle \mu(e^{t_0\sigma_n}
x_0),-i\sigma_n\rangle>0
$$
Assuming $t_n\ra \infty$, we can find $n$ with the above property such that also
$t_n\geq t_0$. By the monotony property Remark \ref{increasing} we would get
$$0< \langle \mu(e^{t_n\sigma_n}
x_0),-i\sigma_n\rangle=\frac{1}{t_n} \langle \mu(e^{s_n}
x_0),-is_n\rangle=-\frac{1}{t_n}\varepsilon_n\|s_n\|^2 \ ,
$$
which is obviously a contradiction.
\\

Now it suffices to note that  the limit of a convergent subsequence of $(s_n)_n$ will
be a solution of the equation $\mu(e^sx_0)=0$.
\qed

\subsection{Maximal weight functions for linear and projective actions}

The explicit form of the maximal weight function associated
with a linear or projective actions is well known (see for instance \cite{Mu1}).  We
explain  below this computation  for completeness.

\subsubsection{Linear  actions}\label{linact}

Let $\rho:G\ra GL(V)$ be a rational linear representation of a connected reductive group
$G$ on a finite dimensional complex vector space $V$.

Let $K$ be a maximal compact subgroup of
$G$ with a fixed invariant inner product on its
Lie algebra $\kg$,  and let
$g$ be a
$K$-invariant Hermitian inner product on $V$.  One has a standard moment
map for the
$K$ action which is given by
$$\mu_0(v)=\rho^*(-\frac{i}{2}v\otimes v^*)\ ,
$$
and any other moment map has the form
$$\mu_\tau=\mu_0-i\tau\ ,
$$
where $\tau\in i z(\kg)$.

Let $\xi\in i\kg$ and let $V=\oplus_{j=1}^k V_j$ be the decomposition of $V$
into eigenspaces of $\xi$. In other words  $\rho_*(\xi) |_{V_j}=\xi_j\id_{V_j}$, where
$\xi_j$ are the pairwise distinct eigenvalues of $\xi$.

Put
$$V^\pm_\xi:=\bigoplus\limits_{\pm\xi_j>0} V_j\ ,\
V_{\pm}^\xi:=\bigoplus\limits_{\pm\xi_j\geq0} V_j\ .
$$

Let $v\in V$. Decompose $v$ as $v=\sum_j v_j$ with $v_j\in V_j$.

One gets
$$\lambda_\tau^\xi(v):=\lim_{t\ra\infty}
\langle\mu_\tau(\rho(e^{t\xi})v),-i\xi\rangle=\left\{
\begin{array}{cc}
+\infty&{\rm if}\ \exists j\in\{1,\dots,k\}\ {\rm such\ that}
\\& \xi_j>0\ {\rm and}\
v_j\ne 0\ ,\\ \\
\langle \tau,\xi\rangle &{\rm otherwise}
\end{array}
\right.
$$

Suppose for simplicity that $\ker(\rho_*)=\{0\}$.

In this case we conclude that a vector   $v\in V$ is $\tau$-stable if and only if for every
$\xi\in i\kg\setminus\{0\}$ with $\langle \tau,\xi\rangle\leq 0$ one has
${\rm pr}_{V^+_\xi}(v)\ne 0$.

Therefore, the non-stable locus is
$$NS_\tau=\union_{\langle \tau,\xi\rangle\leq 0,\xi\ne 0} V^\xi_- \eqno{(*)}
$$

Let us describe this set in the case  when $G$ is a complex torus.   Decompose $V$
as
$$V=\bigoplus_{\chi\in R} V_\chi\ ,
$$
where $R\subset \Hom(i\kg,\R)\subset \g^\vee$ is the finite set of
weights  of the representation $\rho$, and $V_\chi:=\{v\in V\ |\
\rho_*(\gamma)(v)=\chi(\gamma)v\ \forall \gamma\in\g\}$.

For every subset $A\subset R$ put
$$V_A:=\bigoplus_{\chi\in A} V_\chi \ ,$$
Since $G$ is a torus, the space $V^\xi_-$ associated with  any
$\xi\in\i\kg$ has the form $V_A$ for some $A\subset R$. Define
$$C_A:=\{\xi\in i\kg\ |\
V^\xi_-=V_A\}=\{\xi\in i\kg\ |\ \chi(\xi)\geq 0\ \forall\chi\in A,\
\chi(\xi)>0\ \forall \chi\in R\setminus A\}\ .
$$

Note that the sets $C_A$ give a partition of $i\kg$ in pairwise disjoint
polyhedral convex cones.  Using $(*)$ we get
$$NS_\tau=\union_{C_A\cap (H_\tau\setminus\{0\})^{\geq 0}\ne\emptyset}
V_A\ ,
$$
where $H_\tau^{\geq 0}$ denotes the half-space defined by the
inequality $\langle \tau,\xi\rangle\geq 0$.
\begin{co} Suppose that $G$ is a torus.  Then there exists a finite set
$\Phi_\tau\subset i\kg$ such that the following conditions become equivalent:
\begin{enumerate}
\item  $v$ is $\tau$-stable,
\item  $\lambda_\tau^\xi(v)>0$ for all $\xi\in \Phi_\tau$,
\item $v\not\in\union_{\xi\in \Phi} V^\xi_-$.
\end{enumerate}
\end{co}
\pf  The non-empty intersections   $C_A\cap (H_\tau\setminus\{0\})$   define a
finite partition of the pointed half-space
$H_\tau^{\geq 0}\setminus\{0\}$. Take a point in every set of this partition.
\qed

\subsubsection{Projective actions}\label{ProjAct}

Let $\rho:G\ra GL(V)$ be a rational linear representation of a complex reductive group
$G$,  and let
$\alpha:G\times\P(V)\ra\P(V)$ be the induced projective action.

The Fubini-Study K\"ahler form on $\P(V)$ associated with a Hermitian metric $h$ on
$V$ is
$\frac{i}{2\pi}\partial \bar\partial\log\|\cdot\|^2=\frac{1}{4\pi}
dd^c \log\|\cdot\|^2$.

Let $K\subset G$ be a maximal compact subgroup which acts isometrically on the
Hermitian space $(V,h)$. The standard moment map for the $K$-action on $\P(V)$ is
(see \cite{Ki})
$$\mu([v])=\rho^*\left[-\frac{i}{2\pi}\frac{v\otimes\bar v}{\|v\|^2}\right]
$$

The maximal weight of  the Hermitian
endomorphism $f\in\Herm(V)$ with respect to a vector $v\in V$ is defined
$$\lambda_m(f,v):= \max\left\{\lambda\in{\rm Spec}(f)|\ {\rm pr}_{{\rm
Eig}(f,\lambda)}(v)\ne 0\right\}\ .
$$
If $f$ is associated with a one parameter subgroup of $GL(V)$ (i. e. it has
the form $d_1\varphi (1)$, where  $\varphi:\C^*\ra GL(V)\in\Hom(\C^*,GL(V))$),
then all its eigenvalues are integers, so in this case $\lambda_m(f,\cdot)$ takes integer
values. In this case we put $\lambda_m(\varphi,\cdot):=\lambda_m(d_1\varphi
(1),\cdot)$. Writing
$$\varphi(\zeta)=\sum_{n\in R(\varphi)} \zeta^n\id_{V_n^\varphi}
$$
with $V=\bigoplus_{n\in R(\varphi)} V_n^\varphi$, we see that
$$\lambda_m(\varphi,[v])= \max\left\{n\in R(\varphi)|\ {\rm pr}_{V_n^\varphi}(v)\ne
0\right\}\  ,
$$
which is just the maximal weight of  $\varphi$ with
respect to $[v]$ occurring in the algebraic geometric Hilbert criterion.

Let $\xi\in i\kg$.  One has
$$\lambda^\xi([v]):=\lim_{t\ra\infty}
\langle\mu(\rho(e^{t\xi})v),-i\xi\rangle=\frac{1}{2\pi}\lim_{t\ra\infty}\left \langle
\left[\frac{e^{t\xi}v\otimes\overline{e^{t\xi} v}}{\|e^{t\xi}v\|^2}\right], \rho_*(\xi)
\right\rangle= $$
$$=\frac{1}{2\pi} \lambda_m(\rho_*(\xi),v)
$$

Therefore, up to a positive constant, $\lambda^\xi([v])$ is just the maximal
     weight
of the Hermitian endomorphism $\rho_*(\xi)$ with respect to $[v]$.

The point is that, for a rational representation $\rho$, the analytic (semi)sta\-bility of a
point $[v]\in\P(V)$ can be checked by verifying the corresponding inequality only for
one-parameter subgroups
$\C^*\ra G$  (see \cite{Ki}). In other words,       in this case  one  has

\centerline{\it Analytic Stability = GIT Stability = Symplectic Stability\ . }

\section{A ``universal" complex geometric classification problem}

\subsection{Oriented holomorphic pairs}

Let
$$1\map G\textmap{j}\hat G\map G_0\map 1
$$
be an exact sequence of  complex reductive Lie groups;  the closed normal subgroup
$G$ of $\hat G$  will be called the
symmetry group of the variables and the quotient group $G_0$ will be called the
symmetry group of   parameters (see the discussion in the Introduction).

Let $X$ be a connected compact complex manifold,  $\hat Q$   a differentiable
principal $\hat G$-bundle over $X$, and set
$ Q_0:=\hat Q\times_{\hat G}G_0$. Let
$F$ be a K\"ahler   manifold
and
$\hat \alpha$ a holomorphic action of  $\hat G$ on $F$. Consider the associated bundle
$$E:= \hat Q\times_{\hat G}
F\ .$$
We fix a   holomorphic structure $J_0$ on
$ Q_0$, and we denote by ${\cal Q}_0$ the resulting holomorphic
bundle.  We
state the following   universal complex geometric
classification problem:

{\it Classify the pairs $(\hat J,\varphi)$, where $\hat J$ is a
holomorphic structure on the bundle $\hat Q$ which
induces $J_0$ on
$Q_0$ and $\varphi$ is
a $\hat J$-holomorphic section in $E$. The classification is considered
up to isomorphy defined by the natural action of the
complex gauge group}
$${\cal G}:=\Aut_{ Q_0}(\hat Q)=\Gamma(X,\hat Q
\times_{\Ad_{\hat G}}  G  )\ .$$
A pair $(\hat J,\varphi)$ as above will be called {\it a holomorphic
pair of type} $(\hat Q,J_0,\hat\alpha)$.

A possible answer to the classification problem above can
be given by restricting our attention to the subspace of
{\it simple} pairs.  The concept of simple pair can be
defined in our general framework as follows:

For  any pair $(\hat J,\varphi)$ we define the    {\it space of its
infinitesimal automorphisms} (its {\it infinitesimal stabilizer})  by
$$\g_{\hat J,\varphi}:=\{v\in A^0(\hat Q\times_{\ad}
\g)\ \vert\ \bar\partial_{\hat J}v=0,\ v^\#\circ \varphi=0\}\ ,
$$
where $v^\#$ denotes the vertical vector field on $E$ defined by $v$.

\begin{dt} The pair $(\hat J,\varphi)$ is called simple if
$\g_{\hat J,\varphi}=\{0\}\ .
$
\end{dt}

Using non-linear versions of the  methods of [LO], [LT] one should be
able to endow the moduli
space of simple holomorphic  pairs with the structure of a
(possibly non-Hausdorff, possibly singular) complex orbifold (see \cite{Su} for the
case of standard holomorphic pairs).

A more involved approach to the classification problem is based on the
concept of
stability. The data of a Hermitian metric
$g$ on $X$  allows  us to introduce a concept of stability for holomorphic
pairs of a
fixed type $(\hat Q,J_0,\hat \alpha)$, and to construct the corresponding
moduli space  of $g$-stable pairs. This moduli space will be
a Hausdorff open subspace  in the moduli  orbifold of simple
pairs.
\\ \\
{\bf Remark:}  In algebraic geometry and   classical complex geometry one avoids
using   differential geometric  objects (as is  our fixed differentiable bundle $\hat
Q$), so it is worth to note that there is a purely complex geometric (algebraic
geometric) equivalent formulation of our classification problem.

{\it Fix an algebra morphism $\hat\cg:H^*(B\hat G,\Z) \ra H^*(X,\Z)$. Classify the
triples
$(\hat{\cal Q},\iota,\varphi)$, where
$\hat{\cal Q}$ is a
holomorphic $\hat G$ bundle on $X$ with fixed characteristic morphism
$\hat\cg$,
$\iota$ is a holomorphic isomorphism
$\iota:\hat{\cal Q}\times_{\hat G}
G_0\ra{\cal Q}_0$ and $\varphi$ is a holomorphic section of $\hat {\cal
Q}\times_{\hat G} F$.
    The classification is considered up to isomorphism of holomorphic $\hat
G$-bundles.}

A triple $(\hat {\cal Q},\iota, \varphi)$ as above will be  called
a holomorphic triple of type $(\hat\cg,{\cal Q}_0,\hat\alpha)$.  The algebraic version
can be obtained by replacing everywhere the word ``holomorphic" by ``algebraic". This
terminology was used in \cite{OST}.

One can prove easily that a ``moduli set"   of   holomorphic
triples of a fixed type is a discrete disjoint union of  ``moduli sets" of  oriented
holomorphic pairs.

\subsection{The stability condition for universal oriented holomorphic pairs}

\subsubsection{The degree of a  meromorphic $L$-reduction  with respect to  an
ad-invariant linear form}

Let $G$ be a  complex reductive group, and let $L$ be a parabolic subgroup of
$G$. When   $G$ is connected it follows that $L$ has the form
$L=G(s)$, where
$s\in H(G)$.   In the general case we will consider only parabolic subgroups of the the
form $G(s)$ so, for us, a parabolic subgroup will  always be associated to an element in
$s\in H(G)$.    We know that
$L$ fits in the exact sequence
$$1\ra U(L) \ra  L\ra Z(L) \ra 1
$$
which is canonically associated with  $L$, i. e. it does not depend on the presentation
$L=G(s)$. $U(L)$ is just the unipotent  radical of  $L$, and $Z(L)$ is the canonical
reductive quotient of $L$.

 Let $\psi:\lg\ra\C$ be an $\ad_L$-invariant linear form.  Such a form must be a
Lie-algebra morphism\footnote{For a connected Lie group, the converse is
 also true: a
 {\eightmsb C}-valued Lie algebra morphism is an $\ad$-invariant form.}, i. e. it must
fulfill
$$\psi([a,b])=0\ \forall \ a,\ b\in\lg\ .
$$
Recall that
$$\ug(L)=\ug(s):=\bigoplus_{\lambda<0}{\rm Eig}([s,\cdot],\lambda)\ ,\
\zg(L)=\zg(s)=z_\g(s) \ ,
$$
and that $[s,\ug(s)]=\ug(s)$.   This shows that $\psi|_{\ug(s)}=0$. Therefore,
\begin{re}\label{InvForm} Any $\ad_L$-invariant linear form $\psi:\lg\ra\C$ is induced by
an
$\ad_{Z(L)}$-invariant linear form $\zg(L)\ra\C$ on its canonical reductive quotient.
\end{re}

Let now $Z$ by any complex reductive group.  Let $C$ be a maximal compact subgroup
of $Z$ and $T\subset C$ a maximal torus. Consider the root decomposition
$$\zg=\tg^\C\oplus(\bigoplus\limits_{r\in R\setminus\{0\}} \zg_r)
$$
of  the Lie algebra $\zg$.  For any $r\in R\setminus\{0\}$, one has obviously
$\zg_r\subset [\tg^\C,\zg_r]$. Therefore,
\begin{re}\label{restriction}
If $\psi:\zg\ra\C$ is an $\ad_Z$-invariant form, then $\psi|_{\zg_r}=0$ for every
$r\in R\setminus\{0\}$.
\end{re}
\begin{dt} An $\ad$-invariant form $\psi:\zg \ra\C$ will be called real  if for any (and
hence for every) maximal compact subgroup $C$ of $Z$,  the restriction
$\psi|_{i\cg}$ takes real values.
\end{dt}
\begin{dt}\label{InnerProd} A non-degenerate  $\ad_{ G}$-invariant  symmetric
bilinear form
$$h: \g\times  \g\ra\C$$
will be called invariant \ub{complex} \ub{inner} \ub{product} of
\ub{Euclidean} \ub{type} if it induces an Euclidean inner product on any (and hence on
all)  subspaces of the form  $i
\kg$, where $
\kg$ is the Lie algebra of a maximal compact subgroup $K$.
\end{dt}

In the semisimple case, the Killing form is  such an invariant  complex inner product of
Euclidean type.  We can always find such  an complex inner product by choosing an
$\ad_{  K}$-invariant real  inner product  on any subspace $i  \kg$ and extending it in a
$\C$-bilinear way. However, this notion should be regarded as a complex geometric one
(which does not require a distinguished maximal compact subgroup).
\\
\\
{\bf Example:}  Let $h$ be an invariant  complex inner product of Euclidean type on
$\g$.  Any   element
$s\in H(G)$ defines a real   $\ad_{  G(s)}$-invariant form  $h(s):
\g(s)\ra\C$ given by  $u\mapsto h(s,u)$. \\

Note that the restriction of $h$ to $\g(s)$ is  \ub{not} non-degenerate. Its kernel is
$\ug(s)$. Indeed, for any $u,\ v\in\g(s)$, one has
$$h(u,v)=h(\ad_{\exp(ts)}(u),\ad_{\exp(ts)}(v))=h({\rm pr}_{\zg(s)}(u),{\rm
pr}_{\zg(s)}(v))\ ,
$$
because $\ad_{\exp(ts)}(u)\ra {\rm pr}_{\zg(s)}(u)$ as $t\ra\infty$. Therefore, the
formula for $h(s)$ reads
$$h(s)(u)=h(s,{\rm pr}_{\zg(s)}(u))\ ,
$$
which is obviously $\ad_{G(s)}$-invariant, because $G(s)\ra Z(s)$   is a group
morphism,  $h$ is $\ad$-invariant, and $s$ is $\ad_{Z(s)}$-invariant.
\\

By the Chern Weil theory, any
$\ad_Z$-invariant linear form $\psi$ on its Lie algebra
$\zg$ defines  a   characteristic class $c_\psi(W)$ in the 2 -cohomology of the
base of any
$Z$-bundle $W$.   More precisely, if $A$ is a a connection on $W$, one
sets
$$c_\psi(W)=\left[\psi(\frac{i}{2\pi} F_A)\right]_{\rm DR}\ .
$$
For a $L$-bundle $\Lambda$ we put
$$c_\psi(\Lambda):=c_\psi(\Lambda\times_L Z(L))\ .
$$
\begin{dt}\label{ddegree}
The  degree  of a  $Z$-bundle $W$ (respectively of a $L$-bundle $\Lambda$)  on
a compact $n$-dimensional K\"ahler manifold
$(X,g)$ with respect to the real $\ad$-invariant form $\psi:\zg \ra\C$ (respectively
$\psi:\zg(L) \ra\C$)  is defined by
$$\deg_g(W,\psi):=\langle c_\psi(W)\cup [\omega_g]^{n-1},[X]\rangle\
(\deg_g(\Lambda,\psi):=\deg_g(\Lambda\times_L Z(L)))
$$
\end{dt}

In the non-K\"ahler case, one can still define the  degree with respect to $\psi$ of a
holomorphic $Z$ bundle ($L$-bundle), but the obtained invariant will no longer be a
topological invariant.  One has to use the same strategy as in the well known case of line
bundles on Gauduchon manifolds.

 Let $g$ be a Gauduchon metric on
$X$ and let
${\cal W}$ be a {\it holomorphic}
$Z$-bundle on
$X$.
Choose a maximal compact subgroup $C$ of $Z$,  an $C$-reduction $R$ of ${\cal W}$
and let   $A_{R}$ its Chern connection (see section \ref{Chern}). The point is that, as in
the case of line bundles, the following important result holds:
\begin{lm} The real closed  2-form
$\psi(\frac{i}{2\pi} F_{A_R})$ is independent of the chosen $C$-reduction $R$
up to a
$i\partial\bar\partial$-exact form. Therefore, the $i\partial \bar\partial$-cohomology
class defined by $\psi(\frac{i}{2\pi} F_{A_R})$   depends  only on the holomorphic
$Z$-bundle
${\cal W}$.
\end{lm}
\pf Let $R_0$ and $R$ be two $C$-reductions of ${\cal W}$. We can
write
$$R=e^{-\frac{s}{2}}(R_0)\ ,
$$
where $s\in A^0(R_0\times_{\ad} i\cg)$. By Corollary \ref{newcurvature}, one has
$$F_{ A_R}=F_{A_{R_0}}+\bar\partial( e^{-s}\partial_{A_{R_0}} (e^s))\ ,
$$

It is easy to see that
$$\psi(e^{-s}\partial_{A_{R_0}} (e^s))=\partial (\psi(s))\ .
$$
This follows easily using Proposition \ref{LocDiag} and  formula (\ref{ds}), taking into
account Remark \ref{restriction}. Therefore
$$\psi(F_{ A_R})=\psi(F_{A_{R_0}})+\bar\partial\partial \psi(s)\ .
$$
\qed

This allows us to define
\begin{dt}\label{defdeg} Let ${\cal W}$ be a holomorphic $Z$-bundle, ${\cal L}$   a
holomorphic $L$-bundle on a compact $n$-dimensional Gauduchon manifold $(X,g)$, and
let
$\psi:\zg\ra\C$ (respectively $\psi:\zg(L)\ra\C$) be a real $\ad$-invariant form.
Let $C$ be a maximal compact subgroup of $Z$ (respectively $Z(L)$). We define
$$c_\psi({\cal W})=\left[\psi(\frac{i}{2\pi}
F_{A_R})\right]\ \in \ \qmod{Z^{1,1}(X,\R)}{i\partial \bar\partial A^0(X,\R)}\ , \
c_\psi({\cal L}):=c_\psi({\cal L}\times_L Z(L))
$$
$$\deg({\cal W},\psi):=\left\langle c_\psi({\cal
W}),\omega_g^{n-1}\right\rangle=\int\limits_X
\psi(\frac{i}{2\pi} F_{A_R})\wedge\omega_g^{n-1}\ ,\
$$
$$\deg({\cal L},\psi):=\deg({\cal
L}\times_L Z(L),\psi)\ ,
$$
where $A_R$ is the Chern connection of any $C$-reduction  of ${\cal W}$ (
of ${\cal L}\otimes_L Z(L)$).
\end{dt}
\vspace{3mm}

We will need an extension of the degree map to meromorphic $L$-reductions of a
holomorphic $G$-bundle. We begin with the following

\begin{dt}\label{MeromRed} Let $L$ be a parabolic subgroup of a complex reductive
group
$G$ and ${\cal Q}$ a holomorphic $G$-principal bundle over a complex manifold
$X$. A meromorphic reduction of  $G$ to $L$ is a meromorphic section
in the associated
bundle
$${\cal Q}_L:=\qmod{{\cal Q}}{L}={\cal Q}\times_G\left(\qmod{G}{L}\right)\ ,$$
i. e. a closed reduced irreducible subspace $\rho\subset {\cal Q}_L$ which
generically projects isomorphically onto $X$ via the bundle projection
$q_L:{\cal Q}_L\ra X$.
\end{dt}

Since the fibre $G/L$ of ${\cal Q}_L$ is a projective manifold, it follows
easily that the maximal open subset $X_\rho\subset X$ over which $\rho$ is the
graph of a holomorphic section is the complement of a Zariski closed subset of
codimension at least 2.  On $X_\rho$ $\rho$ defines a holomorphic $L$-reduction
${\cal Q}^\rho\subset {\cal Q}|_{X_\rho}$.
\\ \\
{\bf Example:}  Let ${\cal Q}$ be a holomorphic $GL(r,\C)$-bundle and ${\cal E}$
the associated holomorphic vector bundle.  The data of a meromorphic reduction
of
${\cal Q}$ to the parabolic group
$$L:=\left\{\left.\left(\matrix{A&B\cr0&C}  \right)\right|\ A\in
GL(p,\C),\ B\in M(p,r-p),\ C\in GL(r-p,\C)\right\}
$$
is equivalent to the data of a rank $p$ subsheaf ${\cal F}\subset {\cal E}$
{\it with torsion free quotient} ${\cal H}$. Such a subsheaf is necessarily reflexive
and defines a subbundle on the subset $X \setminus {\rm Sing}({\cal H})$.
Given such a subsheaf ${\cal F}$, the  corresponding meromorphic section
$\rho$ is given by
$$\rho:=\overline{\im(f)}
$$
where $f:X \setminus {\rm Sing}({\cal H})\ra {\G}r_p({\cal E})$ is   defined by
the holomorphic subbundle ${\cal F}|_{X \setminus {\rm Sing}({\cal H})}$ of
${\cal E}|_{X
\setminus {\rm Sing}({\cal H})}$.\\

Suppose now that ${\cal Q}$ is a holomorphic $G$-bundle on the compact
Gauduchon   manifold $(X,g)$,  $L$   a parabolic subgroup of $G$,
$\rho\subset {\cal Q}_L$  a meromorphic $L$-reduction of ${\cal Q}$ and
$\psi:\zg(L)\ra\C$ a real
$\ad_{Z(L)}$-invariant form. We denote by $q_L:{\cal Q}_L\ra X$ the natural
projection and by $\iota_\rho:\rho\hookrightarrow {\cal Q}_L$ the inclusion map.

We claim that the assignment
\begin{equation}\label{current}
A^{1,1}({\cal Q}_L)\ni\varphi\mapsto \int\limits_{\rho\setminus{\rm
Sing}(\rho)}
\varphi\wedge q_L^*(\omega_g^{n-1})
\end{equation}
defines a $\partial\bar\partial$-closed current of bidimension $(1,1)$ on ${\cal Q}_L$.
Indeed, let $\Rg\textmap{\alpha}\rho$ any desingularization of $\rho$. Then the
assignment
\begin{equation}\label{desing}A^{1,1}({\cal Q}_L)\ni\varphi\mapsto \int\limits_{\Rg}
 (\iota_\rho\circ \alpha)^*(\varphi)\wedge ( q_L\circ
\alpha)^*(\omega_g^{n-1})
\end{equation}
obviously defines a $\partial\bar\partial$-closed current of bidimension $(1,1)$ on
${\cal Q}_L$, because the pull-back  $( q_L\circ \alpha)^*(\omega_g^{n-1})$ is a
$\partial\bar\partial$-closed form on the smooth manifold $\Rg$.

But the assignment (\ref{desing}) coincides with (\ref{current}) -- which is
independent of the desingularization $\alpha$ -- because
$\rho\setminus{\rm Sing}(\rho)$ is identified via $\alpha$ with an open subset of total
measure of $\Rg$.

We will denote the current defined by the formulae (\ref{current}), (\ref{desing}) by
$q_L^*(\omega_g^{n-1})|_\rho$.

Since, by definition, ${\cal Q}_L$ is the base of the holomorphic $L$-bundle ${\cal
Q}\ra{\cal Q}_L$, we can define
\begin{dt}\label{DefDeg}
$$\deg(\rho,\psi):=\left\langle c_\psi({\cal Q}\ra{\cal
Q}_L),q_L^*(\omega_g^{n-1})|_\rho\right\rangle\ .
$$
\end{dt}

In the presence of a $K$-reduction $P\subset {\cal Q}$ of the $G$-bundle  ${\cal Q}$,
one can give an explicit formula for the degree map.

Note first that the intersection $K_L:=L\cap K$ projects isomorphically onto a maximal
  compact          subgroup of
  $Z(L)$ via the canonical epimorphism $L\ra Z(L)$. This follows from formula
(\ref{intersection}) in Appendix.

Therefore,   $P\ra{\cal Q}_L\simeq P/K_L$ can
be regarded  as a  reduction of the associated
$Z(L)$-bundle ${\cal W}:={\cal Q}\times_L Z(L) \ra{\cal Q}_L$
 to its maximal compact subgroup $K_L$.

It is not difficult to express the Chern connection $B$ of the $K_L$-bundle $P\ra
{\cal Q}_L$ in terms of the Chern connection $A$ of the original  $K$-bundle $P\ra X$.
Denoting by $\omega_A$, $\omega_B$ the corresponding connection forms, one has
$$\omega_B={\rm pr}_{\lg\cap\kg}(\omega_A)
$$
where ${\rm pr}_{\lg\cap\kg}$ denotes the projection $\kg\ra  \kg_L$ defined by
the canonical  $\ad_{K_L}$-invariant decomposition
\begin{equation}\label{dec}\kg=\kg_L\oplus\kg_L^\bot\ ,
\end{equation}
where  $\kg_L^\bot:={\rm pr}_{\kg}(\ug(L))$. Therefore, on $P$
one can write
$$\omega_A=\omega_B+\alpha\ ,\  {\rm
pr}_{\kg_L}\Omega_A=\Omega_B+\frac{1}{2}{\rm
pr}_{\kg_L}[\alpha\wedge \alpha]\ ,
$$
where $\alpha\in A^1(P,\kg_L^\bot)$ is the $\kg_L^\bot$-component of
$\omega_A$.  Regarding the forms in the second formula as bundle valued forms on the
base $Q_L$ one gets
$${\rm pr}_{P\times_{K_L}\kg_L} [q_L^*(
F_A)]=F_B+\frac{1}{2}{\rm pr}_{P\times_{K_L}\kg_L}[a\wedge
a]\ ,
$$
where $a\in A^1({\cal Q}_L, P\times_{K_L}\kg_L^\bot)$ is the bundle valued 1-form on
${\cal Q}_L$ associated with $\alpha$.

Using the definitions above, we get
$$\deg(\rho,\psi)=\int\limits_{\rho\setminus{\rm
Sing}(\rho)}   \psi(\frac{i}{2\pi} F_B)\wedge
 q_L^*(\omega_g^{n-1})
=\int\limits_{X_\rho}\psi(\frac{i}{2\pi}\rho^*(F_B))\wedge\omega_g^{n-1}\ ,
$$
where $\rho^*$ is the pull-back associated with the holomorphic map $X_\rho\ra {\cal
Q}_L$ defined by $\rho$.  The point is that the form  $ \rho^*(F_B)$ is just the
curvature of the connection $A_\rho$ induced by $A$ on the $K_L$-subbundle
$Q^\rho\cap P\subset P|_{X_\rho}$, whereas $a_\rho:=\rho^*(a)$ is just the second
fundamental form of this subbundle.  Therefore we get the following
\begin{pr}\label{DegForm}   If  $g$ be a Gauduchon metric on $X$, then
\begin{equation}\label{DegreeFormula}
\deg_g(\rho,\psi)=\int\limits_{X}\psi(\frac{i}{2\pi}F_{A_\rho})\wedge
\omega_g^{n-1}\ ,
\end{equation}
where $A_\rho$  is the connection induced by the Chern connection $A$ of
any $K$-reduction $P\subset {\cal Q}$ on the
$K_L$-subbundle $P^\rho:=Q^\rho\cap P\subset P|_{X_\rho}$\ .
\end{pr}

  Using the canonical decomposition $\g=\zg(L)\oplus\zg(L)^\bot$ obtained by
complexifying the decomposition (\ref{dec}), one can extend $\psi$ in a natural way to
a linear form  $\psi:\g\ra\C$ (which of course is not $\ad_G$-invariant in general).

Therefore the degree formula becomes
\begin{equation}\label{degree-new}\deg_g(\rho,\psi)=
\int\limits_X\left[\psi(\frac{i}{2\pi}F_{A})-\psi(\frac{i}{4\pi}[a_\rho\wedge
a_\rho])\right]
\wedge\omega_g^{n-1}\ ,
\end{equation}
where $a_\rho$ is the second fundamental form of $Q^\rho\cap P$ in $P|_{X_\rho}$.

The degree map can be extended to the case of an arbitrary Hermitian metric $g$ in
the following way (see \cite{LT} for the case of vector bundles):

Recall first that the maximum principle holds for the operator
$P=i\Lambda_g\bar\partial\partial$ for any Hermitian metric $g$. Consequently, one
gets a direct sum decomposition
$${\cal C}^\infty(X,\R)=\R\oplus P({\cal C}^\infty(X,\R))
$$
which is $L^2$-orthogonal when $g$ is Gauduchon. Dualizing, we get a direct sum
decomposition
$${\cal D}'(X,\R)=\R\oplus P({\cal D}'(X,\R))
$$
Using this decomposition,  we define  the operator   $I_g:{\cal
D}'(X,\R)\ra
\R$ (the generalized integral) by
$$I_g(\varphi):=Vol_g(X) {\rm pr}_\R(\varphi)\ .
$$
$I_g$ vanishes always on $\im(P)$ and coincides with the standard integral   when $g$
is Gauduchon.

\begin{re}\label{HermCase} (The Hermitian case) Let $g$ be a Hermitian metric.
Replacing the usual integral by the operator
$I_g$ in the formulae (\ref{DegreeFormula}), (\ref{degree-new}) one  gets a well defined
holomorphic invariant $\deg_g(\rho,\psi)$ for pairs $(\rho,\psi)$ consisting  of a real
$\ad$-invariant linear form $\psi$ on $\zg(\lg)$ and a meromorphic reduction $\rho$
of ${\cal Q}$ to the parabolic subgroup $L$.

The degree with respect to a Hermitian metric $g$ coincides with the degree with
respect to the unique Gauduchon metric $g'$ in its  conformal class defined
by the condition $\int_X  (g'/g)^{n-1}vol_g=Vol_g(X)$. Therefore, the properties of
the degree map associated with general Hermitian metrics reduce easily to the
Gauduchon case.
\end{re}

\subsubsection{Stability, Semistability, Polystability}\label{SSSP}

We come back to our gauge theoretical problem, so consider  a
differentiable principal
$\hat G$-bundle $\hat Q$ on a compact Hermitian manifold $(X,g)$.

Let $\hat \alpha:\hat G\times F\ra F$ be a holomorphic action, $G$ a closed
normal subgroup of $\hat G$. We choose a   K\"ahler structure  $g_F$ on $F$,  a
maximal compact   subgroup $K$ of
$G$ acting isometrically on $(F,g_F)$, a maximal compact subgroup $\hat K$ of $\hat
G$ which contains $K$, and a $\hat K$-equivariant moment map $\mu$ for the induced
$K$-action.   We will suppose that the Hamiltonian  triple $\sigma=(K,g_F,\mu)$ is
energy complete, so that it defines a generalized maximal weight function
$\lambda:F\times H(G)\ra
\R\cup\{\infty\}$ with  the fundamental properties listed in Proposition
\ref{FundProp} and
\ref{BigGroupEquiv}.  We also recall that $\lambda$ depends only on the conjugacy
class of $\sigma$ (see section \ref{AnStabSympStab}).

Let $E$ be the associated bundle $E:=\hat Q\times_{\hat G} F$. Consider again
$\zeta\in H(G)$   and
a  $\hat G(\zeta)$-reduction $\hat Q^\rho\subset\hat Q|_{X_\rho}$ defined on an
open subset
$X_\rho\subset X$.   For    points $x\in
X_\rho$,   $e\in E_x$
we  put
$$\lambda^\zeta (e,\rho):=\lambda^\zeta(\varphi_q(e))\ ,
$$
where $\varphi_q:E_x\ra F$ is the holomorphic isomorphism defined by an
element $q\in \hat Q^\rho_x$. By Proposition \ref{BigGroupEquiv}  it follows that
the term on
the right does not depend on the choice of $q\in \hat Q^\rho_x$.

We fix an $\ad$-invariant complex inner product $h$ of Euclidean type on $\hat\g$ (see
Definition \ref{InnerProd}).

\begin{dt}\label{stability} Fix a holomorphic structure $J_0$ on $Q_0$.
A holomorphic pair $(\hat J,\varphi)$ of type  $(\hat Q,J_0,\hat \alpha)$   is called
$\sigma$-\ub{semistable} if, for every  $\xi\in H(G)$   and for every meromorphic
reduction $\rho$ of the holomorphic bundle $\hat {\cal Q}_{\hat J}$ to $\hat
G({\xi})$  it holds
$$\frac{2\pi}{(n-1)!}\deg_g(\rho,  {h(\xi)})+\int\limits_{X_{\rho}}
\lambda^\xi (\varphi,\rho)vol_g\geq 0\ .
$$
A pair $(\hat J,\varphi)$ of type  $(\hat Q,J_0,\hat \alpha)$  is called
$\sigma$-\ub{stable} if it is
$\sigma$-semistable and  strict inequality holds for any $\xi \in
H(G)\setminus
 \{0\} $.
\end{dt}

The left hand term in the inequality in Definition \ref{stability} will be called the
 \ub{total} \ub{maximal } \ub{weight} of the pair $(\hat J,\varphi)$ with respect to
the pair $(\xi,\rho)$.

\begin{re}
  In general, the (semi)stability condition   depends on  the Hermitian
metric $g$, the  Hamiltonian triple $\sigma=(K,g_F,\mu)$ and the invariant complex
inner product of Euclidean type $h$.
\end{re}

Note that,  for any $\hat g\in\hat G$,  the translation $R_{\hat g}:\hat{\cal
Q}\ra\hat{\cal Q}$ defines a biholomorphic diffeomorphism  $\hat{\cal Q}/{\hat
G(\psi)}\ra\hat{\cal Q}/ {\hat G(\ad_{\hat g^{-1}}(\psi))}$, hence it assigns to a
$\hat G(\psi)$ meromorphic reduction $\rho$ of $\hat {\cal Q}$ a ${\hat
G(\ad_{g^{-1}}(\psi))}$ meromorphic reduction
 $\rho \hat g$. The corresponding total maximal weights of a pair $(\hat
J,\varphi)$ coincide, so we get the following important

\begin{re}\label{representatives} It suffices to check  the (semi)stability condition for  a
single representative in every equivalence class of $H(G)$ modulo the adjoint action $
\ad_{\hat G}$ of the group $\hat G$.

In particular it suffices to check this inequality for any element $\psi\in i C $,    where
$C\subset \tg$ is a closed Weyl chamber in the Lie algebra $\tg$ of a maximal torus of
a maximal compact subgroup $K\subset G$.
\end{re}

We recall that a Lie subalgebra $\bg\subset\g$ is called reductive if it is the
complexification of the Lie algebra of a compact  subgroup of $G$.

\begin{dt}\label{RedAlg} Let $\hat Q$ be a   $\hat G$-bundle on
a connected compact complex manifold $X$. A finite dimensional Lie subalgebra
$\ug\subset A^0(\hat Q\times_\ad
\g)$ will be called \ub{reductive} if there exists a reductive Lie subalgebra
$\bg\subset\g$ and a    $Z_{\hat G}(\bg)$-reduction $\hat
Z\hookrightarrow\hat Q$ of $\hat Q$ such that $\ug$ is the image of $\bg$
via the obvious  morphism
$$\bg\hookrightarrow  A^0(\hat Z\times_\ad\bg)\ra A^0(\hat Z \times_\ad\g)=
A^0( \hat Q\times_\ad \g)\ .
$$
\end{dt}

Such a reduction gives an isomorphism ${-}_{\hat Z}:\bg\ra \ug$.
Note that the reduction $\hat Z$ can be recovered from the   isomorphism
$b_{\hat Z}$, because
\begin{equation}\label{recover}
\hat Z=\{q\in\hat Q|\ b_{\hat Z}(q)=b\ \forall b\in\bg\}\ .
\end{equation}
In particular one has
\begin{re}\label{HolRed} If $\hat {\cal Q}$ is a holomorphic $\hat G$-bundle on $X$ and
$\ug\subset H^0(X,\hat {\cal Q}\times_\ad \g)$ is a reductive Lie subalgebra
consisting of holomorphic sections, then any
$Z_{\hat G}(\bg)$-reduction $\hat Z\subset\hat {\cal Q}$ which induces an
isomorphism $\bg\ra\ug$ is holomorphic.
\end{re}

\begin{dt}\label{polystable}

A pair $(\hat J,\varphi)$ of type  $(\hat Q,J_0,\hat \alpha)$  is called
$\sigma$-\ub{polystable} if:
\begin{enumerate}
\item its infinitesimal stabilizer $\g_{\hat J,\varphi}$  is a reductive
 subalgebra of
$A^0(\hat Q\times_{\ad}\g)$,\\
\item  it is
$\sigma$-semistable and  the equality
$$\frac{2\pi}{(n-1)!}\deg_g(\rho,  {h(\xi)})+\int\limits_{X_{\rho}}
\lambda^\xi (\varphi,\rho)vol_g=0
$$
in the semistability condition holds  when and only when $\rho$ is induced by a
global \ub{holomorphic}
$ Z_{\hat G}(\xi)$-reduction
$ \hat{\cal Z}(\xi)\hookrightarrow \hat{\cal Q}_{\hat J}$, and the section
$$\xi_{\hat{\cal Z}(\xi)}\in H^0(X,\hat {\cal Q}_{\hat J}\times_{\ad}\g)$$
 defined by $\xi$ via this
reduction belongs to the infinitesimal stabilizer $\g_{\hat J,\varphi}$.
\end{enumerate}
  \end{dt}

\begin{dt}\label{relcent} Let $z(\g,\hat G)\subset z(\g)$ be the centralizer of $\hat
G$ in
$\g$, i. e. the subalgebra consisting of those elements of $\g$ which are fixed under the
adjoint $\hat G$-action. In other words, $z(\g,\hat G)$ is the Lie algebra of
$Z(G,\hat G):=G\cap Z(\hat G)$.
\end{dt}
When $\hat G$ is connected, one   has  $z(\g,\hat G)=z(\g)$.
Note that $z(\g,\hat G)$ can  always be  identified with a finite dimensional abelian
subalgebra of
$A^0(\hat Q\times_\ad\g$).  For pairs $(\hat J,\varphi)$  with infinitesimal
stabilizer   contained in $z(\g,\hat G)$ the polystability condition becomes:
\begin{re}\label{centralcase}  Let $\beta=(\hat J,\varphi)$ be a  pair with
infinitesimal stabilizer  $\g_\beta$ contained in $z(\g,\hat G)$. Then $\beta$ is
polystable if and only if it is semistable and the inequality in
the semistability  condition is strict unless $\xi\in H(G)\cap \g_\beta$, in which
case one must have $\hat G(\xi)=\hat G$, $\hat {\cal Q}_{\hat J}^\rho=\hat {\cal
Q}_{\hat J}$.
\end{re}

We show now that, by a suitable  reduction of the structure group of
the bundle, one can  reduce the study of \ub{any} polystable pair   $\beta=(\hat
J,\varphi)$ to the case  $\g_{\beta}\subset z(\g,\hat G)$ (compare with Step
1 in the proof of Theorem \ref{HilbCrit} in the finite dimensional framework).

Let $\beta$ be a $\sigma$-polystable pair,   $\bg\subset\g$ a
reductive subalgebra and  $\hat Z\subset \hat   Q$ a
$Z_{\hat G}(\bg)$-reduction of $\hat Q$ such that the infinitesimal
stabilizer $\g_{\beta}$ is identified with  $\bg$ via this reduction (see Definition
\ref{RedAlg}).     Put
$$\hat G':=Z_{\hat G}(\bg)\ ,\ G':=Z_G(\bg)=\hat G'\cap G\ ,\
 G'_0:=\qmod{\hat G'}{G'}\subset G_0\ ,\ \hat \alpha':=\hat\alpha|_{\hat G'\times
F}\ .
$$

By Remark \ref{HolRed} it follows that   $\hat Z$ is a holomorphic
reduction of $\hat{\cal Q}_{\hat J}$; we denote by $\hat J'$ the induced holomorphic
structure.

Put $ Z_0:=\hat Z\times_{\hat G'} G'_0$, and note that $Z_0$ is a holomorphic
$G'_0$-subbundle of ${\cal Q}_0$, so it has an induced holomorphic structure $J'_0$.

The restriction of the maximal weight function $\lambda:H(G)\times F\ra
\R\cup\{\infty\}$ to $ H(G')\times F\ra
\R\cup\{\infty\}$ is the maximal weight function $\lambda'$ of a suitable
Hamiltonian triple $(K',g',\mu')$ for  the restricted action
$\alpha':=\alpha|_{F\times\hat G'}$, where  $\mu'$ is equivariant with respect to
a maximal compact subgroup  $\hat K'\supset K'$ of $\hat G'$.

\begin{pr}\label{CentStab} If $\beta=(\hat J,\varphi)$ is $\sigma$-polystable, then the
associated pair
$\beta':=(\hat J',\varphi)$ is a
$\sigma'$-polystable pair of type
$(\hat Z,J'_0,\hat\alpha')$. The infinitesimal stabilizer $\g'_{\beta'}$ of this pair is
$\bg\cap z_\g(\bg)=z(\bg)$ which is    contained in the centralizer
$z(\g',\hat G')$.
\end{pr}
\pf Consider the holomorphic bundles  $\hat {\cal Q}:=\hat{\cal Q}_{\hat J}$, $\hat
{\cal Z}:=\hat {\cal Z}_{\hat J'}$ defined by the holomorphic structures $\hat J$,
$\hat J'$.  It is easy to see that
$(\hat J',\varphi)$ is
$\sigma'$-semistable  because,   for any
$s\in H(G')$  any meromorphic $\hat G'(s)$-reduction $\rho'$ of $\hat{\cal Z}$
extends to a meromorphic $\hat G(s)$-reduction $\rho$ of $\hat {\cal Q}$,
and the two total maximal weights coincide.

The delicate part is the fact that $\beta'$ is $\sigma'$-polystable.  Let
$\rho'$ be a meromorphic  $\hat G'(s)$-reduction of   $\hat{\cal Z} $ such that the
corresponding total maximal weight vanishes, and let $\rho$ be the associated $\hat
G(s)$-reduction of $\hat {\cal Q} $.  We know that $\rho$ is induced by a
global holomorphic $Z_{\hat G}(s)$-reduction ${\cal Z}(s) \subset \hat {\cal Q}$ such
that the section $\sigma:=s_{\hat {\cal Z}(s)}$  defined by $s$ via this reduction
belongs to the infinitesimal stabilizer $\g_{\beta}$.  We claim that, under these
assumptions, one has
\begin{enumerate}
\item $s$ coincides with the element $b\in\bg$ defined by the
equality $\sigma=b_{\hat{\cal Z}}$.
\item  $ \hat G'(s)=\hat G'$, in particular the meromorphic reduction $\rho'$ is just
$\hat {\cal Z}\stackrel{=}{\hookrightarrow} \hat {\cal Z}$.
\end{enumerate}

In order to prove this, note first that
\begin{equation}\label{Z(s)}
\hat {\cal Z}(s)=\{q\in\hat{\cal Q}|\ \sigma(q)=s\}\ .
\end{equation}
On the other hand,   $\hat {\cal Z}(s)|_{X_\rho}\subset\hat
{\cal Q}^\rho\subset \hat {\cal Z}|_{X_\rho}\cdot \hat G(s)$, because,   by
construction, $\hat {\cal Q}^\rho$ is the $\hat G(s)$-saturation of a $\hat
G'(s)$-bundle contained in $\hat {\cal Z}$.   Since $X_\rho$ is dense, the
inclusion $\hat {\cal Z}(s) \subset \hat {\cal Z} \cdot \hat G(s)$ holds on whole
$X$. Therefore, any $q\in \hat Z(s)$ has the form $z g$ with $z\in \hat {\cal Z}$ and
$g\in \hat G(s)$.  We get, for any $q\in \hat {\cal Z}(s)$
$$s=\sigma(q)=\sigma(zg)=\ad_{g^{-1}}(\sigma(z))=\ad_{g^{-1}}(b)\ ,
$$
where $b\in\bg$ is the element which corresponds to $\sigma\in \g_{\hat
J,\varphi}$ via the reduction $\hat {\cal Z}$.

Therefore $b=\ad_g(s)$ with $g\in \hat G(s)$.   Now recall that $\ad_{\hat
G(s)}(s)=s+\ug(s)$.  But $b\in z_{\hat \g}(s)$ (because
$s\in z_{\hat\g}(\bg)$). We get $b=s$. The second statement follows from the first.
\qed

\section{Hermitian-Einstein pairs}

\subsection{The Hermitian-Einstein equation}

Choose a maximal compact subgroup $K$ of $G$ and a maximal compact subgroup
$\hat K$ of $\hat G$ which contains $K$. Let
$\hat P\subset \hat Q$ be a $\hat K$-reduction of $\hat Q$, and let
$P_0$ be the bundle associated with $\hat P$   with
structure group   $K_0:=\qmod{\hat K}{K}$.  Consider the
(unique)   integrable connection  $A_0\in{\cal A}^{1,1}(P_0)$
which corresponds to $J_0$ via the Chern correspondence. Using the Chern
correspondence (see section \ref{Chern}), our  classification problem associated with the
complex  data
$(\hat Q, J_0,\alpha)$ can be reformulated as follows:
\\ \\
{\it Classify the pairs $(\hat A,\varphi)$, where $\hat A$ is an integrable
connection on  $\hat P $ which induces a fixed integrable connection $A_0$ on
$P_0$,
and $\varphi$ is a $J_A$-holomorphic section in the associated bundle $E=\hat
P\times_{\hat K} F$. The classification is considered up to isomorphy defined
by the action of the complex gauge group
$${\cal G}:=\Gamma(X,\hat P\times_{\Ad_{\hat K}}
   G)=\Aut_{Q_0}({\hat Q})\ .$$
This action is  induced -- via the Chern correspondence -- by the
${\cal G}$-action on the space of holomorphic pairs.}

A pair $(\hat A,\varphi)$ as above will be called {\it an integrable
pair of type} $(\hat P,A_0,\alpha)$.

Note that, formally, ${\cal G}$ can be regarded as the  complexification of
the real gauge group
$${\cal K}:=\Gamma(X,\hat P\times_{\Ad_{\hat K}}
   K)=\Aut_{P_0}({\hat P})\ .
$$

We denote  by ${\cal A}_{A_0}(\hat P)$ the space of connections on
$\hat P$ which
induce $A_0$ on $P_0$.

Recall that we have fixed  an ad-invariant  complex inner product of Euclidean  type $h$ on
$\hat
\g$ (see Definition \ref{InnerProd}).
  This inner product defines an orthogonal
   projection ${\rm pr}_{i\kg}:i\hat\kg\ra i\kg$ which induces a bundle projection
$$\hat P\times_{\ad}  i\hat \kg\map  \hat P\times_{\ad}  i  \kg
$$
which will be denoted by the same symbol.

Using the inner product defined by $h$ to $\kg$, we can view
$\mu$ as a map with values in $\kg$.

\begin{dt}\label{HE}  An   pair $(\hat A,\varphi)$ of type $(\hat
P,A_0,\hat\alpha)$
is called
$\mu$-Hermitian-Einstein if it is integrable and   solves the equation
$${\rm pr}_{i\kg}\left[i\Lambda_g F_{\hat A} \right]+i\mu(\varphi)=0\ , \eqno{(\rm
HE)}
$$
which will be  called the generalized Hermitian-Einstein equation associated with the
data  $(\hat  P,A_0,\hat\alpha,\mu,g)$.

Let $(\hat J,\varphi)$ be a holomorphic pair of type $(\hat Q, J_0,\hat\alpha)$. A
$\hat K$-reduction $\hat P$ of $\hat Q$ will be called $\mu$-Hermitian-Einstein  if the
associated pair  $(A_{\hat P,\hat J},\varphi)$ is $\mu$-Hermitian-Einstein.
\end{dt}
Here we denoted by $A_{\hat P,\hat J}$ the Chern connection of the pair
$(\hat P,\hat J)$ (see section \ref{Chern}).

  Note that this equation is  ${\cal G}$-invariant, so that one can
consider the moduli
space of solutions
$${\cal M}^{\rm HE}:=
\qmod{[{\cal A}_{A_0}(\hat P)\times \Gamma(E)]^{\rm HE} }{{\cal G}}\ ,
$$
where $[{\cal A}_{A_0}(\hat P)\times \Gamma(E)]^{\rm HE} \subset[{\cal
A}_{A_0}(\hat P)\times \Gamma(E)]$ is the subspace of (integrable) pairs
satisfying the equation (HE).

Note that any section $v \in \Gamma(X,\hat Q\times_{\ad} \g)$ induces a
section
$v^\#$ in the vertical tangent bundle $V_E$ of $E$.

\begin{dt}  A pair $(\hat A,\varphi)$ of type $(\hat P,A_0,\alpha)$ is called
  irreducible if its infinitesimal stabilizer
$$\kg_{(\hat A,\varphi)}:=\{v\in A^0(\hat P\times_{\ad} \kg)|\   \nabla_{\hat
A}v=0,\ v^\#\circ\varphi=0\}
$$
(which is   a sub-Lie-algebra of $A^0(\hat P\times_{\ad}\kg)$) vanishes.
\end{dt}

\subsection{Pairs which allow Hermitian-Einstein reductions are polystable}

We come back to the assumptions and notations of section \ref{SSSP}. The aim of this
section is to prove the following
\begin{thry}\label{SimImp}   Let $(\hat
J,\varphi)$ be a holomorphic pair of type $(\hat Q,J_0,\alpha)$.
Suppose that there exists a $\hat K$-reduction $\hat P\subset \hat Q$ such that the
associated pair $( A_{\hat P,\hat J},\varphi)$  is
$\mu$-Hermitian-Einstein.
Then $(\hat J,\varphi)$ is
$\sigma$-polystable, and it is stable if only if $(\hat A_{\hat P,\hat J},\varphi)$ is
irreducible.
\end{thry}
\pf

We will assume that $g$ is Gauduchon (see Remark \ref{HermCase}). We check first
the inequality in the definition of $\sigma$-semistability.  By Remark \ref
{representatives} it suffices to check this  inequality for elements
$\xi\in i\kg$.

Let    $\rho$ be
a meromorphic reduction of
$\hat Q$ to
$\hat G(\xi)$, and put $\hat P^\rho:=\hat Q^\rho\cap \hat P$. The inclusion
$\iota:\hat P^{\rho}\hookrightarrow\hat P|_{X_\rho}$ defines a reduction of the bundle
$\hat P|_{X_\rho}$ to the subgroup
$\hat K({\xi}):=\hat G({\xi})\cap
\hat K$ of $K$.  It is easy to see that
$$\hat K({\xi})^\C=Z_{\hat G}(\xi) \ .
$$
 Therefore $\xi$ is invariant under the structure group
of $\hat P^\rho$, so it defines  a section
$s(\rho,\xi)$ in
$A^0(\hat P|_{X_\rho}\times_{\ad} i\kg)$.

Put $\hat A:=A_{\hat P, \hat J}$. We get
$$0=\int_{X_\rho} \left\langle {\rm pr}_{i\kg}\left[i \Lambda_g
F_{\hat A} \right]+i\mu(\varphi),s(\rho,\xi)\right\rangle vol_g=$$
$$=\int_{X_\rho}  \left\langle  i  \Lambda_g F_{\hat A} ,s(\rho,\xi)\right\rangle vol_g+
\int_{X_\rho}  i\mu(\varphi)(s(\rho,\xi)) vol_g\ .
$$

Recall now that the functions $t\mapsto \mu_\xi(e^{t\xi})$ involved in
the definition of $\lambda $ are monotonely increasing, so that, in
every point $x\in X_\rho$ we get
$$i\mu(\varphi)(s(\rho,\xi))(x)\leq\lambda ^{s(\rho,\xi)}(\varphi)(x)\ .
$$
This shows that
$$\int_{X_\rho}  i\mu(\varphi)(s(\rho,\xi)) vol_g\leq
\int\limits_{X_{\rho}}
\lambda^{s(\rho,\xi)}(\varphi)vol_g=\int\limits_{X_{\rho}}
\lambda^{\xi}(\varphi)vol_g\ ,
$$
with equality if and only if $s(\rho,\xi)^\#\circ\varphi = 0$ over
$X_\rho$.

The adjoint bundle $\hat P\times_{\ad} \hat \kg=\hat P^{
\rho}\times_{\ad}\hat \kg$ of
$\hat P$ splits as
$$\hat P^{
\rho}\times_{\ad}\hat \kg=\ad(\hat P^{
\rho})\oplus \left[\hat P^{
\rho}\times_{\ad}[\hat \kg({\xi})]^{\bot}\right]\ .
$$

The connection $\hat A$ on $\hat P$ induces a connection $\hat
A_{\rho}$ on $\hat P^{\rho}$ whose horizontal space at a point
$p\in  \hat P^{\rho}$ is the intersection
$$\left[H_{\hat A,p}\oplus p\cdot [\hat\kg({\xi})]^{\bot}\right]\cap
T_p(\hat P^{\rho})\ ,
$$
where $H_{\hat A,p}$ stands here for the $\hat A$-horizontal space.
The connection form of this   connection  is
$$\omega_{\hat A_{\rho}}={\rm pr}_{\ad(\hat P^{\rho})}
\iota^*(\omega_{\hat A})
$$
  The difference $a:=\hat A-\iota_*(\hat
A_{
\rho})$ is a section of the subbundle $\hat P^{
\rho}\times_{\ad}[\hat
\kg({\xi})]^{\bot}$ of $\ad(\hat P)$ (the second fundamental form of the
subbundle $\hat P^{ \rho}$ of $\hat P$ with respect to the connection $\hat
A$).

Decompose the Lie algebra $\hat \g$ as
$$\hat \g=\hat \g_+\oplus\hat\g_-\oplus\hat\g_0
$$
where $\hat\g_\pm$ is direct sum of the eigenspaces of the
endomorphism
$$[\xi,\cdot] \in\End(\hat\g)$$
  corresponding to the
positive (negative) eigenvalues, and $\hat\g_0=\ker [\xi,\cdot]$.  Note that
$$\hat\g_0=  [\hat \kg({\xi})]^\C,  \
\hat\g_-\oplus\hat\g_0=\hat\g({\xi}),\ \hat\g_+=[\hat
\kg({\xi})]^{\bot}\ ,\ \hat\g_+=\overline{\hat\g}_- \eqno{(1)}
$$
where the conjugation in the last formula is taken with respect to the real
structure $\hat\g=i\hat\kg\otimes\C$ on $\hat\g$.

Therefore, the second fundamental form $a$ decomposes as
$$a=a^+ + a^-
$$
where $a_\pm\in A^1(X,\hat P^\rho\times_\ad\g_\pm)$ and
\begin{equation}\label{ObvId}
a^-=-\bar a^+\ .
\end{equation}
   But $\hat\g_+=T_{[e]}\left[\qmod{\hat G}{\hat G(\xi)}\right]$ and on $X_\rho$
$$\hat P^\rho\times_\ad\hat\g_+=\hat
Q^\rho\times_\ad\hat\g_+=\rho^*\left(T^V\left(\hat Q\times_{\hat
G}\left[\qmod{\hat G}{\hat G_\xi}\right]\right)\right)\ ,
$$
(where $T^V$ denotes as usually the vertical tangent bundle of a fibration).   Therefore
the holomorphy of
$\rho$ is equivalent to the condition
\begin{equation}\label{Hol}
a^+\in A^{10}(X,\hat P^\rho\times_\ad\hat\g_+)
\end{equation}
(compare with Proposition \ref{AlmostHol}).
Write
$$a=\sum_{\lambda\in {\rm Spec}   { [\xi,\cdot] }}a_\lambda$$
 so that
$$a^\pm=\sum_{\lambda\in {\rm Spec}^\pm  { [\xi,\cdot] }}a_\lambda \ ,$$
and take into account that
$$[{\rm Eig}([\xi,\cdot],\lambda),{\rm Eig}([\xi,\cdot],\mu)]\subset {\rm
Eig}([\xi,\cdot],\lambda+\mu)\ ,
$$
$$\overline{{\rm Eig}([\xi,\cdot],\lambda)}={\rm Eig}([\xi,\cdot],-\lambda)\ .
$$
By (\ref{ObvId}), (\ref{Hol}) it follows that
\begin{equation}\label{a-wedge-a}{\rm pr}_{\hat \kg_{\xi}} [a\wedge
a]=\sum_{\lambda\in {\rm Spec}^+ [\xi,\cdot]} [a_\lambda \wedge
a_{-\lambda} ]+\sum_{\lambda\in {\rm Spec}^- [\xi,\cdot]} [a_\lambda \wedge
a_{-\lambda} ]=
$$
$$
=-2\sum_{\lambda\in {\rm Spec}^+
[\xi,\cdot]} [a_\lambda \wedge\overline{a}_{\lambda} ]\ ,
\end{equation}
where $[\cdot\wedge\cdot]$ is the (symmetric) bilinear multiplication
obtained by
multiplying the Lie algebra bracket with the wedge product on 1-forms:
$$[u\wedge v](x,y):=[u(x),v(y)]-[u(y), v(x)]\ .
$$
Comparing the curvatures of the two connections $\hat A$, $\hat
A_\rho$, we get
\begin{equation}\label{curvatures}
{\rm pr}_{\hat \kg_{\xi}}F_{\hat A}=F_{\hat
A_\rho}+{\rm pr}_{\hat \kg_{\xi}}\frac{1}{2}[a\wedge a]\ .
\end{equation}
But, since $\langle\cdot ,\cdot\rangle$ is ad-invariant, we get for any
$u,v\in\hat
\g$ with
$u\in {\rm Eig}([\xi,\cdot],\lambda)$
\begin{equation}\label{prod-wedge}
\langle[u,v],  \xi \rangle^\C =-\langle v,[u, \xi]\rangle^\C=
\langle[\xi,u],v\rangle^\C= \lambda\langle u,v\rangle^\C\ .
\end{equation}
Therefore, by the degree formula given in Proposition \ref{DegForm}, one obtains
\begin{equation}\label{Kob}\int\limits_{X_\rho}  \left\langle  i \Lambda_g F_{\hat
A},s(\rho,\xi))\right\rangle vol_g= \frac{2\pi}{(n-1)!}\deg(\rho,
 {h(\xi)})- \sum_{\lambda\in {\rm Spec}^+
[\xi,\cdot]}\lambda\nr a_\lambda\nr^2_{L^2}\ ,
\end{equation}
hence $\frac{2\pi}{(n-1)!}\deg(\rho,
 {h(\xi}))\geq \int_{X_\rho}  \langle  i \Lambda_g F_{\hat
A},s(\rho,\xi)\rangle vol_g$ with equality if and only if
$\alpha=0$, i. e. if
and only if  the connection $\hat A$ reduces to a $\hat
K({\xi})$-connection, which   implies that $s(\rho,\xi)$   is an $\hat A$-parallel section
of
$\hat P\times_{\ad}i\kg$.

Therefore we have shown that $(\hat A,\varphi)$ is $\sigma$-semistable,   and
that, for $\xi\in\i\kg$, one has
$$  \frac{2\pi}{(n-1)!}\deg_g(\rho,
 {h(\xi)})+\int\limits_{X_{\rho}}
\lambda^\xi (\varphi,\rho)vol_g = 0
$$
if and only if  $\rho$ is induced by an $\hat A$-parallel $\hat
K(\xi)$-reduction of
$\hat P$ and the associated section $s(\rho,\xi)\in A^0(\hat
P\times_{\ad}i\kg)$
belongs to $i \kg_{ \hat A,\varphi }$, where $\kg_{ \hat A,\varphi }\subset A^0(\hat
P\times_{\ad}\kg)$ is the
infinitesimal stabilizer algebra of the pair $(\hat A,\varphi)$ with respect to the   action of
the real gauge group
${\cal K}:=\Gamma(\hat P\times_{\Ad} K)$. But an
$\hat A$-parallel $\hat K(\xi)$-reduction of $\hat P$ defines a holomorphic
$Z_{\hat G}(\xi)$-reduction of $\hat {\cal Q}_J$. Therefore, the last condition
in the definition of polystability is verified.

Using the same method as in the proof of Proposition 1.3
\cite{Te3}
 one can easily check, using the maximum principle for the operator
$P=i\Lambda\bar\partial\partial$, that
$$\g_{ \hat J,\varphi }=\kg_{ \hat A,\varphi }\otimes\C
$$
 (compare with Theorem 2.2.1 in \cite{LT}).  Taking into account
that the elements of
$\kg_{
\hat A,\varphi }$ are $\hat  A$-parallel sections, we get an $\hat A$-parallel
reduction of
$\hat P$ to
$Z_{\hat K}(\cg)$, where $\cg\subset\kg$ is a subalgebra of $\kg$ isomorphic to
$\kg_{ \hat A,\varphi }$.   This gives a holomorphic  $Z_{\hat
G}(\cg\otimes\C)$-reduction of $\hat{\cal Q}_{\hat J}$, so    that
$\g_{ \hat J,\varphi }$ is a reductive subalgebra of
$A^0(\hat Q\times_{\ad}\hat \g)$ in the sense of
Definition \ref{RedAlg}, which proves
the first condition in the definition of polystability.
\qed
\begin{re} The proof of the theorem  shows that,  for a  Hermitian-Einstein
pair $(\hat A,\varphi)$,  the following conditions are equivalent:
\begin{enumerate}
\item  $(\hat A,\varphi)$ is irreducible.
\item  The associated holomorphic pair $(\hat J,\varphi)$ is simple.
\end{enumerate}
\end{re}

\section{Polystable   pairs  \  allow Hermitian--Einstein   reductions}\label{DiffImpl}

\begin{dt}\label{reductions} We will denote by ${\cal R}(\hat Q)$ the space of all $\hat
K$-reductions of $\hat Q$  and, for a fixed $K_0$-reduction $P_0$ of $Q_0$, we will
denote by ${\cal R}_{P_0}(\hat Q)$ the space of all  $\hat K$-reductions of $\hat Q$
which project onto $P_0$.
\end{dt}

Our final purpose is the following theorem:
\begin{thry}\label{purpose} Let $(\hat J,\varphi)$ be a $\sigma$-polystable pair of type
$(\hat Q,{\cal Q}_0,\hat \alpha)$. For any $K_0$-reduction $P_0$ of $Q_0$ there
exists a reduction $\hat P\in{\cal R}_{P_0}(\hat Q)$ which is      $\mu$-
Hermitian-Einstein.
\end{thry}

In the proof, which will be completed at the and of section \ref{DiffImpl}, we will
assume again for simplicity   that
$g$ is a Gauduchon metric; the statement is true for a general Hermitian
metric (see Remark \ref{HermCase}).

\subsection{The perturbed equation}

Let $(\hat {J},\varphi)$ be a $\sigma$-polystable pair of type $(\hat
Q,J_0)$, and let
$P_0$ be a fixed
$K_0$-reduction of $Q_0$. We  seek a
$\hat K$-reduction $\hat P\subset \hat Q$ which projects on $P_0$
such that  the
pair
$(  A_{\hat P,\hat J}, \varphi)$ satisfies the  HE equation
(see Definition \ref{HE})), where $A_{\hat P,\hat J}$ is the Chern connection of the pair
$(\hat P,\hat J)$.

By Proposition \ref{CentStab} we may suppose that\\ \\
{\bf Assumption}: {\it  The  infinitesimal stabilizer $\g_\beta$ of the polystable pair
$\beta=(\hat J,\varphi)$    is contained in the
centralizer $z(\g,\hat G)$.}

This assumption has  the following important consequence:
\begin{re}\label{CplxStab}
If this assumption is satisfied, the infinitesimal stabilizer of $\beta$
$$\kg_\beta:=\g_\beta\cap A^0(\hat P\times_{\ad} \kg)$$
 with respect to the action of the real
gauge group ${\cal K}:=\Gamma(X, \hat P\times_{\Ad} K)$ is independent of the
$\hat K$-reduction $\hat P$, and it can be identified with a subalgebra
of $z(\kg,\hat K)$.
Moreover,
for any $\hat K$-reduction
$\hat P$ of $\hat Q$, one has
\begin{equation}\label{condition}\g_\beta=[\g_\beta\cap
A^0(\hat P^0\times_{\hat K}\kg)]\otimes\C=\kg_\beta\otimes\C\ .
\end{equation}
\end{re}

Indeed,
one has
$$\g_\beta\cap
A^0(\hat P^0\times_{\hat K}\kg)=\g_\beta\cap z(\g,\hat G)\cap
A^0(\hat P\times_{\hat K}\kg)=\g_\beta\cap z(\kg,\hat K)\ .
$$
The connected component $G_\beta^0$ of $e$ in the stabilizer
$$G_\beta\subset {\cal G}:=\Gamma(\hat Q\times_{\rm Ad}G)$$
of $\beta$ is a connected closed complex subgroup of ${\cal G}$ with
Lie algebra
$\g_\beta\subset z(\g,\hat G)$, so it is contained in the connected
component $Z(G,\hat G)^0$ of $e$ in the intersection
$Z(G,\hat G):=G\cap Z(\hat G)$.

Now it is easy to see that, for every connected closed complex subgroup
$H$ of a complex
torus $\Theta$, one has $\hg=[\hg\cap\tg]\otimes\C$, where $\tg$ is
the Lie algebra of
the {\it unique} maximal compact subgroup $T$ of $\Theta$.

Therefore, in our case one gets
$$\g_\beta=[\g_\beta\cap z(\kg,\hat K)]\otimes\C=[\g_\beta\cap
A^0(\hat P^0\times_{\hat K}\kg)]\otimes\C\ .
$$
\qed

We fix a $\hat K$-reduction $\hat P\in{\cal R}_{P_0}(\hat Q)$ and seek a solution
$\hat P\subset
\hat Q$ of our problem of the form
$$\hat P_s=e^{-\frac{s}{2}}(\hat P)\ ,
$$
where $s\in A^0(\hat P \times_{\ad} i\kg)$. Put $h:=e^s$. The new
(pointwise) moment
map
$\mu:\Gamma(\hat Q\times_{\hat\alpha} F)\ra A^0(\hat
P_s\times_{\ad}\kg)$ associated
with the reduction $\hat P_s$ is given by the formula
\begin{equation}\label{newmoment}
\mu_s(\varphi)=\ad_{e^{-\frac{s}{2}}}(\mu (e^{\frac{s}{2}}\varphi))\ ,
\end{equation}
where $\mu $ is the moment map with respect to the initial reduction $\hat P$.

Denote by $A^0$ the Chern connection of the pair $(\hat P ,\hat J)$. The
Hermitian-Einstein equation (Definition \ref{HE}) becomes
\begin{equation}\label{s-eq}
{\rm pr}_{i\kg}\left[i\Lambda_g (F_{A^0}+\bar\partial( e^{-s}\partial_0
(e^s)))\right]+i\ad_{e^{-\frac{s}{2}}}(\mu (e^{\frac{s}{2}}\varphi))=0\ ,
\end{equation}
or, equivalently,
\begin{equation}\label{ss-eq}
{\rm pr}_{i\kg}\ \ad_{e^{\frac{s}{2}}}\left[i\Lambda_g (F_{A^0}+\bar\partial(
e^{-s}\partial_0
(e^s)))\right]+i (\mu (e^{\frac{s}{2}}\varphi))=0\ .
\end{equation}

Our perturbed equation is
\begin{equation}\label{s-eq-p}
{\rm pr}_{i\kg}\left[i\Lambda_g (F_{A^0}+\bar\partial(
e^{-s}\partial_0(e^s))\right]+i\ad_{e^{-\frac{s}{2}}}(\mu (e^{\frac{s
}{2}}\varphi))+\varepsilon
s=0\ ,
\end{equation}
which, since $s$ is $\ad_{e^{\frac{s}{2}}}$ -- invariant, is equivalent to
\begin{equation}\label{ss-eq-p}
\ad_{e^{\frac{s}{2}}}{\rm pr}_{i\kg}\left[i\Lambda_g
(F_{A^0}+\bar\partial( e^{-s}\partial_0(e^s)))\right]+i
(\mu (e^{\frac{s}{2}}\varphi))+\varepsilon s=0\ .
\end{equation}

\begin{pr}\label{initialisation} There exists a $\hat K$-reduction
$\hat P \in{\cal
R}_{P_0}(\hat Q)$
  such that, denoting by $A^0$ the corresponding Chern connection $\hat A_{\hat P,\hat
J}$, one has:
\begin{enumerate}
\item[A.] \begin{equation}\label{step1a}
{\rm pr}_{A^0(X,i\kg_\beta)} \left[{\rm pr}_{i\kg}[i\Lambda_g  F_{\hat A^0}]
+i\mu (\varphi)\right]=0
\end{equation}
\item[B.] Equation (\ref{s-eq-p}) has a solution $\zeta \in A^0(\hat
P^0\times_{\ad}i\kg_\beta^\bot)$ for
$\varepsilon =1$.
\end{enumerate}
\end{pr}
\pf

Step 1. We start with an  arbitrary reduction $\hat P^1\in {\cal
R}_{P_0}(\hat Q)$, and we
seek a smooth function $v\in{\cal C}^\infty(X,i\kg_\beta )\subset {\cal
C}^\infty(X,iz(\kg))$ such that the Chern connection $\hat
A^2:=A_{\hat P^2,\hat J}$  associated with the reduction $\hat
P^2:=e^{-\frac{v}{2}}(\hat P^1)$ satisfies (\ref{step1a}).

By Proposition \ref{newconnection}, Corollary \ref{newcurvature} one has
$$\hat A^2:=\hat A^1+\partial v\ ,\ F_{\hat A^2}:=F_{\hat
A^1}+\bar\partial\partial v\ .$$

The new moment map with respect to the reduction $\hat P^2$ is given by
$$\mu_2(\varphi)=\ad_{e^{-\frac{v}{2}}}\
\mu_1(e^{\frac{v}{2}}\varphi)=\mu_1(\varphi)\ ,
$$
because $v$ is central and belongs pointwise to the infinitesimal
stabilizer of $\varphi$.
Therefore (\ref{step1a}) becomes
\begin{equation}\label{Laplace}
{\rm pr}_{A^0(X,i\kg_\beta)} \left[{\rm pr}_{i\kg}[i\Lambda_g  F_{\hat A^1}]
+i\mu_1(\varphi)\right]+i\Lambda\bar\partial\partial v =0\ .
\end{equation}

On the other hand, applying the polystability condition  to central elements of     the
form  $u\in i\kg_\beta$ one gets that
\begin{equation}\label{orthogonality}
\langle {\rm pr}_{i\kg}[i\Lambda_g  F_{\hat A^1}]
+i\mu_1(\varphi), u\rangle_{L^2}=0 \ \forall u\in i\kg_\beta\
\end{equation}
(compare with Remark \ref{Obv} in the finite dimensional framework).
This implies that the Laplace equation (\ref{Laplace}) has solutions, so we get a
reduction $\hat P^2$ with the desired property.

Step 2.

Put
$$\zeta:=-{\rm pr}_{i\kg}[i\Lambda_g  F_{\hat A^2}] -i\mu_2(\varphi)
$$
and $\hat P :=e^{\frac{\zeta}{2}}(\hat P^2)$. By  Step 1   it follows that
$\zeta\in A^0(\hat P^2\times_{\ad}i\kg_\beta^\bot)$.  Since $e^{\frac{\zeta}{2}}$
commutes with $\zeta$, one can also regard $\zeta$ as an element of $A^0(\hat
P \times_{\ad}i\kg)$  which will be also pointwise orthogonal to $\kg_\beta$.  By
Corollary
\ref{newcurvature} and (\ref{newmoment}) one has
$$-\zeta ={\rm pr}_{i\kg}\left[i\Lambda_g (F_{A^2}
\right]+i (\mu_2( \varphi))=$$
$$={\rm
pr}_{i\kg}\left[i\Lambda_g (F_{A^0}+\bar\partial( e^{-\zeta}\partial_0
(e^{\zeta})))\right]+i\ad_{e^{-\frac{\zeta}{2}}}
(\mu_0(e^{\frac{\zeta}{2}}\varphi))
\ ,
$$
which shows that,  with these choices, B holds.
\qed

\subsection{A priori estimates for the solution $s_\varepsilon$}

\begin{lm}\label{uniform} Fix a $\hat K$-reduction   $\hat P$ as in
Proposition
\ref{initialisation}. Let
$$s \in A^0(\hat P \times_{\ad}i\kg_\beta^\bot)$$
be a section solving the equation
$$
\ad_{e^{\frac{s}{2}}}{\rm pr}_{i\kg}\left[i\Lambda_g
(F_{A^0}+\bar\partial( e^{-s}\partial_0(e^s)))\right]+i
 \mu (e^{\frac{s}{2}}\varphi)+\varepsilon s=0 \ .\eqno{(e_\varepsilon)}
$$
Put
$$k^0:={\rm pr}_{i\kg}\left[i\Lambda_g
(F_{A^0}\right]+i\mu (\varphi) \ .$$
Then one has:

\begin{enumerate}
\item ${1\over 2}\cdot P(\vert s{\vert^2}) +
\varepsilon\cdot\vert s\vert^2 \leq \vert{k^0}\vert\cdot\vert s\vert\ ;$

\item $m: = {\sup_X}\vert s\vert \leq
{1\over\varepsilon}\cdot{\sup_X}\vert{k^0}\vert\ ;$

\item $m \leq C\cdot(\Vert s{\Vert_{L^2}} +
{\sup_X}\vert{k^0}\vert)^2$\ \ where the
constant $C$ only depends on $g$ and the fixed $\hat K$-reduction $\hat P$.
\end{enumerate}
\end{lm}
\pf Taking the  pointwise inner product  with respect to $\hat P $ of both
sides of $(e_\varepsilon)$ with
$s$, one gets
\begin{equation}\label{(eq.s)}
\left({\rm pr}_{i\kg}[i\Lambda_g
 F_{A^0}],s\right)+\left(i\Lambda_g\bar\partial(
e^{-s}\partial_0(e^s))),s \right)+
2\mu ^{-i\frac{s}{2}}(e^{\frac{s}{2}}\varphi)+\varepsilon|s|^2=0\ .
\end{equation}
On the other hand, note that
$$\mu^{-i\frac{s}{2}}(e^{\frac{s}{2}}\varphi)=\mu^{-i\frac{s}{2}}(\varphi)+
\eg^\frac{s}{2}(\varphi)\ ,
$$
where $\eg^s(\varphi)$ is the $\R_{\geq 0}$-valued function which maps
every $x\in X$
to the energy of the curve $[0,1]\ni t\mapsto e^\frac{ts}{2}\varphi(x)\in F_x$.

Therefore,
\begin{equation}\label{FirstIneq}
\left(i\Lambda_g\bar\partial(
e^{-s}\partial_0(e^s))),s \right)+\varepsilon|s|^2=-\left({\rm
pr}_{i\kg}[i\Lambda_g
(F_{A^0}],s\right)-2\mu^{-i\frac{s}{2}}(e^{\frac{s}{2}}\varphi)\leq
$$
$$-\left({\rm pr}_{i\kg}[i\Lambda_g
 F_{A^0}],s\right)-2\mu^{-i\frac{s}{2}}(\varphi)=-\left({\rm
pr}_{i\kg}[i\Lambda_g
 F_{A^0}]  +  i\mu  (\varphi),s\right)=-(k^0,s)\ .
\end{equation}

We introduce now the \ub{real}
\ub{analytic} function  $\eta:\R\ra\R_{>0}$   defined by
$$\eta(t)=\left\{
\begin{array}{ccc}
\sqrt{\frac{1-e^{-t}}{t}}&\rm if& t\ne 0\\
1&\rm if&t=0
\end{array}
\right.
$$
(see section \ref{UsefulFormulae}).
We use   the formalism explained in section \ref{UsefulFormulae}. By
Proposition
\ref{identity} the inequality (\ref{FirstIneq}) becomes
\begin{equation}\label{NewIneq}
\frac{1}{2}
P(|s|^2)+|\eta([s,\cdot])(\partial_0(s))|^2+\varepsilon|s|^2\leq-(k^0,s)\ .
\end{equation}
This proves the claim.\\ \\
2. follows from 1. by the Maximum Principle. \\ \\
3.  From 1. we deduce
$$
    P(\vert s{\vert^2}) \leq 2\cdot\vert{k^0}\vert \vert s\vert \leq
    \vert s\vert^2 + \vert{k^0}\vert^2 \leq
    \vert s\vert^2 + {\sup_X}\vert{k^0}{\vert^2}\ ,
$$
so, since ${\sup_X}\vert{k^0}{\vert^2}$ is a constant, Lemma 3.3.2 [LT] implies
$$\begin{array}{cc}
    m &\leq C\cdot\left(\Vert\phantom{.}\vert s{\vert^2}{\Vert_{L^1}}
+   {\sup_X}\vert{k^0}{\vert^2}\right) =
    C\cdot\left(\Vert s{\Vert^2_{L^2}} +
{\sup_X}\vert{s^0}{\vert^2}\right)\\
    &\leq C\cdot\left(\Vert s{\Vert_{L^2}} +
{\sup_X}\vert{k^0}\vert\right)^2 \ .
\end{array}$$
\qed

\subsection{Solving the equation $(e_\varepsilon)$ for $\varepsilon\in (0,1]$.}

Let
$$l(\varepsilon,s):=\ad_{e^{\frac{s}{2}}}{\rm pr}_{i\kg}\left[i\Lambda_g
(F_{A^0}+\bar\partial( e^{-s}\partial_0(e^s)))\right]+i
 \mu(e^{\frac{s}{2}}\varphi)+\varepsilon s$$
 be the left  hand term of the equation
$(e_\varepsilon)$.

We  state first the following  regularity result:
Choose $p\in \N$ such that the Sobolev space  $L^p_k$ is an
algebra for $k\geq 1$.
\begin{lm}\label{regularity} Let $s\in  L^p_2(\hat P
\times_{\ad}i\kg_\beta^\bot)$ such that $l(\varepsilon,s)=0$ for some
$\varepsilon\in\R$. Then
$s\in {\cal C}^\infty(\hat P
\times_{\ad}i\kg_\beta^\bot).$
\end{lm}
\pf Recall first that, by Proposition \ref{NewConn}, one has
\begin{equation}\label{firstterm}
\ad_{e^{\frac{s}{2}}}
(F_{A^0}+\bar\partial( e^{-s}\partial_0(e^s)) )=F_{A_s}\ ,
\end{equation}
where $A_s:= e^{-\frac{s}{2}}(A^0+ e^{-s}\partial_0(e^s))\in {\cal A}(\hat P)$. The
components of the linear connection defined by $A_s$ on the   bundle  $\hat
P\times_{\ad}\g$ are
\begin{equation}\label{dAs}
\bar\partial_{A_s}=\ad_{e^{\frac{s}{2}}}\circ\bar
\partial_{A^0}\circ\ad_{e^{-\frac{s}{2}}}\ ,\  \partial_{A_s}=
\ad_{e^{-\frac{s}{2}}}\circ
\partial_{A^0}\circ\ad_{e^{\frac{s}{2}}}\ ,
\end{equation}
which can be rewritten as
$$\bar\partial_{A_s}=
\bar\partial-\bar\partial(e^{\frac{s}{2}})e^{-\frac{s}{2}}\ ,\
 \partial_{A_s}=
 \partial_0 +e^{-\frac{s}{2}} \partial_0(e^{\frac{s}{2}}) \ .
$$
Here the forms $\bar\partial(e^{\frac{s}{2}})e^{-\frac{s}{2}}$, $e^{-\frac{s}{2}}
\partial_0(e^{\frac{s}{2}})$ operate on $\hat
P\times_{\ad}\g$ via the adjoint representation $a\mapsto[a,\cdot]$.
On the other hand, by formulae (\ref{d(exp(s))}), (\ref{newd(exp(s))}) in the
Appendix,  these forms  can be expressed in terms of $\bar\partial
s$, $\partial_0 s$ as
\begin{equation}\label{forms}
\bar\partial(e^{\frac{s}{2}})e^{-\frac{s}{2}}=\frac{1}{2}
\psi\left(-\left[\frac{s}{2},\cdot\right]\right) (\bar
\partial s)
\ ,\ e^{-\frac{s}{2}}
\partial_0(e^{\frac{s}{2}})=\frac{1}{2}
\psi\left(\left[\frac{s}{2},\cdot\right]\right) (
\partial_0 s)\ ,
\end{equation}
where $\psi$ is the   positive real analytic function
$$\psi(t)=\left\{\begin{array}{ccl}
\frac{e^t-1}{t}&\rm if& t\ne 0  \\
1&\rm if& t = 0\ .
\end{array}\right.
$$

Let $\Psi: \hat P\times_{\ad} i\kg \ra \Herm(\hat P\times_{\ad}\g)$
be the (non-linear!) real analytic    bundle map defined by
$$a\mapsto \frac{1}{2}
\psi\left(\left[\frac{a}{2},\cdot\right]\right)\ .
$$
We get
\begin{equation}\label{As}\bar\partial_{A_s}=
\bar\partial-\Psi(-s)(\bar\partial s)\ ,\
 \partial_{A_s}=
 \partial_0 +\Psi(s)(\partial_0 s )\ .
\end{equation}
This gives  for the $(1,1)$ curvature form of $A_s$
\begin{equation}\label{Fform} F_{A_s}=F_{A^0}+d_{A^0}\left( \Psi(s)(\partial_0
s )-\Psi(-s)(\bar\partial s)\right)+
$$
$$+\left ( \Psi(s)(\partial_0 s )-\Psi(-s)(\bar\partial
s)\right)\wedge \left( \Psi(s)(\partial_0 s )-\Psi(-s)(\bar\partial
s)\right)  =$$
$$
=2\Psi(s)(\bar\partial\partial_0
s)+(1-\Psi(-s))(F_{A^0})+\bar\partial(\Psi(s))\wedge\partial_0 s
-\partial_0(\Psi(-s))\wedge \bar\partial s
$$
$$
-\Psi(s)(\partial_0 s )\wedge\Psi(-s)(\bar\partial s)-\Psi(-s)(\bar\partial
s)\wedge \Psi(s)(\partial_0 s )\ .
\end{equation}
Therefore the equation $(e_\varepsilon)$ is equivalent to an equation of the form
$$i\Lambda\bar\partial\partial_0 s={\cal F}(\varepsilon,s, d_{A^0} (s))\ .
$$
The left composition lemma applies on the right, more precisely,   ${\cal
F}(\varepsilon,s, d_{A^0} (s))$ belongs to $L^p_{k-1}$ when $s$ belongs to $L^p_k$.
The result follows by standard bootstrapping.
\qed

The left composition lemma gives also the following
\begin{re}\label{lp2diff} The map $l$ extends  to a differentiable  map
$$l:\R\times L^p_{2}(\hat P \times_{\ad}i\kg_\beta^\bot)\ra L^p(\hat
P
\times_{\ad}i\kg_\beta^\bot)\ .
$$
\end{re}

\begin{lm}\label{invertible}  Suppose that $l(\varepsilon,s)=0$ with $s\in
A^0(\hat P
\times_{\ad}i\kg_\beta^\bot)$ and $\varepsilon > 0$.  Then the partial
derivative
$$\frac{\partial}{\partial s}l(e,s):
L^p_{2}(\hat P \times_{\ad}i\kg_\beta^\bot)\ra L^p(\hat P
\times_{\ad}i\kg_\beta^\bot) $$
 is an invertible operator.
\end{lm}
\pf

For $\dot s\in L^p_{2}(\hat P
\times_{\ad}i\kg_\beta^\bot)$, put
$$ u:=(d_{\frac{s}{2}}
\exp)(\frac{\dot s}{2})\ e^{-\frac{s}{2}} \in  \hat P
\times_{\ad}(i\kg_\beta^\bot\otimes\C) \ .
 $$
The components of $u$ in  $\hat P
\times_{\ad}i\kg_\beta^\bot$ and   $\hat P
\times_{\ad} \kg_\beta^\bot$ are respectively
$$u_h=\frac{1}{2}(u+\ad_{e^{-\frac{s}{2}}} u)\ ,\
u_a= \frac{1}{2}(u-\ad_{e^{-\frac{s}{2}}} u)\ .
$$

One has
$$(\frac{\partial}{\partial s}F_{A_s})(\dot s)=d_{A_s}(\frac{\partial}{\partial s}
A_s)=d_{A_s}\left(\partial_{A_s}(\ad_{e^{-\frac{s}{2}}} u)
-\bar\partial_{A_s} u \right)=
$$
$$-\partial_{A_s}\bar\partial_{A_s}
u+\bar\partial_{A_s}\partial_{A_s}(\ad_{e^{-\frac{s}{2}}}u)=
 [u_a,F_{A_s}]+(\bar\partial_{A_s}\partial_{A_s}-
\partial_{A_s}\bar\partial_{A_s})(u_h)\ .
$$
On the other hand, the same computation as in the finite dimensional framework (see
the proof of Proposition \ref{HilbCrit}) give
$$ \frac{\partial}{\partial
s}i\mu(e^{\frac{s}{2}}\varphi)=d(i\mu)(u^\#_{e^{\frac{s}{2}}\varphi})
=d(i\mu)([u^\#_h]_{e^{\frac{s}{2}}\varphi}+[u^\#_a]_{e^{\frac{s}{2}}\varphi})=
$$
$$=
d(i\mu)([u^\#_h]_{e^{\frac{s}{2}}\varphi})+[u_a,i\mu(e^{\frac{s}{2}}\varphi )]\ .
$$
Therefore
\begin{equation}\label{ddsl}\frac{\partial}{\partial s}l(\varepsilon,s)(\dot
s)=
$$
$$i\Lambda_g(\bar\partial_{A_s}\partial_{A_s}-
\partial_{A_s}\bar\partial_{A_s})(u_h)+d(i\mu)([u^\#_h]_{e^{\frac{s}{2}}\varphi})+
[u_a, l(\varepsilon,s)-\varepsilon s]+\varepsilon\dot s\ .
\end{equation}
Recall the following well-known identity which holds
for any unitary  integrable connection $B$ on a Hermitian vector bundle:
$$P(|v|^2):=i\Lambda_g\bar\partial \partial |v|^2= (i\Lambda_g\bar\partial_B
\partial_Bv,v)-(v,i\Lambda_g\partial_B\bar\partial_B v)-|d_B(v)|^2\ .
$$

In our  case, using the Euclidean inner product on  $\hat P\times_{\ad}
i\kg$,  this yields
\begin{equation}\label{pnorm}
P |u_h|^2= (i\Lambda_g(\bar\partial_{A_s}
\partial_{A_s}-
\partial_{A_s}\bar\partial_{A_s})u_h,u_h)   -|d_{A_s}(u_h)|^2\ .
\end{equation}

Therefore, for a pair $(\varepsilon,s)$ with $l(\varepsilon,s)=0$,  and using the
well known   formulae obtained in the finite dimensional case, we get
\begin{equation}
\left\langle \frac{\partial}{\partial s}  l(\varepsilon,s)(\dot
s),u_h\right\rangle_{L^2} =\|d_{A_s} u_h\|^2_{L^2}
+\|[u^\#_h]_{e^{\frac{s}{2}}\varphi}\|^2-
\varepsilon\langle [u_a,  s],u_h\rangle+ 2\varepsilon\langle \frac{\dot
s}{2},u_h\rangle $$
$$\geq\frac{1}{2}
\varepsilon \|\dot s\|_{L^2}^2\
\end{equation}
by Proposition \ref{Differential} in the Appendix.  This proves that $
\frac{\partial}{\partial s} l(\varepsilon,s)$ is an injective operator. It suffices to
notice that its index vanishes.
\qed

\begin{lm} There is a positive constant $C=C(A^0,\varphi)$ such that the
following estimate holds:
$$C\|v\|_{L^2}^2\leq \| \bar\partial_{A^0}
v\|^2_{L^2}+\|v^\#_\varphi\|^2_{L^2}\    \forall v\in A^0(\hat P\times_\ad
(i\kg_\beta^\bot\otimes\C)) \ .
$$
\end{lm}
\pf If not, there would exist a sequence $(v_n)_n$ with $\|v_n\|^2_{L^2}\equiv
1$ such that
$$\| \bar\partial_{A^0}
v_n\|^2_{L^2}+\|(v_n)^\#_\varphi\|^2_{L^2}\ra 0
$$
But this implies that $v_n$ has a subsequence which is weakly convergent in
$L^2_1$ and strongly convergent in $L^2$. The limit, say $v_0$,  satisfies
$\bar\partial_{A^0} (v_0)=0$, hence it is smooth. Noting that
$(\cdot)^\#_\varphi$ can be regarded as a smooth linear bundle  morphism
$$\hat P\times_\ad
(i\kg_\beta^\bot\otimes\C))\map \varphi^*(T^{\rm vert} (E))\ ,
$$
one gets
$$\bar\partial_{A^0} v_0=0\ ,\ (v_0)^\#_{\varphi}=0\ .
$$
Therefore, by Remark \ref{CplxStab},  one has
$v_0\in\g_\beta=\kg_\beta\otimes\C$. But this contradicts $v_0\in A^0(\hat P\times_\ad
(i\kg_\beta^\bot\otimes\C))$.
\qed

Note that   this lemma holds for any holomorphic pair $\beta$ and $\hat K$-reduction
$\hat P$ such that the relation
$\g_\beta=\kg_\beta\otimes\C$ between its infinitesimal stabilizers holds.
\\

For any $s\in A^0(\hat P\times_\ad
(i\kg_\beta^\bot))$ put
$$m_s:= \sup_X|s|\ .$$

\begin{co}\label{RevEst}
There exists a positive constant $C(A^0,\varphi,m_s)$ such that
for all $s\in A^0(\hat P\times_\ad
(i\kg_\beta^\bot))$ and $w\in A^0(\hat P\times_\ad
(i\kg_\beta^\bot\otimes\C))$ it holds
$$C(A^0,\varphi,m_s)\|w\|^2_{L^2}\leq
\|\bar\partial_{A_{s}} w\|^2_{L^2}
+\|[w]^\#_{e^{\frac{s}{2}}\varphi}\|^2\ .
$$
\end{co}
\pf  Defining $v:=\ad_{e^{-\frac{s}{2}}} w$ one has
$$\bar\partial_{A_{s}} w=\ad_{e^{\frac{s}{2}}}(\bar\partial_{A^0} v)\ ,\
[w]^\#_{e^{\frac{s}{2}}\varphi}=(e^{\frac{s}{2}})_{*,\varphi}(v^\#_\varphi)\
.
$$
But one has estimates of the form
$$\|\ad_{e^{\frac{s}{2}}}(\bar\partial_{A^0} v)\|^2_{L^2}\geq C'(m_s)
\|(\bar\partial_{A^0} v)\|^2_{L^2}\ ,\
\|(e^{\frac{s}{2}})_{*,\varphi}(v^\#_\varphi)\|^2_{L^2}\geq
C''(\varphi,m_s)\|(v^\#_\varphi)\|^2_{L^2}
$$
$$\|v\|_{L^2}^2\geq C'''(m_s) \|w\|_{L^2}^2\ .
$$
It suffices to apply the previous lemma.
\qed
\\

Let $(s_\varepsilon)_\varepsilon$ be a smooth family
in $A^0(\hat P\times_\ad i\kg_\beta^\bot)$  such that for all $\varepsilon$
it holds $l(\varepsilon,s_\varepsilon)=0$.  Put
$$u =(\frac{d}{d\varepsilon}
\exp(\frac{s }{2}))e^{-\frac{s }{2}}\ ,\ \dot s := \frac{d}{d\varepsilon} s
 \ .
$$

\begin{lm} One  has an estimate of the  form
\begin{equation}
\label{uhest}\sup_X | u_h|\leq C(m_s)\ .
\end{equation}
\end{lm}
\pf By (\ref{ddsl})  one gets
\begin{equation}\label{CrucId}
0=\frac{\partial}{\partial \varepsilon} l(\varepsilon,
s_\varepsilon)=
$$
$$i\Lambda_g(\bar\partial_{A_s}\partial_{A_s}-
\partial_{A_s}\bar\partial_{A_s})(u_h)+d(i\mu)([u^\#_h]_{e^{\frac{s}{2}}\varphi})-
\varepsilon[u_a,   s]+\varepsilon\dot s +s\ .
\end{equation}

Taking this time inner product with $u_h$ pointwise and using again
(\ref{pnorm}) and the inequalities in Proposition \ref{Differential}, one gets
\begin{equation}\label{pointwise}
0=P ( |u_h|^2)+|d_{A_s} u_h|^2+|
[u^\#_h]_{e^{\frac{s}{2}}\varphi}|^2 -\varepsilon ([u_a,s],u_h)+\varepsilon
(\dot s,u_h)
(\dot s,u_h)+(s,u_h)\geq
$$
$$
P ( |u_h|^2)+|d_{A_s} u_h|^2+|
[u^\#_h]_{e^{\frac{s}{2}}\varphi}|^2+(s,u_h)\ .
\end{equation}

Integrating over $X$   one gets
$$
\|d_{A_s}
u_h\|_{L^2}^2+\|[u^\#_h]_{e^{\frac{s}{2}}\varphi}\|_{L^2}^2\leq
\|s\|_{L^2}\|u_h\|_{L^2}\leq\sqrt{Vol_g(X)}\ m_s  \|u_h\|_{L^2}\ ,
$$
whereas the term on the left is larger than
$C(A^0,\varphi,m_s)\|u_h\|_{L^2}^2$ by Corollary  \ref{RevEst}.  This gives an
$L^2$-estimate  of the form
\begin{equation}\label{EstEst}\| u_h\|_{L^2}\leq C_1(A_0,\varphi, m_s)\ .
\end{equation}

Coming back to (\ref{pointwise}) we see that
$$P(|u_h|^2)\leq  -(s,u_h) \leq |s||u_h|\leq \frac{1}{2}(|u_h|^2+m_s^2)\ .
$$
Using now  Lemma 3.3.2 in \cite{LT} and (\ref{EstEst}), one gets an estimate
of the form
\begin{equation}
\sup_X(|u_h|^2)\leq  C(g)(\|u_h\|_{L^2}^2 +
\frac{1}{2} m_s^2)\leq C(g,A^0,\varphi,m_s)\ .
\end{equation}
\qed

Using again the formulae (\ref{As}) we get
\begin{equation}\label{Pform}\bar\partial_{A_s}\partial_{A_s}-\partial_{A_s}\bar\partial_{A_s}
=\bar\partial\partial_0-\partial_0\bar\partial
$$
$$
-\Psi(-s)(\bar\partial s)\wedge \partial_0-\Psi(-s)(\bar\partial s)\wedge
\Psi(s)(\partial_0 s )-\Psi(s)(\partial_0 s )\wedge\bar\partial +
\Psi(s)(\partial_0 s )\wedge\Psi(-s)(\bar\partial s) $$
$$+\Psi(s)(\bar\partial\partial_0 s)+\Psi(-s)(\partial_0\bar\partial s)+
\partial_0( \Psi(-s))\wedge  \bar\partial s +\bar\partial(
\Psi(s))\wedge  \partial_0 s\ .
\end{equation}

We   plug the identity (\ref{Pform}) in (\ref{CrucId}) and use standard elliptic $L^p$ -
 estimates for the operator
$i\Lambda_g(\bar\partial\partial_0-\partial_0\bar\partial )$.

Note first that, since $\Psi$ is of class ${\cal C}^1$, one has an estimate  of the
form
$$|d_{A^0}( \Psi(s))|\leq C(m_s) |d_{A^0}s| \ . $$
Recalling that, by  (\ref{uhest}),   $\sup_X| u_h|$   is   bounded
by a constant $C(m_s)$, we obtain
\begin{equation}\label{lpestuh}\| u_h\|_{L^p_2}\leq
c(m_s)(1+\|s\|_{L^{2p}_1}\|
u_h\|_{L^{2p}_1}+\|s\|_{L^{2p}_1}^2+\|s\|_{L^{p}_2})+
$$
$$+\|d(i\mu)([u^\#_h]_{e^{\frac{s}{2}}\varphi})-
\varepsilon[u_a,   s]+\varepsilon\dot s +s\|_{L^p}\ .
\end{equation}

On the other hand, by the formulae obtained in the proof of Proposition
\ref{Differential}, one can express $\dot s$ in terms of $u_h$ as
$$\dot s=2\theta \left(\left[\frac{s}{2},\cdot\right]\right)(u_h)\ ,
$$
where  this time $\theta (t)= \frac{2t}{e^t-e^{-t}}$. Since the bundle map $\hat
P\times_{\ad} i\kg \ra \Herm(\hat P\times_{\ad}\g)$ given by $a\mapsto 2
\theta\left(\left[\frac{a}{2},\cdot\right]\right)$ is of class ${\cal C}^2$, one
gets  estimates  of the form
$$ |\dot s|\leq c(m_s) |u_h|\ ,\ \| \dot s\|_{L^p_2}\leq c(m_s) \|
u_h\|_{L^p_2}  \  ,
$$
in particular, from (\ref{uhest}) we deduce
\begin{equation}\label{supdotsX}
\sup_X  |\dot s| \leq c(m_s)\ .
\end{equation}
Expressing $u_h$ and $u_a$ similarly in terms of $\dot s$, one gets
$$  |u_a|\leq c(m_s)
|\dot s|\ ,\ \| u_h\|_{L^{2p}_1}\leq c(m_s)\| \dot s\|_{L^{2p}_1}\ .
$$

Note finally that  (since $u_h$ and $s$ are bounded in terms of $m_s$)
$[u^\#_h]_{e^{\frac{s}{2}}\varphi}$ remains in a compact subset
$K(\varphi,m_s)$ of the vertical tangent bundle $T^{\rm vert}(E)$. Therefore,
one gets
$$\sup_X|d(i\mu)([u^\#_h]_{e^{\frac{s}{2}}\varphi})|\leq c(m_s)\ .
$$
Summarizing, we obtain an estimate  of the form
\begin{equation}\label{lpestdots}\| \dot s\|_{L^p_2}\leq
c(m_s)\left(1+\|s\|_{L^{2p}_1}\|
\dot s\|_{L^{2p}_1}+\|s\|_{L^{2p}_1}^2+\|s\|_{L^{p}_2}\right) \ .
\end{equation}

Using the general inequality
\begin{equation}\label{xlpestim}
\|x\|_{L_1^{2p}} \leq C\cdot\left(\sup_X|x|^\frac{1}{2}\right)\cdot
\|x\|^\frac{1}{2}_{L^p_2} + \|x\|_{L^{2p}}\ ,
\end{equation}
(\cite{Au} Theorem 3.69; the constant $C$ depends only on $p$ and the
dimension of $X$) and (\ref{supdotsX}) one proceeds as in the proof of \cite{LT}
Proposition 3.3.5 to obtain an inequality of the form
\begin{equation}\label{speed}
\|\dot s\|_{L_2^{p}} \leq c(m_s)(1+\|s\|_{L_2^{p}})\ .
\end{equation}
Integrating from $\varepsilon$ to 1 this yields (compare with \cite{LT}, Proposition 3.3.5)

\begin{pr}\label{lp2est} Let $(s_\varepsilon)_{\varepsilon\in (\varepsilon_0,1]}$ be a
smooth family   in $A^0(\hat P\times_\ad i\kg_\beta^\bot)$  such that for all
$\varepsilon$ it holds $l(\varepsilon,s_\varepsilon)=0$.  Then one has an
estimate of the form
\begin{equation}
\|s_\varepsilon\|_{L_2^{p}} \leq e^{c(m_{s_\varepsilon})(1-\varepsilon)}(1 +
\|s_1\|_{L_2^{p}})\ \forall \varepsilon\in (\varepsilon_0,1]\ .
\end{equation}
\end{pr}
\vspace{3mm}

We know that the equation  $l(\cdot, 1)=0$ has a smooth solution $\zeta$
(Proposition \ref{initialisation}). Consider the set
$${\cal S}:=\{\sg: (\varepsilon_\sg,1]\ra A^0(\hat P\times_\ad
i\kg_\beta^\bot) |\ 0\leq\varepsilon_\sg<1,\
l(\varepsilon,\sg(\varepsilon))\equiv 0\ ,\ \sg(1)=\zeta\}\ .
$$

${\cal S}$ is nonempty by Lemma \ref{invertible}, implicit function theorem and
Lemma \ref{regularity}. Writing
$\sg_1\leq  \sg_2$ when
$\varepsilon_{\sg_2}\leq
\varepsilon_{\sg_1}$ and $\sg_1=\sg_2|_{(\varepsilon_{\sg_1},1]}$, ${\cal
S}$ becomes  an inductively ordered set, hence there is a maximal element in
${\cal S}$.
\begin{pr} \hfill{\break}
1. Let $\sg:(\varepsilon_\sg,1]\ra  A^0(\hat P\times_\ad
i\kg_\beta^\bot)\in {\cal S}$ and put $s_\varepsilon:=\sg(\varepsilon)$. If
$$\sup\{\|s_\varepsilon\|_{L^2}| \ \varepsilon\in (\varepsilon_\sg,1]\}<
\infty\ ,
$$
then the strong limit $\lim_{\varepsilon\ra \varepsilon_\sg} s_\varepsilon$
exists in $L^p_2$,  and is a smooth solution of the equation
$l(\varepsilon_\sg,\cdot)=0$.
\vspace{2mm}
\\
2. For any maximal element $\sg\in {\cal S}$, one
has $\varepsilon_\sg=0$.
\end{pr}
\pf \\
1.  By assumption and the third statement  in Lemma \ref{uniform} one gets an
uniform bound
\begin{equation}\label{M}m_{s_\varepsilon}=\sup_X  |s_\varepsilon | \leq  M\ ,\
\forall\
\varepsilon\in (\varepsilon_\sg,1]\ .
\end{equation}

Therefore, the family $(s_\varepsilon)_{\varepsilon\in (\varepsilon_\sg,1]}$ is
bounded in $L^p_2$ by  Proposition \ref{lp2est}.   Using now the inequality
(\ref{speed}), we see that $\frac{d}{d\varepsilon} s_\varepsilon$ is also bounded in
$L^p_2$, hence the map $\sg:(\varepsilon_\sg,1]\ra L^p_2(\hat P\times_{\ad}
i\kg_\beta^\bot)$ is Lipschitz. Therefore  one has a strong limit
$$s=\lim_{\varepsilon\ra \varepsilon_\sg} s_\varepsilon\ ,
$$
which will be a solution of the equation $l(\varepsilon_\sg,\cdot)=0$ by Remark
\ref{lp2diff}, and will be smooth by   Lemma
\ref{regularity}.

2. If $\varepsilon_\sg>0$, the second statement of Lemma  \ref{uniform} gives an
uniform bound of the form (\ref{M}). Therefore, the first part of this proposition
applies and gives a   strong  $L^p_2$ - limit $s=\lim_{\varepsilon\ra \varepsilon_\sg}
s_\varepsilon$ which is smooth and solves the equation
$l(\varepsilon_\sg,\cdot)=0$.

Using Lemma \ref{invertible} and the implicit function theorem we get a smooth
extension of $\sg$ on a larger interval, contradicting maximality.
\qed

This gives immediately the following crucial result:
\begin{thry}\label{boundedcase} \hfill{\break}
1. There exists a smooth family $(s_\varepsilon)_{
\varepsilon\in(0,1]}$ with $l(\varepsilon,s_\varepsilon)\equiv 0$.\\
2. If  $\sup\{\|s_\varepsilon\|_{L^2}|\ \varepsilon\in (0,1]\}<
\infty$, then $s_\varepsilon$ converges strongly in $L^p_2$ to a smooth solution
$s$ of the equation (\ref{s-eq}) as $\varepsilon\ra 0$.
\end{thry}

Therefore, when $\sup\{\|s_\varepsilon\|_{L^2}|\ \varepsilon\in (0,1]\}<
\infty$, one gets  a Hermitian-Einstein reduction   $\hat
P_s:=e^{-\frac{s}{2}}(\hat P)\in {\cal R}_{P_0}(\hat Q)$.

\subsection{Destabilizing the pair in the unbounded case}

\subsubsection{Estimates   in the unbounded case }
\label{Estimates}

  Let $(s_\varepsilon)_{
\varepsilon\in(0,1]}$ be the smooth family given by Theorem \ref{boundedcase}.
Our first aim is to  get a  uniform
$L^2_1$ - bound for  the $L^2$-normed section
$$\sigma_\varepsilon:=
\frac{s_\varepsilon}{\nr s_\varepsilon\nr_{L^2}}\ .$$

 By the third statement  of Lemma \ref{uniform}, we
get    uniform bounds
\begin{equation}\label{SupSigma} \sup_X|
\sigma_\varepsilon|\leq M\ , \ \sup_{x\in X} \{|\lambda|\ |\ \lambda\in{\rm
Spec}([(\sigma_\varepsilon)_x,\cdot])\}\leq M\
\end{equation}
for all  $\varepsilon >0$ for which $\|s_\varepsilon\|_{L^2}\geq 1$.

\begin{re}\label{IneqExp}  Fix $\alpha\in (0,1)$. Then, for any
$u\in\R$, $n>0$ one has
\begin{enumerate}
\item $$  |e^{n u}-1|\geq \left\{
\begin{array}{ccc} \alpha n  |u|&\rm if&n u\geq \log(\alpha)\\
(1-\alpha)  &\rm if&n u\leq \log(\alpha)\ .
\end{array}
\right.
$$
\item  Supposing $|u|\leq M$, we get
$$ |e^{n u}-1|\geq
\left[\min\left(n\alpha,\frac{1-\alpha}{M}\right)\right]|u|\ ,\
\forall
u\in[-M,M]\ .
$$
\item Choosing   $\alpha=\frac{1}{2}$, we get
\begin{enumerate}
\item
$$|e^{n u}-1|\geq \frac{1}{2M}|u|\ ,\ \forall n\geq\frac{1}{M}\
,\ \forall\
u\in[-M,M]\ .
$$
\item \label{useful} $$\frac{1}{2M}\leq n\eta^2(nu) \ ,\ \ \forall n\geq\frac{1}{M}\
,\ \forall\
u\in[-M,M]\ .
$$
\item\label{c} Let $V$ be Hermitian space.  Then
$$\frac{1}{2M }\|x\|^2\leq n\|\eta( nh)(x)\|^2
$$
for all $n\geq \frac{1}{M}$ and for every Hermitian endomorphism $h\in\Herm(V)$
with the property
${\rm Spec}(h)\subset [-M,M]$.
\end{enumerate}
\end{enumerate}
\end{re}
\pf The first statements are straightforward. The last one  follows from the previous ones
and Remark \ref{monotony} in the following way: we take
$$f:=\sqrt{\frac{1}{2M}}\ \id_V\ ,\ g(t):=\sqrt{n}\ \eta(nt)\ ,
$$
and we notice that $|f|_{[-M,M]}\leq |g|_{[-M,M]}$ by  \ref{useful}.
\qed

  We divide (\ref{NewIneq}) by
$\|s\|_{L^2}$, taking into account      Remark
\ref{IneqExp} \ref{c}  and (\ref{SupSigma}). We get, for $\|s\|_{L^2}\geq\max(
\frac{1}{M},1)$,
$$\frac{1}{2M}|\partial_0
\sigma|^2\leq
\|s\| |\eta(\|s\|[\sigma,\cdot])(\partial_0\sigma)|^2=\frac{1}{\|s\|}|\eta(
[s,\cdot](\partial_0 s)|^2 \leq
$$
$$
\leq\frac{1}{\|s\| }\left[-\frac{1}{2}P(|s|^2)-(k^0,s)
\right]\ .
$$

 Integrating on
$X$ the obtained inequality, we get
\begin{equation}\label{BoundSigma}\|\partial_0 \sigma\|^2_{L^2}\leq 2M\|
k^0\|_{L^2}\ .
\end{equation}

Concluding, we can state
\begin{pr}\label{sigma_n}  Suppose that $(s_\varepsilon)_{\varepsilon\in
(0,1]}$ is not $L^2$-bounded. Then    the set
$$\left\{\sigma_\varepsilon\ |\  \|s_\varepsilon\|_{L^2}\geq\max\left(
\frac{1}{M},1\right)\right\}
$$
 is
uniformly bounded in
$L^\infty\cap L^2_1$.  In particular, there exists a  sequence $(\varepsilon_n)_n$,
$\varepsilon_n\ra 0$ such that
\begin{enumerate}
\item $\lim_{n\ra\infty} \|s_{\varepsilon_n}\|_{L^2}=\infty$  .
\item The sequence $(\sigma_{\varepsilon_n})_n$ converges
\begin{itemize}
\item weakly in $L^2_1$ ,
\item strongly in $L^2$ ,
\item almost everywhere on $X$
\end{itemize}
to an  $L^2_1\cap L^\infty$ - section $\sigma$ in the bundle $\hat
P^0\times_{\ad}i\kg_\beta^\bot$.
\end{enumerate}
 \end{pr}

\subsubsection{The properties of the   limit $\sigma$}

The aim of this paragraph is the following
\begin{thry}\label{PropSigma}
\hfill{\break}\vspace{-6mm}
\begin{enumerate}
\item $\|\sigma\|_{L^2}=1$.
\item The map $X\ni x\mapsto [\sigma(x)]\in i\kg/{\ad_{\hat K}}$ which assigns to
every $x\in X$ the conjugacy class of $\sigma(x)$, is constant almost everywhere.
\item ${\rm pr}_{V^+_{\sigma(x)}}(\bar\partial \sigma)=0$ for almost
every $x\in X$,
where $V^+_{\sigma(x)}$ is the direct sum of the eigenspaces
corresponding to strictly
positive eigenvalues of the endomorphism
$[\sigma(x),\cdot]$ on $\g_x:=\hat P_x\times_{\ad}\g$.
\end{enumerate}
\end{thry}

In order to prove this, we will need   several  preparations:

Let $s\in A^0(\hat P\times_{\ad}i\kg)$ be a smooth section. Denote by
$V^0_s$ the
linear fibration $\union_{x\in X} V^0_{s,x}$, where
$V^0_{s,x}=\ker[s_x,\cdot]\subset
\hat P_x\times_{\ad}\g$.

Let $X^s_k$ be the closed subset of $X$ where $\ker[s_x,\cdot]$ has rank $\leq k$.
Then $\ker([s,\cdot])$ defines a continuous subbundle of  $\hat
P\times_{\ad}\g|_{X^s_k\setminus X^s_{k-1}}$.   Therefore, the projection
$p_{V^0_s}$ is continuous on every $X^s_k\setminus X^s_{k-1}$, hence it is
a bounded measurable section.

  With these remarks, we can state

\begin{pr}\label{proj^0}  Let $\hat A^0\in{\cal A}(\hat P)$ and let
$(\sigma_n)_n$,
$\sigma_n\in A^0(\hat P\times_{\ad}i\kg)$  be a sequence of smooth
sections with the
following properties:
\begin{enumerate}
\item  It is   weakly convergent in $L^2_1$.
\item It is bounded in $L^\infty$.
\item It converges strongly in $L^2$.
\item It converges almost everywhere on $X$.
\item  $(\|{\rm pr}_{V^0_{\sigma_n}}[\partial_0\sigma_n]\|_{L^2})_n$
converges
to 0.
\end{enumerate}
Then the weak  limit $\sigma$  of this sequence (which belongs to  $L^2_1\cap
L^\infty$)   defines an almost everywhere  constant map $x\mapsto
[\sigma(x)]\in\kg/\hat K$.
\end{pr}
\pf
Let  Let $\iota_1,\dots,\iota_k\in\C[\hat \g]$ be
homogeneous
generators of the invariant algebra $ \C[\hat \g]^{ \hat G}$.  We will prove that
$\iota_i(\sigma)$ are constant almost everywhere, and the claim will follow from Corollary
\ref{Invariants} in the Appendix (compare with the proof of
\ref{ConstConjClassMap}).

Let $j_i:\hat\g^{d_i}\ra\C$ be the $\ad_{\hat G}$-invariant symmetric
multilinear map
which corresponds to $\iota_i$

Using properties 2. and 3. and 4. one gets easily, by the Lebesgue dominant
convergence theorem, that
$\iota_i(\sigma_n)=j_i(\sigma_n,\dots,\sigma_n)$ converges to $\iota_i(\sigma)$
strongly in $L^2$.

One has
$$\partial[j_i(\sigma_n,\dots,\sigma_n)]=
d_i j_i(\partial_0\sigma_n,\sigma_n,\dots, \sigma_n)\ .
$$

On the other hand, by 2. and 3. the sequence of sections
$$j_i(\cdot,\sigma_n,\dots,\sigma_n)\in
A^0(\hat P\times_{\ad}\Hom(i\kg,\C))$$
is $L^\infty$ - bounded and converges strongly in $L^2$  to
$j_i(\cdot,\sigma,\dots,\sigma)$. Taking into account  1., 2. and 3. we deduce by
Lemma \ref{analysis} in Appendix  that
$$d_i j_i(\partial_0\sigma_n,\sigma_n,\dots, \sigma_n)\ra d_i
j_i(\partial_0\sigma,\sigma,\dots ,\sigma) \ \hbox{weakly in}\ L^2\ .$$

 On the other hand
$d_i j_i(\partial_0\sigma,\sigma,\dots ,\sigma)=\partial [ j_i(\sigma,\dots,\sigma)]$ as
distributions because Leibniz's rule extends to
$L^2_1\cap L^\infty$-sections.\footnote{This can be proved by constructing  a
sequence of smooth sections  converging to $\sigma$ with respect to both
 $L^2_1$ - and
$L^p$ - norms, for $p$ sufficiently large.}

Therefore,
\begin{equation}\label{convergence}
d_i j_i(\partial_0\sigma_n,\sigma_n,\dots, \sigma_n)\ra \partial [ j_i(\sigma,\dots,\sigma)]\
  \ \hbox{weakly in}\ L^2\ .
\end{equation}

On the other hand, taking into account that $j_i$ is $\ad_{\hat
G}$-invariant (hence in
particular $\ad_G$-invariant), it follows easily that, for any $v\in
\g$, the linear
functional  $j_i(\cdot,v,\dots,v)$ vanishes on
$[v,\g]$.  The point is that, when $v\in i\kg$, one has an orthogonal
decomposition
$$\g=[v,\g]\oplus z_\g(v)\ .
$$
This follows from the fact that, in this case, $z_\g(v)$ is just the
complexification of
$z_{\kg}(iv)$.

Therefore one
gets pointwise
$$d_i j_i(\partial_0\sigma_n,\sigma_n,\dots,
\sigma_n) =d_i j_i({\rm
pr}_{V^0_{\sigma_n}}(\partial_0\sigma_n),\sigma_n,\dots,
\sigma_n) \ ,
$$
and the right term converges to 0 in $L^2$ as $n\ra\infty$ by 5. and
2. By (\ref{convergence}) we get that
$\partial [ j_i(\sigma,\dots,\sigma)]=0$ in the week sense, hence
$j_i(\sigma,\dots,\sigma)$ is a constant. The claim follows now from
Corollary \ref{Invariants}
\qed

For an element $u\in i\kg$ we denote by $V^-_u$ the direct sum of
eigenspaces of
$[u,\cdot]\in\Herm(\g)$ corresponding to strictly negative
eigenvalues.  Again, for a
  section $u\in A^0(\hat P\times_{\ad}i\kg)$, the projection ${\rm
pr}_{V^-_u}$ on the
fibration defined by $u$ is in $L^\infty(\Herm (\hat P\times_{\ad} \g))$.

\begin{pr}\label{proj^-} Let $\hat A^0\in{\cal A}(\hat P)$ and let
$(\sigma_n)_n$,
$\sigma_n\in A^0(\hat P\times_{\ad}i\kg)$  be a sequence of smooth
sections with the
following properties
\begin{enumerate}
\item  It is  weakly convergent in $L^2_1$.
\item It is bounded in $L^\infty$.
\item It converges strongly in $L^2$.
\item It converges almost everywhere on $X$.
\item The weak limit $\sigma$ of this sequence defines an (almost everywhere)
constant conjugacy class $[\sigma]\in i\kg/\hat K$.
\item  $(\|{\rm pr}_{V^-_{\sigma_n}}[\partial_0\sigma_n]\|_{L^2})_n$
converges
to 0.
\end{enumerate}
Then the weak   limit $\sigma$   (which belongs to $L^2_1\cap
L^\infty$)  satisfies
$${\rm pr}_{V^-_{\sigma,x}}[\partial_0\sigma]=0\ \hbox{ for almost every}\ x\in
X\ .
$$
\end{pr}
\pf  Let $\sigma_0\in i\kg$ be a representative of $[\sigma]$ and
$\lambda_1<\dots\lambda_k<0$ be the negative eigenvalues of
$[\sigma_0,\cdot]$. Then $[\sigma_x,\cdot]$ has the same eigenvalues with
the same multiplicities for almost every $x\in X$.

Indeed, one can easily see that for  $u\in i\kg$, the eigenvalues of the endomorphism
$[u,\cdot]\in\Herm(\g)$ depend  only  on the class $[u]\in i\kg/\hat K$. This follows
from the fact that, for $\hat k\in \hat K$, one has $[\ad_{\hat k}(s),\cdot]=\ad_{\hat k}
[s,  \ad_{\hat{k}^{-1}}(\cdot)]$.

Let $\eta\in (\lambda_k,0)$ and $\chi:\R\ra[0,1]$ a smooth
increasing function  which is 1 on $(-\infty,\eta]$ and 0 on $[0,\infty)$.

Put
$$q_n:= \chi([\sigma_n,\cdot])\ .
$$

Since $\sigma_n$ converges strongly in $L^2$, it follows that   $q_n$
converges strongly in $L^2$ to $q:=\chi([\sigma,\cdot])={\rm
pr}_{V^-_{\sigma}}$.   This follows from  the Lebesgue dominated
convergence theorem, from 4. and the fact that, for any Hermitian vector space $E$,
  the map
$$\chi:\Herm(E)\map\Herm(E)\ ,\  f\mapsto\chi(f)
$$
is continuous.

Therefore $q_n[\partial_0\sigma_n]$ converges  to
$q(\partial_0(\sigma))={\rm pr}_{V^-_{\sigma}}[\partial_0\sigma]$ as
distributions. On the other hand
$$\nr q_n[\partial_0\sigma_n]\nr_{L^2}\leq \| {\rm
pr}_{V^-_{\sigma_n}}[\partial_0\sigma_n]\|_{L^2}\ ,
$$
which converges to 0 by 6.

Therefore ${\rm
pr}_{V^-_{\sigma}}[\partial_0\sigma]=0$ as a distribution, hence as an $L^2$
- section as well.
\qed

We can now give the \\ \\
{\bf Proof} (of Theorem \ref{PropSigma}):

We make use of Proposition \ref{sigma_n} and   put $s_n:=s_{\varepsilon_n}$.
$\sigma_n:=\sigma_{\varepsilon_n}$.  Since $\|\sigma_n\|_{L^2}=1$, and
$\sigma_n\ra\sigma$ strongly in $L^2$, we get $\|\sigma\|_{L^2}=1$ as claimed.

For the proof of 2. and 3. we claim first that
\begin{equation}\label{claim} (\|{\rm
pr}_{V_-^{\sigma_n}}[\partial_0\sigma_n]\|_{L^2})_n
\ra 0\ ,
\end{equation}
where , for $s\in A^0(\hat P\times_{\ad}i\kg)$,  we denoted by $V_-^{s}$
the fibration
whose fibre in $x\in X$ is the direct  sum of all eigenspaces
corresponding to all non-positive
eigenvalues of $[s_x,\cdot]$.

This follows easily in the following way: using the inequality
$$\frac{1-e^{-t}}{t}\geq 1 \ {\rm for}\ t\leq 0\ ,
$$
we see that $\eta\geq\chi_{-}$, where $\chi_-:=\chi_{\R_{\leq 0}}$ is the
characteristic function of the set $\R_{\leq 0}$.  By Remark \ref{monotony} and
(\ref{NewIneq}) one gets
$$\left|{\rm
pr}_{V_-^{\sigma_n}}[\partial_0\sigma_n]\right|^2_x=
\frac{1}{\|s_{n}\|^2}\left|{\rm
pr}_{V_-^{s_n}}[\partial_0
s_n]\right|^2_x=\frac{1}{\|s_{n}\|^2}\left|(\chi_-([s_n,\cdot])(\partial_0 s_n)\right|^2_x\leq
$$
$$\leq \frac{1}{\|s_{n}\|^2}\left|(\eta([s_n,\cdot])(\partial_0 s_n)\right|^2_x
  \leq \frac{1}{\|s_{n}\|^2}
\left[-\frac{1}{2} P(|s_n|^2) -\varepsilon|s_n|^2-(k^0(x),s_n)\right](x)\ .
$$
 It suffices to integrate  the obtained inequality over $X$, and let
$\varepsilon\ra 0$ taking into account that $\|s_n\|\ra\infty$.

In particular, it follows that
$$(\|{\rm
pr}_{V_{\sigma_n}^0}[\partial_0\sigma_n]\|_{L^2})_n
\ra 0\ ,$$
hence, by Propositions \ref{proj^0}, the limit $\sigma$ defines an almost everywhere
constant map $X\ra \i\kg/{\hat K}$.  Therefore, we can apply Proposition
\ref{proj^-}, because (\ref{claim}) obviously  implies  $ \|{\rm
pr}_{V^-_{\sigma_n}}[\partial_0\sigma_n]\|_{L^2}\ra 0$ .
\qed

\subsubsection{The limit $\sigma$ destabilizes}

Let $\sigma$ be the limit given by Proposition $\ref{sigma_n}$.

\begin{pr}\label{rho-s} Let $[\sigma]\subset i\kg$ be the element in
$i\kg/\hat K$ defined by the (almost everywhere constant) map $x\mapsto
[\sigma_x]$ (modulo
$\ad_{\hat K}$) and let $\sigma_0\in i\kg_\beta^\bot\setminus\{0\}$ be any
representative   of
$[\sigma]$.
 Then:
\begin{enumerate}
\item The $L^2_1\cap L^\infty$ section  $\sigma$ is smooth and takes
values in  $[\sigma]$ on the complement
$X_\rho$ of a Zariski subset of  codimension at least 2 such that the subspace
$$ \{q\in \hat {\cal Q}|_{X_\rho}, \ \sigma(q)\in \hat \g_{\sigma_0}\}
$$
extends to a meromorphic reduction $\hat {\cal
Q}^\rho\hookrightarrow\hat{\cal Q}$
of
$\hat {\cal Q}$ to
$\hat G_{\sigma_0}$ which is holomorphic over $X_\rho$.
\item The section  $s(\sigma_0,\rho)$ associated with $\sigma_0$ and the  $Z_{\hat
K}(\sigma_0)$-reduction
$\hat{\cal Q}^\rho\cap \hat P$ of $\hat P|_{X_\rho}$ coincides with $\sigma$. More
precisely,
$\sigma(q)=\sigma_0$ for all $q\in\hat {\cal Q}^\rho\cap \hat P$.
\end{enumerate}
\end{pr}
\pf This follows from  Proposition \ref{PropSigma} and Proposition \ref{AlmostHol}.
There is one simple detail to be explained.   The   weak holomorphy condition in Proposition
\ref{AlmostHol} reads in our case
$$ {\rm pr}_{V^+_{\hat\sigma (x)}}(\bar\partial_x \hat \sigma)=0\ \hbox{for almost
every}\ x\in X
$$
where $\hat\sigma$ is the section in $\hat P\times_{\ad} i\hat \g=i\ad(\hat P)$ defined
by
$\sigma$, and $V^+_{\hat\sigma (x)}$ is the direct sum of strictly positive eigenspaces of
$[\hat \sigma(x),\cdot]\in\Herm(\hat\g_x)$.

But, since
$$\bar\partial \hat\sigma\in L^2(\Lambda^{01}(\hat P\times_{\ad}\g))\
,\ V^+_{\hat\sigma (x)}=V^+_{\sigma (x)}\subset  \g_x\ ,
$$
this condition is equivalent to
$$ {\rm pr}_{V^+_{\sigma (x)}}(\bar\partial_x  \sigma)=0\ \hbox{for almost
every}\ x\in X\ .
$$
\qed

The following  result will end the proof of  Theorem \ref{purpose}.

\begin{thry}\label{destab} The pair $(\rho,\frac{\sigma_0}{2})$ destabilizes the
pair
$(\hat J,\varphi)$.
\end{thry}

The first step is to give an explicit formula for the total maximal weight

$$\frac{2\pi}{(n-1)!}\deg_g\left(\rho,
h(\frac{\sigma_0}{2})\right)+\int\limits_{X_{0}}
\lambda^{\frac{\sigma_0}{2}} (\varphi,\rho)vol_g\ .
$$

For the first term, we obtain  by (\ref{degree-new})
$$\frac{2\pi}{(n-1)!}\deg_g\left(\rho,  h(\frac{\sigma_0}{2}
)\right)=\frac{1}{2(n-1)!}\int_X h\left( \left[i F_{\hat
A^0} - \frac{i}{2}[a\wedge a]\right],\sigma \right)
\wedge\omega_g^{n-1}=$$
$$\frac{1}{2}\int_X h\left(\Lambda_g \left[i F_{\hat
A^0} - \frac{i}{2}[a\wedge a]\right],\sigma \right)
\wedge vol_g\ ,
$$
where $a$ is the second fundamental form of the $\hat K_{\sigma_0} =\hat
G(\sigma_0)\cap\hat K$-bundle $\hat {\cal Q}^\rho\cap\hat P $ in $\hat
P |_{X_\rho}$.

By (\ref{a-wedge-a}) we have
$${\rm pr}_{\hat \kg_{\sigma_0}} [a\wedge
a]=
-2\sum_{\lambda\in {\rm Spec}^+
 {[\sigma_0,\cdot] }} [a_\lambda \wedge\overline{a}_{\lambda} ]\ ,
$$
so using the known formula   (\ref{prod-wedge}), we get
\begin{equation}\label{long}
h(\frac{i}{2}\Lambda_g{\rm pr}_{\hat \kg_{\sigma_0}} [a\wedge
a], \sigma)=-\sum_{\lambda\in {\rm Spec}^+
 {[\sigma_0,\cdot] }}i\Lambda_gh\left([a _\lambda\wedge \bar a _\lambda],
\sigma
\right) =$$
$$=
-\sum_{\lambda\in {\rm Spec}^+
 {[\sigma_0,\cdot] }}i\lambda  \Lambda_g h(a_\lambda\wedge\bar
a_\lambda\rangle= -\sum_{\lambda\in {\rm Spec}^+
 {[\sigma_0,\cdot] }}\lambda| a_\lambda|^2 \ .
\end{equation}

The
(1,0)-form
$\partial_0\sigma$
 has non-trivial components only on the eigenspaces of
$[\sigma_0,\cdot]$ associated with strictly positive eigenvalues. Indeed, since the
conjugacy class of $\sigma$ is constant, $\partial_0\sigma$ has no
$\ker{[\sigma_0,\cdot]}$ - component by Proposition
\ref{ConstConjClassMap}. On the other hand, by the holomorphy   of $\hat{\cal
Q}^\rho|_{X_\rho}$  it follows that $a^{10}_\lambda=0$  for $\lambda\in {\rm Spec}^-
 {[\sigma_0,\cdot] } $  (see formula (\ref{Hol}), Proposition \ref{HolTangSp}).

Therefore
\begin{equation}\label{d0sigma}\partial_0 \sigma = \sum_{\lambda\in {\rm Spec}^+
 {[\sigma_0,\cdot] }} [a _\lambda,\sigma ]=-\sum_{\lambda\in {\rm Spec}^+
 {[\sigma_0,\cdot] }} \lambda  a _\lambda \ .
\end{equation}
Let $\phi$ be the real function defined by
$$\phi(\tau)=\left\{
\begin{array}{ccc}
\frac{1}{\sqrt{\tau}}&\rm if&\tau>0\\
0&\rm if&\tau\leq 0\ .
\end{array}\right.
$$
We get
$$\phi([\sigma,\cdot])(\partial_0 \sigma  )=-\sum_{\lambda\in {\rm Spec}^+
 {[\sigma_0,\cdot] }} \sqrt{\lambda}\   a _\lambda \ ,
$$
so, comparing with (\ref{long}), one can write
$$h(\frac{i}{2}\Lambda_g{\rm pr}_{\hat \kg_{\sigma_0}} [a\wedge
a], \sigma )=-\left|\phi([\sigma_0,\cdot])(\partial_0 \sigma  )\right|^2 \ .
$$
Since the (1,0)-form
$\partial_0\sigma$ has non-trivial components only on the eigenspaces of
$[\sigma,\cdot]$ associated with strictly positive eigenvalues,   the right hand term is
indeed smooth on $X_\rho$.

Therefore, we obtain for the total maximal weight  of the pair $(\hat J,\varphi)$ the
formula
$$\lambda((\hat
J,\varphi),(\rho,\frac{\sigma_0}{2})):=\frac{2\pi}{(n-1)!}\deg_g(\rho,
h(\frac{\sigma_0}{2}))+\int\limits_{X_{\rho}}
\lambda^{\frac{\sigma_0}{2}} (\varphi,\rho)vol_g=$$
$$=\frac{1}{2}\left[\langle
 i\Lambda_g F_{\hat
A^0},\sigma\rangle_{L^2}+\|\phi([\sigma,\cdot])(\partial_0 \sigma )\|_{L^2}^2
\right] + \int\limits_{X_\rho}
\lim_{t\ra\infty}\mu^{\frac{-i\sigma}{2}}(e^{t\frac{\sigma}{2}}\varphi)\ .
$$
\begin{lm}\label{limit} Let $(\rho,\sigma_0)$ be any pair consisting of an element
$\sigma_0\in i\kg$ and a meromorphic reduction of $\hat {\cal Q}$ to $\hat
G(\sigma_0)$.  Fix a
$\hat K$-reduction
$\hat P$ of $\hat {\cal Q}$, let $\sigma$ be the section in $A^0(\hat P\times_{\ad}i\kg)$
which corresponds to $\sigma_0$ and this reduction.
Let $\hat A^0\in{\cal A}(\hat P)$ be the Chern connection of $\hat P$ and
$d_{\hat A^0}=\bar\partial+\partial_0$ the associated covariant derivative. Then
$$\lambda((\hat J,\varphi),(\rho,\frac{\sigma_0}{2}))=$$
$$\frac{1}{2}\langle
 i\Lambda_g F_{\hat
A^0},\sigma\rangle_{L^2}+\lim_{t\ra\infty}\int\limits_{X_\rho}\left[
\frac{1}{2}t |\eta([t\sigma,\cdot])(\partial_0 \sigma )|^2 +
\mu ^{-i\frac{\sigma}{2}}(e^{t\frac{\sigma}{2}}\varphi)\right]vol_g\ .
$$
\end{lm}
\pf It suffices to prove that
$$\|\phi([\sigma,\cdot])(\partial_0(\sigma))\|_{L^2}^2=\lim_{t\ra\infty}
\int\limits_{X_\rho}  t\left|\eta([t\sigma,\cdot])(\partial_0 \sigma )\right|^2
vol_g\ .
$$

For $t>0$, consider the real analytic function $\eta_t$ defined by
\begin{equation}\label{eta_t}\eta_t(\tau)=\sqrt{t}\ \eta(t\tau)=\left\{
\begin{array}{ccc}
\sqrt{\frac{1-e^{-t\tau}}{\tau}}&\rm if& \tau\ne 0\\
\sqrt{t}&\rm if&\tau=0\ .
\end{array}\right.
\end{equation}

Recall that the conjugacy class of $\sigma$ is constant, so the eigenvalues of
$[\sigma,\cdot]$ are also constant (see the proof of Proposition \ref{proj^-}). Let
$\lambda_0$ be the first positive eigenvalue of $[\sigma,\cdot]$.  Note that, by
Proposition \ref{ConstConjClassMap} and   Theorem \ref{HolTangSp}, the form
$\partial_0(\sigma)|_{X_\rho}$ has non-trivial components only along the \ub{strictly}
positive eigenspaces  of
$[\sigma,\cdot]$.  Therefore, pointwise on $X_\rho$ one has
$$
 \left|
t\left|\eta([t\sigma,\cdot])(\partial_0 \sigma )\right|^2-
 |\phi([\sigma,\cdot])(\partial_0 \sigma ) |^2 \right| \leq
 \sup_{[\lambda_0,\infty]}|\eta_t^2-\phi^2||\partial_0\sigma|^2\ .
$$

The point is now that
 $\eta_t^2$ converges uniformly to $\phi^2$ on the interval
$[\lambda_0,\infty)$ as $t\ra\infty$. Integrating both sides on $X_\rho$ we get
the result.
\qed

Let $\sigma_n=\frac{1}{\|s_n\|_{L^2}}s_n$ be the sequence  with limit $\sigma$
given by   Proposition \ref{sigma_n}.

\begin{lm}\label{liminf}  Fix $t_0>0$. Then
$$\langle i\Lambda_g F_{\hat
A^0},\sigma\rangle_{L^2}+\int\limits_{X_\rho}\left[
\frac{1}{2}t_0 |\eta([t_0\sigma,\cdot])(\partial_0 \sigma )|^2 +
\mu ^{-i\frac{\sigma}{2}}(e^{t_0\frac{\sigma}{2}}\varphi)\right]vol_g\leq$$
$$
\liminf_{n\ra\infty}\left\{\langle i\Lambda_g F_{\hat
A^0},\sigma_n\rangle_{L^2}+\int\limits_{X_\rho}\left[
\frac{1}{2}t_0 |\eta([t_0\sigma_n,\cdot])(\partial_0 \sigma_n )|^2 +
\mu ^{-i\frac{\sigma_n}{2}}(e^{t_0\frac{\sigma_n}{2}}\varphi)\right]vol_g\right\}.
$$
\end{lm}
\pf    First of all notice that
\begin{equation}\label{ObvLim}
\langle i\Lambda_g F_{\hat
A^0},\sigma_n\rangle_{L^2}\ra \langle i\Lambda_g F_{\hat
A^0},\sigma\rangle_{L^2}\ ,
\end{equation}
because $(\sigma_n)_n$ converges strongly in $L^2$ to $\sigma$.

Second, one has
$$\mu ^{-i\frac{\sigma}{2}}(e^{t_0\frac{\sigma}{2}}\varphi)=
\langle \mu   ( \varphi),  -i\frac{\sigma}{2}\rangle+
[\mu ^{-i\frac{\sigma}{2}}(e^{t_0\frac{\sigma}{2}}\varphi)-
\mu ^{-i\frac{\sigma}{2}}( \varphi)]\ ,$$
$$\mu ^{-i\frac{\sigma_n}{2}}(e^{t_0\frac{\sigma_n}{2}}\varphi)=
\langle \mu   ( \varphi),  -i\frac{\sigma_n}{2}\rangle+
[\mu ^{-i\frac{\sigma_n}{2}}(e^{t_0\frac{\sigma_n}{2}}\varphi)-
\mu ^{-i\frac{\sigma_n}{2}}( \varphi)]\ .$$

Since  $\sigma_n\ra \sigma$ strongly in $L^2$, it follows that
\begin{equation}\label{FirstLim}\int\limits_{X_\rho}\langle \mu  ( \varphi),
-i\frac{\sigma_n}{2}\rangle vol_g\ra
\int\limits_{X_\rho}\langle
\mu  (
\varphi),  -i\frac{\sigma}{2}\rangle vol_g\ .
\end{equation}

The two terms $[\mu
^{-i\frac{\sigma_n}{2}}(e^{t_0\frac{\sigma_n}{2}}\varphi)-
\mu ^{-i\frac{\sigma_n}{2}}( \varphi)]$,
$[\mu ^{-i\frac{\sigma}{2}}(e^{t_0\frac{\sigma}{2}}\varphi)-
\mu ^{-i\frac{\sigma}{2}}( \varphi)]$ are {\it
pointwise} non-negative, because  the value in every point $x\in X$ can be identified with
the energy of a curve in the fibre $F_x$.   On  the other hand, the first term converges
almost everywhere to the second as $n\ra\infty$.  By   Fatou's lemma, it  follows
that
\begin{equation}\label{SecondLim}
\int\limits_{X_\rho}[\mu ^{-i\frac{\sigma}{2}}(e^{t_0\frac{\sigma}{2}}\varphi)-
\mu ^{-i\frac{\sigma}{2}}( \varphi)]vol_g\leq
$$
$$
\leq\liminf_{n\ra\infty}
\int\limits_{X_\rho}[\mu^{-i\frac{\sigma_n}{2}}(e^{t_0\frac{\sigma_n}{2}}\varphi)-
\mu ^{-i\frac{\sigma_n}{2}}( \varphi)]\ .
\end{equation}

Finally we state that
\begin{equation}\label{ThirdLim}
\int\limits_{X_\rho}\left[
  |\eta([t_0\sigma,\cdot])(\partial_0 \sigma )|^2  \right]vol_g\leq
\liminf_{t\ra\infty}\int\limits_{X_\rho}\left[
  |\eta([t_0\sigma_n,\cdot])(\partial_0 \sigma_n )|^2  \right]vol_g\ .
\end{equation}

Indeed, since $(\eta([t_0\sigma_n,\cdot]))_n$ converges strongly in $L^2$ to
$\eta([t_0\sigma_n,\cdot])$ and is $L^\infty$-bounded,  and
since $\partial_0(\sigma)$ converges weakly in $L^2$ to $\partial_0\sigma$, it follows
that $\eta([t_0\sigma_n,\cdot])(\partial_0 \sigma_n )$ converges weakly to
$\eta([t_0\sigma,\cdot])(\partial_0 \sigma )$ in $L^2$ (see    Lemma
\ref{analysis} in  the  Appendix).

Now the inequality (\ref{ThirdLim}) follows from the well known semicontinuity property
of the norm  with respect to the weak convergence in a Hilbert space.
The claim of our lemma follows from  (\ref{ObvLim}), (\ref{FirstLim}), (\ref{SecondLim})
and  (\ref{ThirdLim}).
\qed

\begin{lm}\label{MonotonyLemma}  Let $\sg\in A^0(\hat P\times_{\ad}i\kg)$.  Then
the maps
$$t_0\mapsto   t_0 |\eta([t_0\sg,\cdot])(\partial_0 \sg )|^2 \ , $$
$$t_0\mapsto
\mu_0^{-i\frac{\sg}{2}}(e^{t_0\frac{\sg}{2}}\varphi) $$
are pointwise increasing.
\end{lm}
\pf The second map is obviously pointwise increasing (see Remark \ref{increasing}).
For the first, it suffices to notice that the function $\eta_{t_0}$ defined by
(\ref{eta_t}) takes positive values and is increasing with respect to $t_0$. The
result follows from Remark
\ref{monotony}.
We can now give the\\ \\
{\bf Proof} (of Theorem \ref{destab}):

  If not, one would have
$$\lambda((\hat J,\varphi),(\rho,\frac{\sigma_0}{2}))>0\ .
$$
By  Lemma \ref{limit}, it would follow  that there exists $t_0>0$ such that
$$U:=\frac{1}{2}\langle
 i\Lambda_g F_{\hat
A^0},\sigma\rangle_{L^2}+\ \int\limits_{X_\rho}\left[
\frac{1}{2}t_0 |\eta([t_0\sigma,\cdot])(\partial_0 \sigma )|^2 +
\mu ^{-i\frac{\sigma}{2}}(e^{t_0\frac{\sigma}{2}}\varphi)\right]vol_g>0\ .
$$

By Lemma \ref{liminf} we get that for \ub{all} sufficiently large $n\in\N$, we have
$$\frac{1}{2}\langle
 i\Lambda_g F_{\hat
A^0},\sigma_n\rangle_{L^2}+\int\limits_{X}\left[
\frac{1}{2}t_0 |\eta([t_0\sigma_n,\cdot])(\partial_0 \sigma_n )|^2 +
\mu^{-i\frac{\sigma_n}{2}}(e^{t_0\frac{\sigma_n}{2}}\varphi)\right]vol_g>
\frac{U}{2}\ .$$
Choose $n$ sufficiently large such that $\|s_n\|_{L^2}> t_0$.  By the monotony Lemma
\ref{MonotonyLemma}, we get that
$$\frac{1}{2}\langle
 i\Lambda_g F_{\hat
A^0},\sigma_n\rangle_{L^2}+ $$ $$+\int\limits_{X}\left[
\frac{1}{2}\|s_n\| |\eta([\|s_n\|\sigma_n,\cdot])(\partial_0 \sigma_n )|^2 +
\mu^{-i\frac{\sigma_n}{2}}(e^{\|s_n\|\frac{\sigma_n}{2}}\varphi)\right]vol_g>
\frac{U}{2}\ ,$$
or, equivalently,
$$\frac{1}{\|s_n\|}\left\{\frac{1}{2}\langle
 i\Lambda_g F_{\hat
A^0},s_n\rangle_{L^2}+\int\limits_{X}\left[
\frac{1}{2}  |\eta([s_n,\cdot])(\partial_0 s_n )|^2 +
\mu^{-i\frac{s_n}{2}}(e^{ \frac{s_n}{2}}\varphi)\right]vol_g\right\}>
\frac{U}{2}  .$$

But, integrating (\ref{(eq.s)}) on $X$ and taking into  account  the identity given by
Proposition \ref{identity}, we get
$$\frac{1}{\|s_n\|}\left\{\frac{1}{2}\langle
 i\Lambda_g F_{\hat
A^0},s_n\rangle_{L^2}+\int\limits_{X}\left[
\frac{1}{2}  |\eta([s_n,\cdot])(\partial_0 s_n )|^2 +
\mu^{-i\frac{s_n}{2}}(e^{ \frac{s_n}{2}}\varphi)\right]vol_g\right\}=$$
$$=
-\frac{\varepsilon\|s_n\|}{2}\ .
$$
\qed

We now can finally give the proof of our main result.
\\ \\
{\bf Proof} (of Theorem \ref{purpose}):   Theorem \ref{boundedcase} yields a
Hermitian-Einstein reduction $\hat P_s$, unless the smooth family
$(s_\varepsilon)_{\varepsilon \in (0,1]}$ given by this theorem is not
$L^2$-bounded. In this case   Proposition \ref{sigma_n} applies and one obtains a
pair $(\rho,\frac{\sigma_0}{2})$ with
$$\lambda((\hat
J,\varphi),(\rho,\frac{\sigma_0}{2}))\leq 0$$ by Theorem
\ref{destab}.   But, using Remark \ref{centralcase} and the fact
that $\sigma_0\in i\kg_\beta^\bot\setminus\{0\}$, this contradicts
polystability.
\qed

\section{Examples and Applications}\label{App}

Our main theorem generalize previous results on ``universal"
Kobayashi-Hitchin correspondences in two directions
\begin{enumerate}
\item It holds when the base manifold is an arbitrary Hermitian  manifold,
\item It holds for {\it oriented} pairs, giving the appropriate {\it   complex  geometric}
stability and polystability conditions for this class of objects.
\end{enumerate}

In this chapter we will illustrate in concrete situations the importance of
both directions.  Using the correspondence between polystable oriented pairs
and  solutions of the generalized Hermitian-Einstein equations, we construct
a canonical isomorphism  between the corresponding moduli spaces, the {\it
Kobayashi-Hitchin isomorphism}. Using this isomorphism and ideas
from Donaldson theory, we endow the  moduli spaces with canonical Hermitian
metrics. We discuss in more detail moduli spaces of oriented connections,
Douady Quot-spaces and moduli spaces of non-abelian monopoles in the
non-K\"aherian framework.

\subsection{Oriented holomorphic principal bundles and oriented connections}
\subsubsection{The Kobayashi-Hitchin correspondence for oriented holomorphic principal
bundles}

Let $\hat Q$ be a principal $\hat G$ bundle over  a compact complex manifold and ${\cal
Q}_0$ a fixed bundle holomorphic structure of the $G$-quotient $Q_0$. Our moduli problem
asks:\\

{\it Classify bundle holomorphic structures $\hat {\cal Q}$ on $\hat Q$ which induce ${\cal Q}_0$
on $Q_0$, modulo the  gauge group ${\cal G}:=\Aut_{Q_0}(\hat Q)=\Gamma(\hat Q\times_{\Ad}
G)$.}

Such a holomorphic structure will be naturally called a ${\cal Q}_0$-{\it oriented}
holomorphic structure on $\hat Q$
\\ \\
{\bf Example:}  Choosing $\hat G=GL(r,\C)$, $G=SL(r)$ our problem becomes: Let $E$
be a differentiable vector bundle vector bundle over
$X$ and ${\cal L}$ a holomorphic structure on its determinant line bundle $L$. Classify all
holomorphic structures ${\cal E}$ on $E$ such that $\det({\cal E})={\cal L}$   modulo
the gauge group
${\cal G}:=\Gamma(X,SL(E))$.

Note first that the ${\cal G}$-stabilizer of a holomorphic structure
$\hat {\cal Q}$ always contains  the relative centralizer   $Z_G(\hat
G)=\{g\in G|\Ad_{\gamma}(g)=g\ \forall\gamma\in\hat G\}$, so  it
is natural to relax the {\it simplicity condition} for an oriented holomorphic structure
$\hat{\cal Q}$ by asking  its infinitesimal stabilizer to be  minimal, i.e.
$\g_{\hat {\cal Q}}=z(\g,\hat G)$, (see Definition \ref{relcent}).

Our bundle classification problem can be formally regarded as a coupled problem with
fibre $F=\{*\}$.  Let $T$ be the unique maximal compact subgroup of the
abelian reductive complex group $Z_G(\hat G)$   and $\tg:={\rm Lie}(T)$. Let   $h$ be an
$\ad$-invariant complex inner product of Euclidean type on $\hat \g$. The choice of a
Hamiltonian system  $(K\subset\hat K,g_F,\mu)$ with a $\hat K$-equivariant moment
map $\mu$ reduces to just the choice of an element  $it\in  \tg$.  Therefore, our
stability condition becomes\\

{\it  A ${\cal Q}_0$-oriented holomorphic structure   $\hat {\cal Q}$ on $\hat
Q$ is $t$-stable if its infinitesimal stabilizer $\g_{\hat {\cal Q}}$ coincides with
$z(\g,\hat G)$ and the following holds: for every vector $\xi\in H(G)$ and
meromorphic
$\hat G(\xi)$-reduction $\rho$ of  $\hat {\cal Q}$, one has
\begin{equation}\label{sss}
\frac{2\pi}{(n-1)!}\deg(\rho,h(\xi))  \geq  h(t,\xi) Vol_g(X)
\end{equation}
with strict inequality when $\xi\not\in z(\g,\hat G)$.
}

A vector $\xi\in H(G)$ decomposes as $\xi=\xi_0+ \zeta$, where $\xi_0\in
z(\g,\hat G)^\bot$ and $\zeta\in z(\g,\hat G)$. The  linear form $h(\zeta)$  is
$\ad_{\hat G}$-invariant, so $\deg(\rho,h(\zeta))$ can be identified with  a
{\it holomorphic   invariant}   $\deg(\hat {\cal Q},h(\zeta))$ of $\hat {\cal Q}$ (see
Definition \ref{defdeg}).    Since
$h$ is non-degenerate, there exists a unique element $\tau(\hat {\cal Q})\in i\tg $
such that
$$\deg(\hat {\cal Q},h(\zeta))=h(\tau(\hat {\cal Q}),\zeta) \  \forall \zeta\in
i\tg=H(Z_G(\hat G))\ .
$$
The holomorphic invariant $\tau$ plays the role of the degree in the
classical case of (non-oriented) bundles; it is obviously a topological invariant when $g$
is K\"ahlerian.  Since the map, $\zeta\mapsto \deg(\hat {\cal Q},h(\zeta))$ is linear on
$i\tg$, the inequality (\ref{sss}) implies
$$t=\frac{2\pi \tau(\hat {\cal Q})}{(n-1)!Vol_g(X)}\ .
$$
which means that {\it an oriented holomorphic structure $\hat {\cal Q}$ can be
$t$-stable only for a single constant $t$, which   is a holomorphic invariant of $\hat{\cal
Q}$} (and a topological invariant of $\hat Q$ when $g$ is K\"ahler).  This
motivates the following
\begin{dt} A ${\cal Q}_0$-oriented holomorphic structure   $\hat {\cal Q}$ on $\hat
Q$ is called:
\begin{enumerate}
\item semistable, if  $\deg(\rho,h(\xi))  \geq 0$ for every vector $\xi\in
H(G)\cap z(\g,\hat G)^\bot$ and meromorphic reduction $\rho$ of $\hat {\cal Q}$ to
$\hat G(\xi)$,
\item  stable, if it  is semistable, its infinitesimal stabilizer
$\g_{\hat {\cal Q}}$ coincides with
$z(\g,\hat G)$ and the semistability inequality above is strict when
$\xi\ne 0$,
\item polystable, if it is semistable, its infinitesimal stabilizer $\g_{\hat {\cal Q}}$ is
reductive in the sense of Definition \ref{RedAlg}, and  the semistability inequality is
equality  when and only when $\rho$ is induced by a global holomorphic $Z_{\hat
G}(\xi)$-reduction and the section in $H^0(\hat {\cal Q}\times_\ad\g)$ defined by
$\xi$ via this reduction belongs to   $\g_{\hat {\cal Q}}$.
\end{enumerate}
\end{dt}

The condition $\xi\in  \zg^\bot$ becomes simply ``$\xi $ is semisimple"  when
$G=\hat G$ is connected.

\begin{re} Using the fact that, for a fixed parabolic subgroup $L\subset \hat G$ and
maximal compact subgroup
$K\subset G$,  the map $\xi\mapsto\deg(\rho,h(\xi))$ is affine on the convex set
$\{s\in i\kg|\ \hat G(s)\supset L\}$,  it follows easily that it suffices to check the
conditions   above only for vectors
$\xi$, for which $\hat G(\xi)$ is maximal in the set
$${\cal P}ar_G({\hat G}) :=\{ \hat G(\xi)|\ \xi\in H(G)\}\ .
$$
In this way, one recovers the usual (semi, poly) - stability condition for vector
bundles in the classical cases $(G,\hat G)=(GL(r,\C),GL(r,\C))$, $(G,\hat
G)=(SL(r,\C),GL(r,\C))$.
\end{re}

 Our main theorem becomes
\begin{thry}\label{KHB} Fix any  $K_0$ reduction $P_0\subset Q_0$. A
${\cal Q}_0$-oriented holomorphic structure $\hat{\cal Q}$ on $\hat Q$ is
polystable if and only if there exists a $\hat K$ reduction $\hat P$ of
$\hat Q$ inducing $P_0$ on $Q_0$ such that the Chern connection $\hat
A$ of the pair $(\hat{\cal Q},\hat P)$ satisfies the Hermitian-Einstein
equation
\begin{equation}
\label{HEB}
{\rm pr}_{i\kg}[i\Lambda_g F_{\hat A}]=\frac{2\pi
\tau(\hat {\cal Q})}{(n-1)!Vol_g(X)}\ .
\end{equation}
\end{thry}

Note again that the right hand term (the ``Einstein constant") of this equation is a
topological invariant of $\hat Q$ in the K\"ahlerian case and a holomorphic invariant of
$\hat{\cal Q}$ in the Gauduchon case.

\subsubsection{Isomorphism of moduli spaces}\label{IsoBdls}

As explained in the introduction, a Kobayashi-Hitchin correspondence always yields an
isomorphism between a moduli space of solutions of a Hermitian-Einstein type
equation   (which is a very difficult non-linear PDE system) and a moduli stable of
polystable objects, which in many cases can be described with complex geometric
methods. This idea was first used in Donaldson theory with spectacular success.

The formalism explained in \cite{DK} and \cite{LT} for comparing moduli spaces     can
be easily extended for  general oriented  bundles as follows.\\

Let $\bar{\cal A}:=\bar{\cal A}_{J_0}(\hat Q)$ be the space of oriented
almost complex structures
$\hat J$ on
$\hat Q$ which project  onto $J_0$ (the almost complex structure of
${\cal Q}_0$).
$\bar{\cal A}$ is an affine space with $A^{0,1}(\hat Q\times_\ad \g)$ as model
vector space.

Let ${\cal H}={\cal H}_{J_0}(\hat Q)\subset \bar{\cal A}$ be the subset
of integrable almost complex structures,  ${\cal H}^{\rm s}$ the open
subset of simple structures and ${\cal H}^{\rm st}_g\subset {\cal
H}^{\rm s}$ the subspace of stable structures.

On the differential geometric side we choose a $\hat K$-reduction $\hat
P$ of $\hat Q$, we denote by $P_0$ the induced $K_0$-reduction of
$Q_0$ and  by $A_0$ the Chern connection of the pair $(P_0,J_0)$. Let
${\cal A}:={\cal A}_{A_0}(\hat P)$ be the space of $A_0$-{\it oriented
connections}, i.e. the space of connections on
$\hat P$ which project onto the fixed connection $A_0$. A connection $\hat
A\in {\cal A}$ will be called {\it irreducible} if its infinitesimal
stabilizer with respect to the
${\cal K}$-action is minimal, i.e. it coincides with the
centralizer $z(\kg,\hat K)$ of $\hat K$ in $\kg$. As usually, we will
use the supscript
$(\cdot )^*$ to denote the open irreducible part and the supscript
$(\cdot)^{\rm HE}$ to denote the Hermitian-Einstein part, i.e. the locus
of {\it integrable} solutions of equation (\ref{HEB}).

One can introduce a   stronger irreducibility condition as
follows: a connection $\hat A\in{\cal A}$ is called {\it strongly
irreducible} if its stabilizer ${\cal K}_A\subset {\cal K}$ is minimal,
i.e. it coincides with the centralizer $Z_{K}(\hat K)$. Since this
condition is  not easily checked, and since we want to avoid unimportant
technical complications, we will impose a natural condition on our pair of
Lie groups
$(K,\hat K)$ which assures that the two irreducibility conditions  coincide,
and which holds for the standard pairs of classical groups.  One can easily
show

\begin{pr} Suppose that $K$ is connected and that for every $k\in K$
there exists $\lambda\in\kg$ which is invariant under $Z_{\hat K}(k)$
and such that
$\exp(\lambda)=k$. Then every irreducible connection $\hat A\in{\cal
A}$ is strongly irreducible.
\end{pr}

For the remainder of this section, we will always assume that the
pair $(K,\hat K)$ satisfies the hypothesis of this proposition.

The  fundamental objects we study in this section are the moduli spaces
$${\cal M}^{\rm s}:=\qmod{{\cal H}^{\rm s}}{{\cal G}}\ ,\ {\cal M}^{\rm
st}_g:=\qmod{{\cal H}^{\rm st}_g}{{\cal G}}\ , [{\cal
M}^{HE}]^*:=\qmod{[{\cal A}^{\rm HE}]^*}{{\cal K}}\ .
$$

 The same methods as in
\cite{Do1},
\cite{LT} show that Theorem
\ref{KHB} can be refined and that one has a canonical isomorphism of
moduli spaces:
\begin{thry}\begin{enumerate}
\item ${\cal M}^{\rm s}$ has a natural structure of a complex analytic
space (in general non-Hausdorff) of finite dimension and $[{\cal
M}^{HE}]^*$ has a natural structure of a Hausdorff real analytic
space of finite dimension.
\item The Chern map $\hat A\mapsto J^{\hat A}$ induces a real analytic
open embedding
$$[{\cal M}^{HE}]^*\textmap{KH}  {\cal M}^{\rm s}
$$
whose image is ${\cal M}^{\rm st}_g$. In particular, this space is open
in ${\cal M}^{\rm s}$, is Hausdorff and is canonically isomorphic via
$KH$ to $[{\cal M}^{HE}]^*$.
\end{enumerate}
\end{thry}

\begin{re} \begin{enumerate}
\item  In general, the moduli space $[{\cal
M}^{\rm HE}]^*$ can  contain orbifold type singularities if it
contains points which are not   strongly irreducible. Our condition on
the pair
$(K,\hat K)$ avoids this complication.

\item One can show that the isomorphism $[{\cal
M}^{HE}]^*\textmap{KH}  {\cal M}^{\rm st}$ extends to a ho\-meo\-morphism
$[{\cal M}^{HE}]\textmap{KH}  {\cal M}^{\rm pst}_g$, where
$${\cal M}^{\rm
pst}_g:=\qmod{{\cal H}^{\rm pst}_g}{{\cal G}}
$$
is the moduli space of polystable holomorphic structures. This space,
in general, does not have  a complex space structure; moreover the
local structure around a non-stable polystable point can be very
complicated (see e.g. \cite{Te4}).
\end{enumerate}
\end{re}

\paragraph{Why are moduli spaces of oriented connections important?} A clear,
simple motivation for considering moduli spaces of {\it oriented} connections
is the following: extending the  Witten conjecture to non-simply connected
4-manifolds leads naturally to Donaldson moduli spaces of oriented,
projectively ASD $U(2)$-connections. This spaces cannot be identified in
general with moduli spaces of $PU(2)$-connections. We explain in detail this
phenomenon:

Let $E$ be a Hermitian 2-bundle on an oriented Riemannian 4-manifold $M$,
 $P_E$ its frame bundle, and $L$ its determinant line bundle. Fix
a unitary connection
$A_0\in{\cal A}(L)$,   denote by ${\cal A}_{A_0}$   the space of
$A_0$-oriented connections, i.e. the space of unitary connections on
$E$ which induce $A_0$ on $L$. Consider the  subspace
$${\cal A}_{A_0}^{\rm ASD}:=\{A\in{\cal A}_{A_0}(E)\ |\ (F_A^0)_+=0\}$$
of {\it projectively  ASD  oriented connections}. Here $F_{A}^0$ stands for
the trace-free part of the curvature, which can be identified with the
curvature of the induced $PU(2)$-connection $\bar A$ on the associated
$PU(2)$-bundle
$\bar P:=\qmod{P_E}{S^1}$. Consider the  moduli space
$${\cal M}^{\rm ASD}_{A_0}(E):=\qmod{{\cal A}_{A_0}^{\rm ASD}}{{\cal K}}
$$
of projectively ASD oriented connections on $E$ modulo the gauge group
$${\cal
K}:=\Gamma(X,SU(E))\ .$$
 The space  ${\cal A}^{\rm
ASD}_{A_0}(E)$ can be identified with the space ${\cal A}^{\rm ASD}(\bar
P)$ of ASD connections on $\hat P$, but the natural  map ${\cal K}\to {\rm
Aut}(\bar P)$ between the two gauge groups is in general not surjective.
 There  is a natural map
$${\cal M}^{\rm ASD}_{A_0}(E)\map {\cal M}^{ASD}(\bar P)
$$
which is always surjective and is an isomorphism when $H^1(M,\Z_2)=0$. In the
general case, one has a natural action of $H^1(M,\Z_2)$ on ${\cal
M}^{\rm ASD}_{A_0}(E)$ given by tensorizing with flat $\{\pm
1\}$-connections, and
${\cal M}^{\rm ASD}(\hat P)$ is just the $H^1(M,\Z_2)$-quotient
(\cite{LT},
\cite{Te4}).
This natural
$H^1(M,\Z_2)$-action might have fixed points \cite{Te4},
so the space  ${\cal M}^{ASD}(\bar P)$ is in general more singular than
${\cal M}^{\rm ASD}_{A_0}(E)$. But there is another important  reason why the
latter space is a more convenient object: in order to embed a Donaldson
moduli space in a moduli space of non-abelian monopoles (this is the
starting point of the ``cobordism strategy" for the proof of the Witten
conjecture), one needs coupled Dirac operators (see section
\ref{NAG}), so one needs   unitary connections.   There is no obvious way to
define a coupled Dirac operator associated with a $PU(2)$-connection.

\subsubsection{The
canonical Hermitian metric on moduli spaces of oriented connections}\label{bcm}

We write $Z$ and $z$ for $Z_{K}(\hat K)$ and $z(\kg,\hat K)$ respectively. Let

$${\cal K}_r:=\qmod{{\cal K}}{Z}\ ,\ A^0_r(\hat
P\times_\ad\kg):=z^{\bot}\subset A^0 (\hat
P\times_\ad\kg)={\rm Lie}({\cal K})
$$
be the {\it
reduced gauge group} and its Lie algebra. We will also need the
orthogonal projection
$$p_r:A^0 (\hat
P\times_\ad\kg)\map A^0_r (\hat
P\times_\ad\kg)\ .
$$
 This group acts freely on our configuration space of
oriented irreducible connections
${\cal A}^*$ from the right, so ${\cal A}^*$ can be regarded as a
principal ${\cal K}_r$-bundle over the quotient ${\cal B}^*:=\qmod{{\cal
A}^*}{{\cal K}_r}$.

Write $\hat K_r:=\qmod{\hat K}{Z}$ and consider the principal
$\hat K_r$-bundle
$$\hat P_r:=\qmod{\hat P}{Z}:=\hat P\times_{\hat K} \hat K_r\ .$$

The associated bundle
$$\hat\P_r:={\cal A}^*\times_{{\cal K}_r}\hat P_r
$$
with standard fibre $\hat P_r$ over ${\cal B}^*$ can be alternatively
regarded as a principal $\hat K_r$-bundle over ${\cal B}^*\times X$ and
will be  called the {\it universal} $\hat K_r$-bundle.

We will introduce a canonical
connection on this principal bundle  (compare with \cite{DK},
section 5.2). Since our base metric is not K\"ahler (but only Gauduchon)
one should not use the
$L^2$-orthogonal distribution used in classical gauge theory (see
\cite{DK} p. 196).

Following the strategy in \cite{LT} section 5.3  we put, for any $A\in
{\cal A}^*$,
$$\delta_{\hat A}:= \Lambda_g\circ d^c_{\hat A}: A^1(\hat
P\times_\ad\kg)\map A^0(\hat P\times_\ad\kg)\ , $$
$$q_{\hat A}:=p_r\circ
\delta_{\hat A}\circ d_{\hat A}:A^0_r(\hat P\times_{\ad}\kg)\to
A^0_r(\hat P\times_{\ad}\kg) \ .
$$
The operator $\delta_{\hat A}$ plays the role of the $L^2$-adjoint of
$d_{\hat A}$ in the K\"ahlerian framework, whereas $q_{\hat A}$ replaces the
standard Laplacian.

The same methods as in \cite{LT} show that
\begin{pr}\label{conn1} After suitable Sobolev completions the following holds: The assignment
$$\hat A\mapsto \ker(p_r\circ \delta_{\hat A})
$$
defines a real analytic connection $\Gamma$ in the ${\cal K}_r$-principal
bundle
${\cal A}^*\to  {\cal B}^*$.
\end{pr}
\pf  We recall briefly the argument in \cite{LT} in this more general framework. First  note that, in general,  for any Euclidean connection $B$ on a Euclidean real bundle $W$ over $X$, the following   identity holds pointwise:
\begin{equation}
\frac{1}{2} \Lambda_g d^c d\vert v\vert^2=(\Lambda_g d_B^cd_B v,v)-\vert d_B v|^2
\end{equation}

We show first that
\begin{equation}\label{slice}
\ker(p_r\circ \delta_{\hat A})\cap T_{\hat A}[ \hat A\cdot{\cal K}_r]=0
\end{equation}
for every   connection $\hat A\in{\cal A}^*$. Indeed, for a an element $u\in A^0_r(\hat P\times_\ad \kg)={\rm Lie}({\cal K}_r)$, one has
$$p_r \delta_{\hat A} d_{\hat A} u=p_r  \Lambda_g d_{\hat A}^c d_{\hat A} u\  .
$$
The vanishing of this implies $\Lambda_g d_{\hat A}^cd_{\hat A}  u =\zeta \in z$, so
$$\frac{1}{2} \Lambda d^c d\vert u\vert^2+\vert d_{\hat A} u|^2= (\Lambda_g d_{\hat A}^cd_{\hat A}  u,u)=  (\zeta, u)\ .
$$
Integrating on $X$ and taking into account that  $g$ is a  Gauduchon metric,
one gets $\| d_{\hat A} u\|^2_{L^2}=0$. Since our oriented connection is
irreducible, this implies $u=0$, so (\ref{slice}) holds.   But this argument
also shows that the operator  $q_{\hat  A}$ is injective; but it is easy to see
that it is an index 0 Fredholm operator, so it is also surjective. Any
$\alpha\in A^1(\hat P\times_\ad\kg)$ can be written as
$$\alpha= d_{\hat A}(q_{\hat A}^{-1}(p_r\delta_{\hat A}(\alpha)))+ (\alpha-d_{\hat A}(q_{\hat A}^{-1}(p_r\delta_{\hat A}(\alpha))))\  ,
$$
and it is obvious that the first term belongs to  $T_{\hat A}[ \hat A\cdot{\cal K}_r]$, whereas the second to $\ker(p_r\circ \delta_{\hat A})$. Therefore, $\Gamma$ is a horizontal distribution, and it is
easily seen to be ${\cal K}_r$-invariant.
\qed
\begin{re} The curvature $F_\Gamma$ of the connection $\Gamma$ is given by
the formula
$$F_\Gamma(\alpha,\beta)= -q_{\hat A}^{-1} p_r\Lambda_g[\alpha^c\wedge
\beta]\ .
$$
for $\Gamma$-horizontal vectors $\alpha$, $\beta$.
\end{re}
(see \cite{LT} p. 177). Note that  $\Lambda_g[\alpha^c\wedge
\beta]$ is obtained by coupling the inner product $g$ on 1-forms with the Lie bracket
on the Lie algebra $\kg$. Therefore it coincides with the term $\{a,b\}$ appearing in \cite{DK}.

The connection $\Gamma$ can be used to define a {\it canonical Hermitian metric}
$H$  on the smooth part $[{\cal M}^{\rm HE}]^*_{\rm reg}$ of the moduli space
$[{\cal M}^{\rm HE}]^*\simeq{\cal M}^{\rm st}_g$. Recall  that we have fixed an
$\ad$-invariant complex inner product of Euclidean type $h$ on the complex
Lie algebra $\hat\g$.  We put
$$\Omega(a,b):= -\int_X h(\alpha\wedge \beta)\wedge\omega_g^{n-1}\ ,
$$
where $a$, $b$ are tangent vectors  to the moduli space and  $\alpha$,
$\beta$ are $\Gamma$-horizontal lifts.  This form has type $(1,1)$ with
respect to the complex structure induced from the canonical complex
structure ${\cal J}$ of $A^1(\hat P\times_\ad \kg)\simeq A^{0,1}(\hat
Q\times_\ad \g)$ on the ${\cal J}$-invariant spaces $\Gamma_{\hat A}\cap
T_{\hat A}({\cal A}^{\rm HE})$. It  is also easy to see that this form
is positive with respect to ${\cal J}$, so it defines a Hermitian metric $H$
as claimed.

Note that using the standard  $L^2$-orthogonal distribution  (as in
\cite{DK} section 5.3) will in general only yield a Riemannian metric, not
a Hermitian one when $g$ is not K\"ahler.

On the product ${\cal A}^*\times \hat P_r$ we   define now a distribution
$D$ which is horizontal with respect to the projection  ${\cal A}^*\times
\hat P_r\to {\cal B}^*\times X$. We put
$$D_{(\hat A,p)}:=\Gamma_{\hat A}\times \hat A_p\subset
T_{(\hat A,p)}({\cal A}^*\times
\hat P_r)\ ,
$$
where $\Gamma_{\hat A}$ is the horizontal space of $\Gamma$ at $\hat A$
and
$\hat A_p$ is the horizontal space at $p$ of the connection  $\hat A_r$
induced by $\hat A$ on
$\hat P_r$. This distribution is invariant with respect to the diagonal
$\hat {\cal K}_r$-action and the right $\hat K_r$-action on $\hat P_r$,
hence it descends to a   connection $\A$ on  the universal $\hat K_r$-bundle
$\hat\P_r$.

The curvature $\F$ of this connection is (see \cite{DK}, p. 197 and
\cite{LT} p. 177 for the Gauduchon case):
\begin{equation}\label{univcurv}
\F_{[\hat A,p]}(u,v)=F_{\hat A_r,p}(u,v)\ ,\ \F_{[\hat
A,p]}(\alpha,v)=a(v)\ ,\
\F_{[\hat A,p]}(\alpha,\beta)=F_\Gamma(\alpha,\beta)
\end{equation}
for tangent vectors $u$, $v\in T_x(X)$ and $\alpha$, $\beta\in \Gamma_{\hat A}$.
\begin{re}\label{integrable} The restriction of   $\A$ to $\resto{\P}{[{\cal
M}^{\rm HE}]^*_{\rm reg}\times X}$  is an integrable connection.
\end{re}

Using our $\ad$-invariant inner product $h$, we get an
associated  {\it closed characteristic form} $h(\F,\F)\in A^4( {\cal B}^*
\times X)$.

For any  $\nu\in  A^{2n-2}(X)$ and $a$, $b\in T_{[\hat A]}({\cal
B}^*)$  we get as in \cite{DK} p. 179
$$[p_{ {\cal B}^*}]_*[h(\F,\F)\wedge
p_X^*(\nu)](a,b)=\int\limits_X\left[- h(\alpha\wedge \beta) +2
h\left(F_{\hat A_r},F_\Gamma(\alpha,\beta)\right)\right]\wedge\nu
$$
$$
=\int\limits_X\left[- h(\alpha\wedge \beta) +2
h\left({\rm pr}_{\kg}(F_{\hat
A_r}),F_\Gamma(\alpha,\beta)\right)\right]\wedge\nu\ ,
$$
where  $\alpha$ $\beta\in\Gamma_{\hat A}$ are horizontal lifts of $a$ and $b$
at $\hat A$,
$p_{{\cal B}^*}$,
$p_X$ stand for the respective projections, and $[p_{{\cal B}^*}]_*$ denotes
fibre integration along the $X$-fibres.

Therefore, taking $\nu=\omega_g^{n-1}$ and using  the Hermitian-Einstein
equation (\ref{HEB}), on the smooth part
$[{\cal M}^{\rm HE}]^*_{\rm reg}$ of the moduli space it holds
\begin{equation}\resto{[p_{ {\cal B}^*}]_*[h(\F,\F)\wedge
p_X^*(\omega_g^{n-1})]}{[{\cal M}^{\rm HE}]^*_{\rm reg}} =\Omega
\end{equation}

Taking into account that $\resto{h(\F,\F)}{[{\cal M}^{\rm
HE}]^*_{\rm reg}\times X}$ is a closed $(2,2)$ form (see Remark
\ref{integrable}),  we get
\begin{thry} The $(1,1)$-form $\Omega$ associated with the canonical
Hermitian metric $H$ on the  smooth part  $[{\cal
M}^{\rm st}_g]_{\rm reg}$ of the moduli space ${\cal
M}^{\rm st}_g$  has the following property:
\begin{enumerate}
\item it is K\"ahler, when $d(\omega_g^{n-1})=0$ (i.e. when $g$ is {\it semi-K\"ahler}),
\item it satisfies the identity
$\partial\bar\partial \Omega=0$
when $g$ is Gauduchon.
\end{enumerate}
\end{thry}

Note that, for the second statement, one needs that $\resto{h(\F,\F)}{[{\cal
M}^{\rm HE}]^*_{\rm reg}\times X}$ is $\partial$ and $\bar\partial$-closed, so
it is  important to know that this form has pure type $(2,2)$.

This result shows that the moduli spaces of stable oriented bundles
over a compact complex manifold have  interesting differential geometric
properties (see Remark \ref{special} below).

\subsection{Moduli spaces of oriented pairs}

\subsubsection{Isomorphy of moduli spaces}

The construction of the moduli spaces of (oriented) holomorphic pairs in the
non-linear case (when the action $\hat \alpha$ is not a linear representation)
is substantially more difficult than in the classical moduli theories.

Indeed, in the non-linear case, it is quite difficult  to construct  the Sobolev
completions of the space of sections in the associated bundle $E:=\hat
Q\times_{\hat G} F$ and to endow these completions with the structure of a
Banach manifold (see \cite{Pa}).

Since here we are only interested in geometric applications, we will skip these
analytical difficulties, which will be addressed in   future work.

Fixing a holomorphic structure $J_0$ on $Q_0$, our complex geometric configuration space is
$$\bar{\cal A}:=\bar{\cal A}_{J_0}(\hat Q)\times \Gamma(X,E)
$$
with    the  natural right  action of the complex gauge group ${\cal G}=\Gamma(X,\hat Q\times_{\Ad} G)$.

Let ${\cal H}\subset\bar{\cal A}$  be the subspace of integrable pairs, ${\cal H}^{\rm s}$ the open subset of {\it strongly simple} integrable pairs, i.e. of pairs with trivial stabilizer. Note that the stabilizer of a simple integrable pair can be in general non-trivial; in order to avoid complications caused by orbifold-type singularities in the moduli spaces, we will restrict ourselves to the strongly simple points.   Similarly, $\bar{\cal A}^{\rm s}$ will denote the subspace of $\bar{\cal A}$ consisting of pairs with trivial stabilizer. The quotient
$${\cal M}^{\rm s}:=\qmod{{\cal H}^{\rm s}}{\cal G}
$$
has the structure of a finite dimensional, in general non-Hausdorff complex  space, which will be called the moduli space of strongly simple oriented holomorphic pairs. This space can be regarded as the vanishing locus of a   section in a vector bundle over $\bar{\cal B}^{\rm s}:=\qmod{\bar{\cal A}^{\rm s}}{\cal G}$ whose fibre over $[\hat J,\varphi]$ is
$$A^{0,2}(\hat Q\times_\ad \g)\times A^{0,1}(\varphi^*(V))\ ,
$$
where $V\to E$ is the vertical tangent bundle of $E$. This section is induced by the left hand side of the integrability condition.

Fix now a $\hat K$-reduction $\hat P$ of $\hat Q$, let $P_0$ be the induced reduction of $Q_0$ and $A_0$ the Chern connection of the pair $(P_0,J_0)$. The configuration space is now
$${\cal A}:={\cal A}_{A_0}(\hat P)\times \Gamma(X,E)
$$
acted on by the gauge group ${\cal K} :=\Gamma(\hat P\times_\ad K)$.
We will indicate by $(\ )^*$ the respective parts of strongly irreducible pairs (i.e. of pairs with  trivial stabilizer), and by $(\ )^{\rm HE}$ the subspaces of integrable pairs solving the generalized Hermitian-Einstein equation
$${\rm pr}_{i\kg}\left[i\Lambda_g F_{\hat A}\right]+ i\mu(\varphi)=0\ .
$$

We are interested in the  moduli spaces
$${\cal M}^{\rm HE}:=\qmod{{\cal A^{\rm HE}}}{\cal K}\ ,\ [{\cal M}^{\rm HE}]^*:=\qmod{[{\cal A^{\rm HE}]^*}}{\cal K}\ .
$$
Both moduli spaces are Hausdorff,  and the second moduli space can be described as the vanishing locus of a section $s$ in a bundle $S$ over the infinite dimensional quotient ${\cal B}^*:=\qmod{{\cal A}^*}{\cal K}$. Note however, that this space has no natural real-analytic structure in the general non-linear framework (even when $\mu$ is a real analytic map on $F$). The point is that the bundle $S$ will  only be a ${\cal C}^k$-bundle, for a regularity class $k$ depending on the chosen Sobolev completion of $\Gamma(X,E)$.

The Kobayashi-Hitchin correspondence becomes
\begin{thry}\begin{enumerate}

\item The   map $(\hat A,\varphi)\mapsto (J^{\hat A},\varphi)$ induces an open  embedding
$$[{\cal M}^{HE}]^*\textmap{KH}  {\cal M}^{\rm s}
$$
whose image is ${\cal M}^{\rm st}_g$. In particular, this space is open
in ${\cal M}^{\rm s}$, is Hausdorff and is canonically isomorphic via
$KH$ to $[{\cal M}^{HE}]^*$.
\end{enumerate}
\end{thry}

As in the bundle case, it holds:
\begin{re}
The isomorphism $[{\cal
M}^{HE}]^*\textmap{KH}  {\cal M}^{\rm st}$ extends to a ho\-meo\-mor\-phism
$[{\cal M}^{HE}]\textmap{KH}  {\cal M}^{\rm pst}_g$, where ${\cal M}^{\rm pst}_g$ is the moduli space of polystable oriented holomorphic pairs.
 \end{re}

\subsubsection{The canonical Hermitian metric on moduli spaces of oriented pairs}

\paragraph{The canonical Hermitian metric on the moduli space.}

As in section \ref{bcm}, one can define a natural connection on the   principal  ${\cal K}$-bundle ${\cal A}^*\to{\cal B}^*$.  Let $E$ be the associated bundle $E:=\hat P\times_{\hat \alpha} F$ and $V\to E$ the
vertical tangent bundle of $E$.  In this section we will always assume that the
K\"ahlerian structure $(J_F,\omega_F)$ on $F$ is $\hat K$-invariant. We will
denote by the same symbols the induced complex, respectively symplectic
structure on the fibres of $E$ and $V$.

 The tangent space of ${\cal A}^*$ at a pair $\ag=(\hat A,\varphi)$ is the product $A^1(\hat P\times_\ad\kg)\times A^0(\varphi^*(V))$, whereas the tangent space to the orbit $\ag\cdot{\cal K}$ is
$T_\ag({\ag\cdot{\cal K}})=\im(\dg_\ag)$
for the operator
$$\dg_\ag:A^0(\hat P\times_\ad\kg)\to A^1(\hat P\times_\ad\kg)\times A^0(\varphi^*(V)) \  ,\  \dg_\ag(k):=(d_{\hat A} k,\varphi^*(k^\#))\ ,$$
where  $k^\#$ denotes  the section of $V$ defined by $k$.

We define for $\ag=(\hat A,\varphi)\in{\cal A}^*$  the operators
$$\delta_\ag:A^1(\hat P\times_\ad\kg)\times A^0(\varphi^*(V))\map A^0(\hat P\times_\ad\kg)\  , \  \delta_\ag(\alpha,\psi)=\Lambda_g d_{\hat A}^c(\alpha)+d_\varphi \mu(J_F \psi)\  ,
$$
where  in the right hand term, $\mu$ is regarded as a map  $\Gamma(X,E)\to A^0(\hat P\times_\ad\kg)$ and $d_\varphi$ is the differential of this map at $\varphi$, and
$$q_\ag:A^0(\hat P\times_\ad\kg)\map A^0(\hat P\times_\ad\kg)\ ,\ q_\ag:=\delta_\ag\circ\dg_\ag\  .
$$

The precise result is the following:
\begin{pr}\label{conn2} The assignment
$
\ag\mapsto \ker(\delta_\ag)
$
defines a connection $\Gamma$ on the principal ${\cal K}$-bundle ${\cal A}^*\to {\cal B}^*$.
\end{pr}
\pf As in the proof of Proposition \ref{conn1} it suffices to show that the operator $q_\ag$ is injective. Using the fact $\mu$ can be regarded as a moment map for the  ${\cal K}$-action on the space of sections $\Gamma(X,E)$, we get pointwise in $X$
$$(q_\ag(k),k)=(\Lambda_g d_{\hat A}^c d_{\hat A} k+d_\varphi \mu(J_F \varphi^*(k^\#)),k)=$$
$$=(\Lambda_g d_{\hat A}^c d_{\hat A} k,k)+\omega_F(\varphi^*(k^\#) ,J_F \varphi^*(k^\#))=\frac{1}{2}\Lambda_gd^c d |k|^2 +|d_{\hat A} k|^2+|\varphi^*(k^\#)|^2\ .
$$
Hence, if $q_\ag(k)=0$, then integrating on $X$ we get $d_{\hat A} k=0$, $\varphi^*(k^\#)=0$, so $k=0$, because the oriented pair $(\hat A,\varphi)$ is irreducible.
\qed

Again this connection can be used to define a  canonical Hermitian metric
$H$  on the smooth part of the moduli space $[{\cal M}^{\rm
HE}]^*\simeq{\cal M}^{\rm st}_g$.   We put
$$\Omega(a,b):= -\int_X h(\alpha\wedge
\beta)\wedge\omega_g^{n-1}+(n-1)!\int_X
\omega_F(\psi,\chi) vol_g ,
$$
where $a$, $b$ are tangent vectors  to the moduli space and  $(\alpha,\psi)$,
$(\beta,\chi)$ are $\Gamma$-horizontal lifts at pair $\ag=(\hat A,\varphi)$.
Note again that our definition of the connection $\Gamma$ implies that this
form has type
$(1,1)$ with respect to the complex structure induced by the natural complex
structure
${\cal J}$ on ${\cal A}^*$ on the ${\cal J}$-invariant spaces $\Gamma_{\ag}\cap
T_{\ag}({\cal A}^{\rm HE})$.

\paragraph{Generalized symplectic reduction.}

Let $K$ be a Lie group, $\pi:P\to B$ a    principal $K$-bundle with a
connection $A$ and let $\eta$ be a closed $K$-invariant 2-form on $P$. We
will {\it not} assume that $\eta$ is non-degenerate. Let $m$ be a {\it
formal  moment map} for $\eta$, i.e. a $K$-equivariant  map $m:P\to \kg^\vee$
satisfying the usual condition
\begin{equation}
\langle dm, \xi\rangle=-\iota_{\xi^\#}\eta\  \forall \xi\in\kg\ .
\end{equation}

The sign on the right takes into account that   $K$ acts one on $P$
from the right. One can project
$\eta$ on the base
$B$ using restriction to the
$A$-horizontal  spaces, but the obtained 2-from $\eta_B$ will not be  closed
in general.

Let $F_A$ be the curvature 2-form of $A$. We define the 2-form $\langle
m,F_A\rangle$ on the base $B$ by
\begin{equation}
\langle m,F_A \rangle (a,b):=\langle m(p), F_{A,p}(\tilde a,\tilde b)\rangle\ ,
\end{equation}
where $a$, $b\in T_x(B)$, $\pi(p)=x$ and $\tilde a$, $\tilde b$  are horizontal
lifts of $a$ and $b$ at $p$.
\begin{pr} It holds
\begin{equation}\label{reduction}
d(\eta_B-2\langle m,F_A \rangle)=0\ .
\end{equation}
In particular,  if $0$ is a regular value of $m$, then the restriction of
$\eta_B$ to the generalized symplectic reduction $\qmod{m^{-1}(0)}{K}$ is
closed.
\end{pr}
\pf

For vector fields $a$, $b$ on $B$ one has the formula
$$[\tilde a,\tilde b]=\tilde{[a,b]}-F_A(\tilde a,\tilde b)^\#\ ,
$$
where $\tilde{\ }$ stands for horizontal lifts.
Taking into  account that $\eta$ is closed, one gets the following  formula for
$d\eta_B$:
$$d\eta_B(a,b,c)=-\eta(F_A(\tilde a,\tilde b)^\#,\tilde c)+\eta(F_A(\tilde
a,\tilde c)^\#,\tilde b)-\eta(F_A(\tilde b,\tilde c)^\#,\tilde a)=
$$
$$ \langle dm(\tilde c), F_A(\tilde a,\tilde b)\rangle -\langle dm(\tilde b),
F_A(\tilde a,\tilde c)\rangle+\langle dm(\tilde a),F_A(\tilde b,\tilde
c)\rangle= 2d(\langle m,F_A\rangle)(a,b,c) $$ by the Bianchi identity.
\qed

There is a   simple way to endow   $P$ with a pair $(\eta,m)$ consisting of
a $K$-invariant closed 2-form and a formal moment map $m$ for $\eta$:
choosing a section in a bundle associated with a Hamiltonian $K$-action on a
a symplectic fibre $F$:
\begin{re}\label{usefulRe} Let $(F,\omega_F)$ be a symplectic manifold endowed
with a left symplectic $K$-action, $\mu$ a moment map for this action  and
$f:P\to F$ a $K$-equivariant map.  Then $\eta:=f^*(\omega_F)$ is closed
and
$K$-invariant, and $m:=\mu\circ f$ is a formal moment map for $\eta$.
\end{re}
Note that the equivariance property of $f$ implies  the following
transformation rule for the fundamental vector fields:
$$f_*(u^\#_p)=-u^\#_{f(p)}\ \forall p\in P\ ,\ \forall u\in\kg\ .
$$

\paragraph{The universal bundle, the universal connection, the universal section.}

We regard the associated bundle
$$\hat\P:={\cal A}^*\times_{\cal K} \hat P
$$
over ${\cal B}^*$   as a $\hat K$-bundle over  ${\cal B}^*\times X$.  The map
$$\Phi:\hat\P\map F\ ,\  \Phi[\hat A,\varphi, p]:=\varphi(p)
$$
is well-defined and $\hat K$-equivariant, so it can be regarded as section in
{\it the universal associated bundle}
$$\E:=\hat\P\times_{\hat K} F\ .
$$

Using the horizontal spaces of   $\Gamma$ and of the variable connection $\hat A$, one constructs in the same way
as in section \ref{bcm}    a universal connection $\A$ on $\hat \P$.

The curvature $\F$ of this connection is given by formulae similar to
(\ref{univcurv})
$$\F_{[\hat A,\varphi,p]}(u,v)=F_{\hat A,p}(u,v)\ ,\ \F_{[\hat
A,\varphi,p]}((\alpha,\psi),v)=\alpha(v)\ ,$$
$$
\F_{[\hat
A,\varphi,p]}((\alpha,\psi),(\beta,\chi))=F_\Gamma((\alpha,\psi),
(\beta,\chi))
$$
for tangent vectors $u$, $v\in T_x(X)$ and horizontal vectors $(\alpha,\psi)$,
$(\beta,\chi)\in
\Gamma_{(\hat A,\varphi)}$.

\begin{re}\label{intpair} The restriction of the universal pair  $(\A,\Phi)$ to
$[{\cal M}^{\rm HE}]^*_{\rm reg}\times X$ is an integrable pair.
\end{re}

 For a  $(2n-2)$-form $\nu$ on $X$ we get as in the
previous section
\begin{equation}\label{firstid}
[p_{ {\cal B}^*}]_*[h(\F,\F)\wedge
p_X^*(\nu)](a,b)
=$$
$$=\int\limits_X\left[ -h(\alpha\wedge \beta) +2
h\left({\rm pr}_{\kg}(F_{\hat
A}),F_\Gamma((\alpha,\psi),(\beta,\chi))\right)\right]\wedge\nu \ ,
\end{equation}
where $a$, $b\in T_{[\hat A,\varphi]}({\cal B}^*)$ and $(\alpha,\psi)$,
$(\beta,\chi)$ are $\Gamma$-horizontal lifts of these vectors at $(\hat
A,\varphi)$. Note that, for $\nu=\omega_g^{n-1}$, the first term on the right coincides with the first term of $\Omega(a,b)$.
\\

Consider now the $\hat K$-invariant  closed form $\eta:=\Phi^*(\omega_F)$ on
$\hat\P$.   This form admits a formal moment map in
the sense of the previous paragraph, namely $m:=\mu\circ \Phi$ (see Remark
\ref{usefulRe}).

The projection of $\eta$ on the basis ${\cal
B}^*\times X$ is given by
$$\eta_{{\cal B}^* \times X}((a,u),(b,v))=\omega_F(\psi(x),\chi(x))\ ,
$$
where $u$, $v\in T_x(X)$, and $\psi$, $\chi\in \varphi^*(V)$ are the second
components of the $\Gamma$-horizontal lifts  $(\alpha,\psi)$, $(\beta,\chi)$ of
$a$ and $b$ in $(\hat A,\varphi)$.

Applying formula (\ref{reduction}) and taking into account that $h$ is
negative definite on $\kg$, we obtain
$$d[\eta_{{\cal B}^* \times X} + 2 h(m,\F)]=0\ .
$$
In particular, for every $2n$-form   $\sigma$ on $X$, the form obtained by fibre
integration along the $X$ fibres
$$[p_{{\cal B}^*}]_*\left[(\eta_{{\cal B}^* \times X} + 2 h(m,\F))\wedge
p_X^*(\sigma)\right]
$$
will be closed on ${\cal B}^*$. But obviously $h(m,\F)\wedge
p_X^*(\sigma)=h(m,F_\Gamma)\wedge p_X^*(\sigma)$. Therefore, we obtain
\begin{lm}\label{cclosed} For every $\sigma\in A^{2n}(X)$ the form
$\Omega_\sigma$ on
${\cal B}^*$ defined by
\begin{equation}\label{secondid}
\Omega_\sigma(a,b)= \int_X
\left[\omega_F(\psi,\chi) +2h(\mu(\varphi),
F_\Gamma((\alpha,\psi),(\beta,\chi))\right] \sigma
\end{equation}
is closed on ${\cal B}^*$.
\end{lm}

Note that the first term of $(n-1)!\Omega_{vol_g}(a,b)$ is just the second term  of $\Omega(a,b)$ (which was missing in (\ref{firstid})).

Specializing   (\ref{firstid}) and  (\ref{secondid}) to the case $\sigma=vol_g=\frac{1}{n!}\omega^n_g$,
$\nu=\omega_g^{n-1}$ and taking the sum of these formulae, we get
\begin{equation}
\left\{[p_{ {\cal B}^*}]_*[h(\F,\F)\wedge
p_X^*(\nu)]  +(n-1)!\Omega_{vol_g}\right\}(a,b) =$$
$$
=\Omega(a,b)+\int\limits_X  2
h\left(\omega_g^{n-1}\wedge{\rm pr}_{\kg}(F_{\hat
A})+(n-1)!vol_g\ \mu(\varphi)\ ,\ F_\Gamma((\alpha,\psi),(\beta,\chi))\right)=$$
$$=\Omega(a,b)+(n-1)!\int\limits_X  2
h\left(\Lambda_g{\rm pr}_{\kg}(F_{\hat
A})+ \mu(\varphi)\ ,\  F_\Gamma((\alpha,\psi),(\beta,\chi))\right)vol_g
\end{equation}

Restricting to the moduli space, one gets
\begin{equation}\resto{[p_{ {\cal B}^*}]_*[h(\F,\F)\wedge
p_X^*(\omega_g^{n-1})]}{[{\cal M}^{\rm HE}]^*_{\rm
reg}}+(n-1)!\resto{\Omega_{vol_g}}{[{\cal M}^{\rm HE}]^*_{\rm
reg}} =\Omega\ .
\end{equation}
Since the characteristic form $h(\F,\F)$ is closed of pure type (2,2) (see
Remark \ref{intpair}) and
$\Omega_{vol_g}$ is closed by Lemma
\ref{cclosed}, this shows
\begin{thry}\label{canmet} The $(1,1)$-form $\Omega$ associated with the canonical
Hermitian metric $H$ on the  smooth part $[{\cal
M}^{\rm st}_g]_{\rm reg}$ of the moduli space ${\cal
M}^{\rm st}_g$  has the following property:
\begin{enumerate}
\item it is K\"ahler, when $d(\omega_g^{n-1})=0$ (i.e. when $g$ is {\it semi-K\"ahler}),
\item it satisfies the identity
$\partial\bar\partial \Omega=0$,
when $g$ is Gauduchon.
\end{enumerate}
\end{thry}
\begin{re}\label{special} We believe that, for a given compact complex manifold
$Y$ of dimension $d\geq 3$, the existence of Hermitian metric $h$ satisfying  $\partial\bar\partial
\omega_h=0$ is a very restrictive condition. Therefore,  smooth compact
complex subspaces of any  moduli space ${\cal
M}^{\rm st}_g$ of oriented holomorphic pairs (or  oriented holomorphic
bundles) have very special differential geometric properties.
\end{re}

\subsubsection{Application: Canonical metrics on Douady Quot spaces}

Let ${\cal E}_0$  be a fixed holomorphic vector bundle of rank $r_0$ on $X$ and $E$
a differentiable rank $r$ vector bundle. Denote by $Quot^E_{{\cal E}_0}$ the
Quot space of quotients $q:{\cal E}_0\to {\cal Q}$ of ${\cal E}_0$ with locally
free kernel of differentiable type $E$. This quot space, which is an open subspace of the Douady space  $Quot_{{\cal E}_0}$ of all quotients of ${\cal E}_0$, has a gauge theoretical
description \cite{OT3}: it can be identified with the moduli space of pairs
$({\cal E},\varphi)$, where ${\cal E}$ is a holomorphic structure on $E$ and
$\varphi:{\cal E}\to {\cal E}_0$ is a holomorphic {\it sheaf monomorphism}, modulo the gauge group ${\cal G}:=\Gamma(X,GL(E))$.

However, it is not clear at all whether the condition ``$\varphi$ is a sheaf monomorphism" is a
stability condition.

Using our general formalism explained in the introduction, the   moduli problem  for holomorphic pairs $({\cal E},\varphi)$ with $\varphi:{\cal E}\to {\cal E}_0$  corresponds to the system
$$\hat G=GL(r,\C)\times GL(r_0,\C)  ,\ G=GL(r,\C)  , G_0=GL(r_0,\C)\ ,
F=\Hom(\C^r,\C^{r_0})\ .
$$
  A moment map for the $K=U(r)$-action on  $F$ has the form
$$\mu_t(f)=\frac{i}{2 }f^*\circ f-it\id, \ t\in\R\ .$$
Using the explicit form of the maximal weight function associated with a  linear representation given in section \ref{linact}, one obtains the following
stability condition   corresponding to $\mu_t$. Put
$\tau:=\frac{(n-1)!Vol_g(X)}{2\pi}t$.
\begin{dt} A holomorphic pair  $({\cal E},\varphi)$ is called $\tau$-stable if for every non-trivial subsheaf ${\cal F}$ of ${\cal E}$, the following holds:
\begin{enumerate}
\item $\mu_g(\qmod{\cal E}{\cal F})> -\tau$ if $\rk({\cal F})<r$.
\item $\mu_g({\cal F})<-\tau$ if ${\cal F}\subset\ker(\varphi)$.
\end{enumerate}
\end{dt}
 Here $\mu_g$ denotes the $g$-slope ${\deg_g}/{\rk}$. Note that the  $g$-degree
$\deg_g$ is {\it not } a topological invariant for general Gauduchon metrics $g$
(see  \cite{LT}).

 Our next purpose is to compare the Quot space $Quot^E_{{\cal E}_0}$ with a moduli space of $\tau$-stable pairs. The result is
 \begin{thry}  Fix $E$ and ${\cal E}_0$ as above.
 \begin{enumerate}
 \item Suppose that $g$ is semi-K\"ahler. Then for sufficiently  large $\tau\in\R$, the following  are equivalent:
\begin{enumerate}
\item  $({\cal E},\varphi)$ is $\tau$-stable.
\item $\varphi$ is injective.
\end{enumerate}
\item Suppose that $g$ is Gauduchon.  Then
 for any $s\in\R$ there exists $\tau(d,{\cal E}_0)\in\R$ such that for every  holomorphic pair $({\cal E},\varphi)$, $\varphi:{\cal E}\to {\cal E}_0$ with $\deg_g({\cal E})\geq s$ and every $\tau\geq \tau(s,{\cal E}_0)$ the following  are equivalent:
\begin{enumerate}
\item  $({\cal E},\varphi)$ is $\tau$-stable.
\item $\varphi$ is injective.
\end{enumerate}
\end{enumerate}
 \end{thry}
 \pf   For the semi-K\"ahlerian case the proof can be found in \cite{OT3}.  The second statement follows using the same methods,  but taking into account that the $g$-degree is no longer a topological invariant in the Gauduchon case.
 \qed

  Let
$$d_g:Quot^E_{{\cal E}_0}\map \R
 $$
be  the real analytic map defined   by $d_g:({\cal E},\varphi):=\deg_g({\cal E})$.  The
second statement of the theorem above  shows that the open subspaces
$$[Quot^E_{{\cal E}_0}]_s:=d_g^{-1}(-s,\infty)$$
of $Quot^E_{{\cal E}_0}$ can be identified with moduli spaces of stable oriented
pairs.
 By Theorem \ref{canmet}, we obtain
 \begin{co}
 \begin{enumerate}
 \item Suppose that $g$ is semi-K\"ahler. Then the smooth part of the quot space
$Quot^E_{{\cal E}_0}$ admits a canonical K\"ahler metric.
 \item Suppose that $g$ is Gauduchon.  Then, for any $s\in\R$,   the smooth part of the
truncated quot space $[Quot^E_{{\cal E}_0}]_s$ admits a canonical Hermitian metric whose
K\"ahler form $\Omega$ satisfies the identity $\partial\bar\partial\Omega=0$.
 \end{enumerate}
  \end{co}

This result can be specialized to Douady spaces  spaces of effective
divisors.  For any class $m\in\H^2(X,\Z)$, let ${\cal D}ou(m)$ be the Douady
space of effective divisors in $X$ representing the homology class $PD(m)$.
{\it On   Gauduchon manifolds, these spaces are in general not compact}.
${\cal D}ou(m)$ can be identified with the quot space $Quot_{{\cal
O}_X}^{M^{\vee}}$, where $M$ is a differentiable line bundle with Chern
class $m$.

For every $s\in \R$, we define the {\it truncated} Douady space   ${\cal
D}ou(m)_s$ by
$${\cal D}ou(m)_s:=\{ D\in {\cal D}ou(m)\ | \  Vol_g(D)=\deg_g({\cal O}(D))<
s\}.
$$
These spaces are relatively compact in ${\cal D}ou(m)$, by Bishop's
compactness theorem. One has
$${\cal D}ou(m)_s=[Quot_{{\cal O}_X}^{M^{\vee}}]_s\ ,$$
so our result gives a canonical K\"ahlerian  metric  on  every ${\cal
D}ou(m)$ when $g$ is semi-K\"ahler, and a canonical metric satisfying
$\partial\bar\partial\Omega=0$ on every truncated space ${\cal D}ou(m)_s$.

 \begin{co}\label{SWM} On a Gauduchon surface $X$ with odd $b_1(X)$, the
moduli space of irreducible solutions of a (suitably perturbed) abelian
Seiberg-Witten equations is a disjoint union of two truncated Douady
spaces of effective divisors (see
\cite{OT7}). Therefore, one gets   canonical Hermitian metrics with the
mentioned property on all these moduli spaces of irreducible abelian
monopoles.
\end{co}

\begin{re}
 It might be interesting to  study the behavior of the obtained metrics on
$[Quot^E_{{\cal E}_0}]_s$ as $s\rightarrow \infty$.
\end{re}

\subsection{Non-abelian monopoles on Gauduchon surfaces}\label{NASW}

\subsubsection{Oriented rank 2 holomorphic pairs}

 Let $E$ be a rank 2 complex vector bundle over a compact Gauduchon
manifold $(X,g)$, and let ${\cal L}$ be a fixed holomorphic  structure on its
determinant line bundle $L:=\det(E)$.

 We are interested in the classification of {\it oriented holomorphic  pairs
of type} $(E,{\cal L})$, i.e. of pairs $({\cal E},\varphi)$ consisting of a
holomorphic structure ${\cal E}$ on
$E$ inducing ${\cal L}$ on $L$ and a holomorphic section $\varphi\in
H^0(X,{\cal E})$. This complex geometric  moduli problem corresponds -- in
our general setting -- to the system
 $$G=SL(2)\ ,\  \hat G=GL(2)\ ,\ G_0=\C^*\ ,\ F=\C^2\ .
 $$
 \begin{re} Non-oriented pairs have been considered by  many authors   (see
e.g.  \cite{Bra}, \cite{Th}, \cite{HL}). Oriented holomorphic  pairs have
been  first introduced in \cite{OT2}  and
\cite{Te2}  in the K\"ahlerian framework  in order  to give a complex
geometric interpretation for the moduli spaces of $PU(2)$-monopoles.
 \end{re}

The stability condition for oriented rank 2 holomorphic pairs is quite simple:

\begin{dt}
An oriented pair $({\cal E},\varphi)$ of type $({ E},{\cal L})$
is called
 \begin{enumerate}
 \item {\it stable} if one of the following conditions holds:
\begin{enumerate}
\item $\varphi=0$ and ${\cal E}$ is a $g$-stable holomorphic bundle;
\item $\varphi\ne 0$  and     $\mu_g({\cal
O}_X(D_\varphi))<\mu_g({\cal E})=\frac{1}{2}\deg_g({\cal L})$, where $D_\varphi$ is the divisorial component of the
vanishing locus $Z(\varphi)$;
\end{enumerate}   \item   $g$-{\it polystable} if  it is $g$-stable, or $\varphi=0$  and
${\cal E}$ is a $g$-polystable holomorphic bundle.
\end{enumerate}
\end{dt}
Note that,  contrary to the case of non-oriented pairs,  this stability condition is not parameter dependent  (compare to the previous
section in this article and with \cite{Bra}).
\begin{re} The stability condition for higher rank oriented pairs is much more difficult (see \cite{OST} for the algebraic case).
\end{re}

We will denote by ${\cal M}^{\rm st}(E,{\cal L})$ the moduli space of stable oriented pairs of type $(E,{\cal L})$. The stabilizer of a stable  oriented  pair of the form  $({\cal E},0)$   with respect to the natural action of the complex gauge group ${\cal G}:=\Gamma(X,SL(E))$
is $\Z_2$, so these moduli spaces have at least $\Z_2$-orbifold
singularities along  the locus $\varphi=0$.

Moduli spaces of oriented rank 2 holomorphic pairs are very interesting geometric objects for the following reason:
\begin{re} The moduli space ${\cal M}^{\rm st}(E,{\cal L})$ comes with a natural $\C^*$-action, given by
$(\zeta,({\cal E},\varphi))\mapsto ({\cal E},\zeta\varphi)$.  In the terminology of \cite{OST}, \cite{OT6},   ${\cal M}^{\rm st}(E,{\cal L})$ is a ``master space" associated to the coupling of holomorphic bundles with holomorphic sections.
\end{re}

The fixed point locus of this $\C^*$-action is described as follows:

\begin{lm} Suppose that $(X,g)$ is a Gauduchon surface, and denote
$$l:=c_1(E)\ ,\ \lambda:=\frac{1}{2}\deg_g({\cal L})\ ,\ c:=c_2(E)\ .$$
Then the  fixed point locus of the $\C^*$-action on  ${\cal M}^{\rm
st}(E,{\cal L})$ decomposes as
$${\cal M}^{\rm st}(E,{\cal L})^{\C^*}={\cal M}^{\rm st}_{\cal L}(E)\coprod \left[\coprod_{m(l-m)=c}{\cal D}ou(m)_{\lambda}\right]\ ,$$
where ${\cal M}^{\rm st}_{\cal L}(E)$ denotes the moduli space of $g$-stable holomorphic structures on $E$ inducing ${\cal L}$ on $L$.
\end{lm}

The truncated Douady spaces ${\cal D}ou(m)_{\lambda}$ have been defined in the previous section.
On higher dimensional manifolds one must take the union of ${\cal D}ou(m)_{\lambda}$, over all classes $m$ such that the line bundle $M$ of Chern class $m$ is a topological summand of $E$.
\\ \\
\pf The part ${\cal M}^{\rm st}_{\cal L}(E)$   is  the locus  $\varphi=0$, which is obviously part of the fixed point locus. A truncated Douady space ${\cal D}ou(m)_{\lambda}$ with $m(l-m)=c$ is embedded in ${\cal M}^{\rm st}(E,{\cal L})^{\C^*}$ as follows:  To every $D\in {\cal D}ou(m)_{\lambda}$ we associate the stable pair $({\cal O}_X(D)\oplus[{\cal L}\otimes {\cal O}_X(-D)],\varphi_D)$, where $\varphi_D$ is the canonical section of ${\cal O}_X(D)$.
It is easy to see that any fixed point of the $\C^*$-action has one of these two types.
\qed

\begin{re}  When $X$ is a surface, the space ${\cal M}^{\rm st}_{\cal L}(E)$
can be identified -- via the Kobayashi-Hitchin isomorphism --  with a
Donaldson moduli space of oriented,  projectively ASD
$U(2)$-connections. When the base is simply connected, this space can be
further identified with a moduli space of
$PU(2)$-instantons.

The  spaces ${\cal D}ou(m)_{\lambda}$ are related to Seiberg-Witten moduli spaces, as pointed out in Corollary  \ref{SWM}.  We will see that the total space ${\cal M}^{\rm st}(E,{\cal L})$ is part of a moduli space of non-abelian monopoles.

On surfaces, the moduli spaces ${\cal M}^{\rm st}(E,{\cal L})$ have natural
compactifications, which are complex schemes when the base surface is
algebraic.

Any non-constant $\C^*$-orbit in  ${\cal M}^{\rm st}(E,{\cal L})$ connects either a point of Donaldson type to a point of Seiberg-Witten type, or two points of Seiberg-Witten  type (one of which can belong to the compactification
of ${\cal M}^{\rm st}(E,{\cal L})$).
\end{re}
\vspace{3mm}

The complex geometric moduli problem for oriented holomorphic  pairs has the
following differential geometric analogue.

Fix a Hermitian metric $h$ on $E$ and let $A_0$ be the Chern  connection of
the pair  $({\cal L},\det(h))$. Our configuration space is the space ${\cal
A}:={\cal A}_{A_0}(E)\times A^0(E)$, where ${\cal A}_{A_0}(E)$ is the space
of unitary connections on $E$ which induce $A_0$ on $L$. The moduli space
of {\it oriented vortices} is defined by
$${\cal M}^{\rm HE}(E,A_0):=\qmod{{\cal A}^{\rm HE}}{{\cal K}}\ ,$$
where ${\cal A}^{\rm HE}$ is the space of integrable solutions of the
{\it projective vortex equation} (see \cite{Te2} for the K\"ahlerian case):
$$i\Lambda_g F_A^0 +\frac{1}{2}(\varphi\bar\varphi)_0=0\ ,  $$
and ${\cal K}$ is the gauge group $\Gamma(X,SU(E))$. Here $F_A^0$ denotes the trace-free part of the curvature, which can be identified with the curvature of the
associated $PU(2)$-connection.  A pair $(A,\varphi)$ will be called
irreducible if its infinitesimal stabilizer vanishes. The stabilizer of an
irreducible pair is either trivial or $\Z_2$.

In our case, the Kobayashi-Hitchin isomorphism gives:
\begin{thry}\label{ORVOR} Let $(X,g)$ be a Gauduchon manifold, $(E,h)$ a
Hermitian rank 2 bundle on $X$, ${\cal L}$ a holomorphic structure on
$L:=\det(E)$, and
$A_0$   the Chern connection of the pair $({\cal L}, \det(h))$.  The map
$(A,\varphi)\mapsto(\bar\partial_A,\varphi)$ induces a real analytic
isomorphism of moduli spaces
$$[{\cal M}^{\rm HE}(E,A_0)]^*\textmap{KH}{\cal M}^{\rm st}(E,{\cal L})\ .
$$
This isomorphism extends to a homeomorphism  ${\cal M}^{\rm
HE}(E,A_0)\simeq{\cal M}^{\rm pst}(E,{\cal L})$.
\end{thry}

\begin{re} \begin{enumerate}
\item With an appropriate  definition of (poly)stability for oriented
pairs (see \cite{OST} for the algebraic   and \cite{OT8} for  the
K\"ahlerian case), the theorem holds for arbitrary rank $r$.
\item On  surfaces, the   spaces  ${\cal M}^{\rm HE}(E,A_0)$ have
Uhlenbeck-type compactifications, which are complex schemes when the base is
algebraic.
\end{enumerate}
\end{re}

\subsubsection{Complex geometric description   of non-abelian
monopoles}\label{NAG}

It is known (see \cite{FL1}, \cite{FL2}, \cite{OT2} \cite{PT}, \cite{Te5})
that the most natural way to prove the Witten conjecture relating the
Donaldson invariants to the Seiberg-Witten invariants is to consider
certain  moduli spaces of non-abelian monopoles, which contain both
Donaldson moduli spaces and moduli spaces of abelian monopoles. Substantial
progress towards a complete proof of the Witten conjecture in full
generality (for arbitrary simple type manifolds) was obtained in
\cite{FL2}. The relation between the two type of differential topological
type of invariants has more than theoretical importance: it plays an
important role in the recent  spectacular results of Kronheimer-Mrowka in
3-dimensional topology
\cite{KM}.

On the other hand, the moduli spaces of non-abelian monopoles are interesting, difficult
and mysterious objects, so we believe that, independently of the Witten conjecture, it
is important and interesting to have methods for describing such
moduli spaces explicitly. The Kobayashi-Hitchin correspondence yields such a
method even on Gauduchon surfaces.

Let $M$ be a closed  oriented Riemannian 4-manifold and $\tau$  a $Spin^c$-structure on
$M$. Denote by $\Sigma^\pm$, $\Sigma:=\Sigma^+\oplus\Sigma^-$ the spinor bundles of
$\tau$, by $D:=\det(\Sigma^\pm)$ its determinant line bundle and by
$\gamma:\Lambda^1_M\to\End(\Sigma)$ the Clifford map.

Fix a Hermitian 2-bundle $E$ of determinant $L$ on $M$. We choose parameter
connections  $A_0$ and $B_0$ on $L$ and  $D$.

The configuration space for the non-abelian monopole equations is
$${\cal A}:={\cal A}_{A_0}(E)\times A^0(\Sigma^+\otimes E) \ ,
$$
and the gauge group is ${\cal K}:=\Gamma(M, SU(E))$. Fix a complex form $\beta\in A^1(M,\C)$. The monopole equations
associated with the data $(\tau, B_0,\beta,A_0)$ read
\begin{equation}\label{SWE}
\left\{\begin{array}{ccc}
\Dr_A\Psi +\gamma(\beta)&=&0\\
\gamma(F_A^0)&=&(\Psi\otimes\bar \Psi)_0
\end{array}\right.
\end{equation}
where   the term $(\Psi\otimes\bar \Psi)_0$ stands for the projection of
$$\Psi\otimes\bar \Psi\in A^0(\End(\Sigma)\otimes\End(E))$$
 on
$A^0(\End_0(\Sigma)\otimes\End_0(E))$, and $\Dr_A$ is the coupled Dirac
operator associated with the connection $A$ and  the abelian connection
$B_0$ on $D$. It is easy to see that this projection belongs to
$${\rm Herm}_0(\Sigma^+)\otimes{\rm Herm}_0(E)\textmap{\simeq\gamma^{-1}}
\Lambda^2_+\otimes su(E)\ .$$

Standards methods  (\cite{FL1}, \cite{Te5}) show that the moduli space
$${\cal M}:=\qmod{{\cal A}^{SW}}{{\cal
K}}$$
 of solutions of this equation is a finite dimensional real  analytic space
which has a natural Uhlenbeck-type compactification. This is constructed  in
the same way as for the moduli spaces of instantons in Donaldson theory,  by
adding a {\it finite} number of strata of the form $S^k(X)\times {\cal
M}_k$, where ${\cal M}_k$ are Seiberg-Witten moduli spaces associated with
bundles $E_k$ having $\det(E_k)=\det(E)$ and $c_2(E_k)=c_2(E)-k$, $k\geq 1$
(see \cite{FL1}, \cite{Te5}).
\vspace{3mm}

The moduli space ${\cal M}$ comes with a natural $S^1$-action (given by
scalar multiplication  of the spinor component) and one can easily see that
the fixed point locus ${\cal M}^{S^1}$ decomposes as a union of the
Donaldson moduli space ${\cal D}:={\cal M}^{\rm ASD}_{A_0}(E)$ of oriented
projectively ASD connections (see section \ref{IsoBdls}) and the union
${\cal S}{\cal W}$ of a finite number of   Seiberg-Witten moduli
spaces of abelian monopoles. The fundamental idea of the ``cobordism
strategy" (\cite{FL1},
\cite{OT2}
\cite{PT},
\cite{Te5}) is to remove a small neighborhood $U$ of ${\cal M}^{S^1}$ in
${\cal M}$ and to use the quotient
$\qmod{[{\cal M}\setminus U]}{S^1}$ as a ``homology equivalence"  between
${\cal D}$ and ${\cal S}{\cal W}$. The problem
is that new Seiberg Witten moduli spaces can appear in the lower strata of
the  Uhlenbeck compactification, and the structure of this compactification
is very complicated \cite{FL1}, \cite{FL2}.

Let now $(X,g)$ be a Gauduchon surface and  $\tau_{\rm can}$ the canonical $Spin^c$ structure  of $X$. The spinor bundles of $\tau_{\rm can}$ are:
$$\Sigma^+_{\rm can}=\Lambda^{0,0}\oplus\Lambda^{0,2}\ ,\  \Sigma^-_{\rm can}=\Lambda^{0,1}\ ,
$$
so that the determinant line bundle $D_{\rm can}$ is $K_X^\vee=\wedge^2(T^{1,0}_X)$. Denote by
$B_{\rm can}$ the connection on $D_{\rm can}$ induced by the Chern connection
of the pair $({\cal T}_X,g)$.

A spinor $\Psi\in\Sigma_{\rm can}^+\otimes E$ decomposes as $\Psi=\varphi+\alpha$, where $\varphi\in A^0(E)$ and $\alpha\in A^{0,2}(E)=A^0(E\otimes K_X^\vee)$.

We will suppose  that  the Chern class $c_1(E)$ belongs to the Neron-Severi group  $NS(X)$ of $X$, and that the  parameter connection $A_0$ is integrable.

Our perturbation 1-form will be
$$\beta_{\rm can}:=-\frac{1}{4}\theta_g\ ,
$$
where $\theta_g$ is the torsion form of $g$, defined by the equality
$d\omega_g=\omega_g\wedge\theta_g$. This particular choice is very
important:  the corresponding  perturbed coupled Dirac operator
$\Dr_{A}+\gamma(\beta_{\rm can})$ coincides with the
``Dolbeault-Dirac"-operator
$${\cal D}_A:=\sqrt{2}(\bar\partial_{A}+\bar\partial_{A}^*)\ .$$
(\cite{Bi}).  As in the  K\"ahlerian case \cite{Te2}, one can prove the
following
  {\it decoupling theorem}, which extends to the non-abelian
Seiberg-Witten theory the well-known Witten  decoupling theorem for abelian
monopoles on K\"ahlerian surfaces \cite{W}:
\begin{thry} A pair
$$(A,\varphi+\alpha)\in {\cal A}_{A_0}(E)\times\left[A^0(E)\oplus
A^{0,2}(E)\right]$$
solves the  monopole equations associated with the data  $(\tau_{\rm
can},B_{\rm can},\beta_{\rm can},A_0)$ if  and only if
$A$ is integrable and one of the following holds:
\begin{enumerate}
\item
$ \alpha=0,\ \bar\partial_A\varphi=0\ \ \ and\ \ \ i\Lambda_g
F_A^0+\frac{1}{2}(\varphi\bar\varphi)_0=0$,
\item $ \varphi=0,\  \partial_A\alpha=0\ \ and\ \ i\Lambda_g
F_A^0-\frac{1}{2}*(\alpha\wedge\bar\alpha)_0=0$.
\end{enumerate}
\end{thry}

This theorem shows that the moduli space  ${\cal M}$ of solutions of
the non-abelian Seiberg-Witten equations associated with the data $(\tau_{\rm
can},B_{\rm can},\beta_{\rm can},A_0)$ decomposes as union of two closed
subspaces  ${\cal M}_I$,
${\cal M}_{II}$ intersecting along the Donaldson moduli space ${\cal
M}_{A_0}(E)$ of oriented, projectively
ASD   connections on $E$.

The first part ${\cal M}_I$ is just the moduli space ${\cal M}^{\rm
HE}(E,A_0)$ of oriented vortices, so it has a complex geometric
interpretation by Theorem
\ref{ORVOR}. The second part ${\cal M}_{II}$ can be   identified with
${\cal M}^{\rm HE}(E^\vee\otimes K_X,A_0^\vee\otimes [B_{\rm
can}^\vee]^{\otimes2})$  via the map
$$(A,\alpha)\mapsto (A^{\vee}\otimes B_{\rm
can}^\vee,\bar\alpha)\ .
$$
Therefore, we get the following complex geometric description of the
moduli space of non-abelian monopoles associated with the considered data:
\begin{thry}\label{MII}
The moduli space ${\cal M}$ of solutions of the monopole equations
associated with the data  $(\tau_{\rm can},B_{\rm can},\beta_{\rm
can},A_0)$ decomposes as a union of two closed subspaces
$${\cal M}={\cal M}_I\cup {\cal M}_{II}\ .
$$
These spaces can be identified with   moduli spaces of polystable
rank 2 oriented holomorphic pairs as follows:
$${\cal M}_I\simeq {\cal M}^{\rm pst}_{{\cal L}}(E)\ ,\ {\cal M}_{II}\simeq
{\cal M}^{\rm pst}_{{\cal L}^\vee\otimes{\cal K}_X^2}(E^\vee\otimes K_X)\ .
$$
The intersection ${\cal M}_I\cap {\cal M}_{II}$ is a Donaldson moduli space
of oriented projectively ASD unitary connections, and   can be identified
with the moduli space of polystable oriented holomorphic structures ${\cal
M}^{\rm pst}_{{\cal L}}(E)$.
\end{thry}

\begin{re} Using Theorem \ref{MII} one can describe explicitly moduli
spaces of non-abelian monopoles  as $S^1$-spaces. An explicit  example
in the  algebraic case -- which illustrates very clearly the cobordism
strategy -- is given in \cite{Te2}.

In the non-algebraic framework, one has to deal with a substantial
difficulty: the appearance of non-filtrable bundles (i.e. of bundles which
admit non proper subsheaf of lower rank).  For such bundles there exist no
general construction and classification methods.

Note however that this difficulty concerns only the ``Donaldson part" of
the moduli space (the intersection ${\cal M}_I\cap {\cal M}_{II}$)
because, for an oriented holomorphic pair $({\cal E},\varphi)$ with
$\varphi\ne 0$, the bundle ${\cal E}$ will be automatically filtrable.
\end{re}

\section{Appendix}

\subsection{Chern connections}\label{Chern}

Let $G$ be a complex Lie group and $Q$    a
principal $G$-bundle over
a  complex manifold $M$. By definition, an (almost)
holomorphic structure in
$  Q$ is an (almost)
complex structure in the total space of $Q$ which
makes the projection on $M$ (almost) holomorphic
and the action $G\times  Q\ra  Q$
holomorphic. In other words, an (almost) holomorphic structure in
the principal bundle $ Q$ is a $G$-invariant
(almost) complex structure on its total space which
agrees with the  complex structure induced by the complex structure in the
Lie
algebra
$\g$ on the vertical subbundle $V_Q\subset T_{Q}$,
and with
the complex structure in $M$ on the quotient bundle
$\qmod{T_{Q}}{V_{Q}}$.

A holomorphic structure on ${Q}$ can be
alternatively defined as the data of a maximal system of
trivializations   with holomorphic
transition functions.\\

If $K$ is a maximal compact subgroup of a reductive group $G$   and $P\subset Q$
is a $K$-reduction of $Q$, then the {\it Chern correspondence}   identifies
the  space  $\bar {\cal
A}({ Q})$ of almost  holomorphic    structures on  the
bundle ${Q}$ with the space   ${\cal A}({P})$) of   connections  on ${P}$.

The correspondence is defined as follows:  if $J$ is an almost complex
structure on the $G$-bundle $Q$, then the horizontal spaces of the associated
connection $A\in{\cal A}(P)$ are given by:
$$T_p^{\rm horiz}(P)=T_p(P)\cap J(T_p(P))\ ,\ p\in P\ .
$$
The tangent space $T_p(Q)$ decomposes as
$$T_p(Q)=T_p(P)\oplus T_p(P)^\bot=T_p^{\rm horiz}(P)\oplus T_p^{\rm
vert}(P)\oplus T_p(P)^\bot\ ,
$$
where $T_p(P)^\bot=\{a^\#_p|\ a\in\i\kg\}$.

The Chern correspondence identifies the space ${\cal H}(Q)$ of  holomorphic
structures on
$Q$ with the space
$${\cal A}^{1,1}(P):=\{A\in{\cal A}(P)|\ F_A\in{\cal A}^{1,1}(\ad(P))\}$$ of
{\it integrable} connections on
$P$.

  We will denote by
$J^A$ the (almost) holomorphic structure corresponding to the connection
$A$ and by $A_J$ (or, more precisely $A_{P,J}$) the connection
corresponding to $J$.

The aim of  this section is to study the dependence of this connection  and its
curvature on $P$, when $J$ is fixed.

\begin{pr}\label{ChernConn} The vertical component of a vector $w\in
T_p(P)$ with
respect to the Chern connection $A_{P,J}$ is given by  the formula
$${\rm vert}(w)=-J\left[{\rm
pr}_{T_p(P)^\bot}(J(w))\right]\ .
$$
\end{pr}
\pf  One has
$$J(T_p(P)^\bot)=T_p^{\rm vert}(P)\ ,\ J(T_p^{\rm horiz}(P))=T_p^{\rm
horiz}(P)\
,
$$
hence
$$w={\rm horiz}(w)+{\rm vert}(w)\ ;\ J(w)=J{\rm horiz}(w)+J{\rm vert}(w)\ ,
$$
so that
$$J{\rm vert}(w)={\rm pr}_{T_p(P)^\bot}(J(w))\ ;\ {\rm vert}(w)=-J\left[{\rm
pr}_{T_p(P)^\bot}(J(w))\right]\ .
$$
\qed

Let $P^0$ be a fixed $K$-reduction of $Q$,
and let $A^0$ be the Chern connection of the pair $( P^0,J)$. One can
write
\begin{equation}\label{newbundle}
  P=e^{-\frac{s}{2}} ( P^0)\ ,
\end{equation}
where $s\in A^0( P^0\times_{\ad}i\kg)$.

\begin{pr}\label{newconnection}
Denote by $A^s$ the Chern connection   of the pair $(P,J)$.
Regarding both
connections $A^s$ and the initial connection
$  A^0:=A_{  P^0, J}$  as connections in $Q$,  the relation
between these two connections is
$$  A^s= A^0+h^{-1}\partial_0 h\ ,
$$
where $h:=e^s$, and the term  $h^{-1}\partial_0 h$ on the right is
defined in the
following way: regard
$h$ as a $  G$-equivariant map $\chi:Q\ra G$, and define the
$\ad$-tensorial $(1,0)$ form $\chi^{-1}\partial_0\chi$ on $  Q$ by setting
$$\chi^{-1}\partial_0\chi(v)=\chi^{-1}\partial \chi({\rm horiz}_{
A^0}(v))=[L_{\chi(q)^{-1}}]_*\partial
\chi({\rm horiz}_{ A^0}(v))
$$
for all  $q\in  Q,  v\in T_q(  Q)$. This tensorial form can be identified with
a form in $h^{-1}\partial_0 h\in A^{1,0}(X,\ad (  Q))$.
\end{pr}

\pf  Put $k:=e^{-\frac{s}{2}}$. Choose $p\in  P^0$, and let
$q:=k(p)=p\kappa(p)\in\
P^s$.   We want to compute the difference
$\alpha:=\theta_{  A^s}-\theta_{  A^0}$ between the two connection forms in the
point $q$.

Fix an $A^0$-horizontal vector $w\in T_{p}(P^0)$.  The vector
$$[R_{\kappa(p)}]_*(w)\in
T_q(R_{\kappa(p)}(P^0))\subset T_q(Q)$$
 will be horizontal with respect to the same
connection regarded as a connection on $Q$.

On the other hand,
$$k_*(w)=[R_{\kappa(p)}]_*(w)+ q\ \kappa(p)^{-1}d\kappa(w)\in T_q(P^s)\ .
$$
Since $J(w)$  belongs to $T_p(P^0)$ too, one has
$$k_*(Jw)=[R_{\kappa(p)}]_*(Jw)+ q\ \kappa(p)^{-1}d\kappa(Jw)\in T_q(P^s)\ ,
$$
hence
$$[R_{\kappa(p)}]_*(Jw)\equiv - q\ \kappa(p)^{-1}d\kappa(Jw)\ {\rm mod}\
T_q(P^s)
$$
and
$$(\theta_{A^s}-\theta_{A^0})(k_*(w))=\theta_{A^s}([R_{\kappa(p)}]_*(w)+ q\
\kappa(p)^{-1}d\kappa(w))- \kappa(p)^{-1}d \kappa(w)\ .
$$
But
$${\rm vert}_{A^s}([R_{\kappa(p)}]_*(w)+ q\
\kappa(p)^{-1}d\kappa(w))=$$
$$=-J\ {\rm pr}_{T_q(P^s)^\bot}(([R_{\kappa(p)}]_*(Jw)+ Jq\
\kappa(p)^{-1}d\kappa(w))=
$$
$$=-J\ {\rm pr}_{T_q(P^s)^\bot}(- q\ \kappa(p)^{-1}d\kappa(Jw)+ Jq\
\kappa(p)^{-1}d\kappa(w))=
$$
$$=q\{-J_\g\ {\rm pr}_{\kg^\bot}(-  \ \kappa(p)^{-1}d\kappa(Jw)+   \
J_\g(\kappa(p)^{-1}d\kappa(w))\}
$$
Therefore, since $\theta_{A^s}-\theta_{A^0}$ is a tensorial form of
type $\ad$ and
$[R_{\kappa(p)}]_*(w)- k_*(w)$ is a vertical vector, one gets
$$\ad_{\kappa(p)^{-1}}(\theta_{A^s}-\theta_{A^0})(w)=
(\theta_{A^s}-\theta_{A^0})([R_{\kappa(p)}]_*(w))=(\theta_{A^s}-\theta
_{A^0})(k_*(w))=$$
$$=-i\ {\rm pr}_{\kg^\bot}\left[-  \
\kappa(p)^{-1}d\kappa(Jw)+   \ i(\kappa(p)^{-1}d\kappa(w)\right]-
\kappa(p)^{-1}d\kappa(w)=
$$
$$=-i\ {\rm pr}_{i\kg}\left[i(i
\kappa(p)^{-1}d\kappa(Jw)+   \  \kappa(p)^{-1}d\kappa(w))\right]-
\kappa(p)^{-1}d\kappa(w)=
$$
$$=-i\ {\rm pr}_{i\kg}\left[2i\kappa(p)^{-1}\bar\partial\kappa(w)\right]-
\kappa(p)^{-1}d\kappa(w)=
$$
$$=-i\ {\rm pr}_{i\kg}\left[2i\kappa(p)^{-1}\bar\partial\kappa(w)\right]-
\kappa(p)^{-1}d\kappa(w)=
$$
$$=-i
\left[i\kappa(p)^{-1}\bar\partial\kappa(w) -
i\partial \kappa(w)\kappa(p)^{-1}\right]-
\kappa(p)^{-1}d\kappa(w)=
$$
$$=  - \partial \kappa(w)\kappa(p)^{-1} -
\kappa(p)^{-1}\partial \kappa(w) \ ,
$$
so that
$$(\theta_{A^s}-\theta_{A^0})(w)=\ad_{\kappa(p)}[ - \partial
\kappa(w)\kappa(p)^{-1}
-
\kappa(p)^{-1}\partial \kappa(w)]=$$
$$[ - \kappa(p)\partial \kappa(w)\kappa(p)^{-2} -
\partial \kappa(w)\kappa(p)^{-1}]\ .
$$
But
$$\chi^{-1}d\chi=\kappa^2\partial(\kappa^{-2})=  -
\kappa (
\partial \kappa)\kappa^{-2} - (
\partial \kappa)\kappa^{-1}\ .
$$
\qed

\begin{co}\label{newcurvature} For the curvature of
$A^s:=A_{e^{-\frac{s}{2}}(P^0),J}$
one gets the formula
$$F_{ A^s}=F_{A^0}+\bar\partial( h^{-1}\partial_0 h)\ ,
$$
where $\bar\partial$ stands for the Dolbeault operator associated with the
fixed holomorphic structure $J$ on $Q$.
\end{co}

Sometimes it is more convenient to work with connections on the fixed
$K$-principal
bundle $P^0$.  Therefore, we will also consider the connection
$$A_s:=[e^{-\frac{s}{2}}]^*(A^s)\in{\cal A}(P^0)\ .
$$
\begin{pr}\label{NewConn} The connection form and the curvature form of the
connection $A_s$ are
given by the formulae
$$A_s-A^0=k^{-1}\bar\partial_0k- (\partial_0
k)k^{-1}=l^{-1}(\partial_0l)-(\bar\partial_0 l)l^{-1}\ ,
$$
$$F_{A_s}= \ad_{l} (F_{A^s})=\ad_{l}(F_{A^0}+\bar\partial(
h^{-1}\partial_0 h)) \ ,
$$
where $k:=e^{-s}$ and $l:=k^{-1}=\sqrt{h}$.
\end{pr}
\pf
$$A_s-A^0=k^{-1}(d_{A^s}k)+h^{-1}\partial_0 h=$$
$$k^{-1}(d_0k)+k^{-1}
[-k(\partial_0 k)k^{-2} - (\partial_0k)k^{-1},k]- k(\partial_0 k)k^{-2} - (
\partial_0k)k^{-1}=
$$
$$
=k^{-1}(d_0k)- (\partial_0 k)k^{-1} - k^{-1}(
\partial_0k)+ k(\partial_0 k)k^{-2}+( \partial_0k)k^{-1}
- k(\partial_0 k)k^{-2} - (
\partial_0k)k^{-1}
$$
$$=k^{-1}(\partial_0k)+k^{-1}(\bar\partial_0k)- (\partial_0 k)k^{-1} - k^{-1}(
\partial_0k)
=k^{-1}\bar\partial_0k- (\partial_0 k)k^{-1}\ .
$$
\qed

\subsection{The  orbits of the adjoint action. Sections in the adjoint bundle}
\label{AdjAct}

Let   $K$ be a compact Lie group and  $T\subset
K$   a fixed maximal torus.

We need first to understand the structure of the quotient $\kg/K$, where $K$
acts on its Lie algebra $\kg$ via the adjoint representation.

 First of all
note that any element $\varphi\in\kg$ is conjugate to an element in $\tg$.
Moreover, if two elements $\varphi_1$, $\varphi_2\in\tg$ are conjugate, i.
e.
$\ad_k(\varphi_1)=\varphi_2$ for some $k\in K$, then $\ad_k(T)$ and $T$ are
both maximal tori of $K$ which are contained in the connected centralizer
$Z_K(\varphi_2)$. Therefore, these tori are conjugate in $Z_K(\varphi_2)$,
hence there exists $u\in Z_K(\varphi_2)$ such that $\ad_u(\ad_k(T))=T$. This
shows that $\kappa:=uk$ belongs to the normalizer $N_K(T)$ of $T$ in $K$, hence
$\ad_\kappa(\varphi_1)=\ad_u(\varphi_2)=\varphi_2$, so that $\varphi_1$,
$\varphi_2$ are congruent modulo $N_K(T)$. Concluding, we get an identification
$$\qmod{\kg}{K}=\qmod{\tg}{N_K(T)}=\qmod{\tg}{W_K(T)}\ .
$$
When $K$ is connected  the quotient  $\qmod{\tg}{W_K(T)}$ can be
identified with
any closed Weyl chamber $C\subset \tg$. Indeed, it suffices to note that two distinct
elements on the boundary of $C$ cannot be congruent modulo $W_K(T)$.

To prove this, let $x,\ y\in  C\setminus \cringle{C}$ be two points which are congruent
modulo $W_K(T)$. Let
$w\in W_K(T)$ such that $w(x)=y$. Then $y\in C\cap w(C)$.  Let $c$ be the
common face
  of $C$ and $w(C)$ of minimal dimension which contains $y$.

Since the chambers $C$ and $w(C)$ have a common face $c$, one can  find a
sequence of chambers $(C_i)_{1\leq i\leq k}$ such that $C_1=C$, $C_k=w(C)$ and
$C_i$, $C_{i+1}$ have a common codimension 1 face which contains $c$.  It follows
that $C_i$ is mapped  to $C_{i+1}$ by the reflexion with respect to the hyperplane
$H_i$ generated by this common codimension 1 face.   Therefore $C$ is mapped to
$w(C)$ by a composition of reflexions with respect to hyperplanes which contain $c$.
Since the operation of the Weyl group is free transitive on the set of chambers, it
follows that $w$ coincides with this composition, hence it leaves $y$ fixed.

Therefore,
\begin{re}\label{ConjClasses} Suppose that $K$ is connected, and let $C\subset \tg$
be a  closed Weyl
chamber. The natural map
$$C\map \qmod{\kg}{K}
$$
is a homeomorphism.
\end{re}

Let $K$ be a connected compact Lie group  and $P$   a principal
$K$-bundle over an
arbitrary manifold $B$. We fix  a closed Weyl
chamber $C\subset \tg$.  $C$ decomposes as a disjoint union of convex
sets, the open faces
of $C$. Each open face is the interior of a face of $C$ (in the linear
subspace generated by this face).
\begin{lm}\label{ConjClassMap}  Let $\varphi\in
\Gamma(X,\ad(P))$ be a section regarded as a $K$-equivariant map $P\rightarrow \kg$.
Suppose that the map
$$B\textmap{[\varphi]} \kg/K\stackrel{\simeq}{\map}C$$
 takes values in a
fixed open face
$c$ of $C$.

Then:
\begin{enumerate}
\item The space $\varphi^{-1}(c)/T\subset P/T$
is a locally trivial fibre bundle with standard fibre  $Z_K(c)/T$ over $B$.
\item  Let $A$ be any connection on $P$. The map $B\textmap{[\varphi]} \kg/K$ is
locally constant if and only if the  projection of  $d_{A,x}(\varphi)$ on
$z_{\kg_x}(\varphi_x)\subset \kg_x:=P_x\times_\ad \kg$  vanishes for every $x\in
B$.
\end{enumerate}
\end{lm}
\pf 1. Note first that centralizer $Z_K(\gamma)$ of any element
$\gamma\in c$ is the
same,  hence
$Z_K(\gamma)=Z_K(c)$ for every $\gamma\in c$. Indeed, since $K$ is
connected, any
centralizer of the form
$Z_K(u)$ with
$u\in\kg$  is connected, so it suffices to look at the Lie algebras of these
centralizers. But the Lie algebra $z_k(u)$ has the following simple
description in
terms of the root set $R\subset\Hom_\R(\tg,i\R)$:
$$z_k(u)=\left[\bigoplus\limits_{\alpha\in R,\
\alpha(u)=0}\kg_\alpha^\C\right]\cap\kg\ .
$$
Therefore it suffices to note that the map
$$\gamma\mapsto\{\alpha\in R|\ \alpha(\gamma)=0\} \subset R
$$
is constant on the open face $c$.

The map
$$\qmod{K}{Z_K(c)}\times c\map \kg\ ,\ ([k],\gamma)\mapsto \ad_k(\gamma)$$
is an embedding, so its image $\ad_K(c)$ is a submanifold of $\kg$.

 The   restriction of the map
$$\varphi:P\ra \ad_K(c) $$ to any fibre $P_x$ is
obviously a submersion, hence $\varphi$ is also a submersion. The manifold
$\varphi^{-1}(c)$  is obviously a principal $Z_K(c)$-bundle over $B$.
Therefore, the quotient
$\varphi^{-1}(c)/T$ is just the associated bundle with fibre $Z_K(c)/T$.

2.  Let $s\in \ad_K(c)$. The  Lie algebra $\kg$ decomposes as
\begin{equation}\label{deco}
\kg=z_\kg(s)\oplus[\kg,s]\ ,
\end{equation}
whereas the tangent space at $s$ of  the submanifold
$\ad_K(c)\subset \kg$  decomposes as
$$T_s(\ad_K(c))=[\kg,s]\oplus \left[z_\kg(s)\cap
T_s(\ad_K(c))\right]=T_s(\ad_k(s))\oplus \left[z_\kg(s)\cap
T_s(\ad_K(c))\right].
$$

Therefore, a map $\sigma:U\ra c$ on a manifold $U$ defines a locally constant map
$[\sigma]:U\ra
\kg/K$ if and only if ${\rm pr}_{z_\kg(\sigma_x)}(d_x\sigma)=0$ for all $x\in U$.

In our case, we deduce that the map $P\ra \kg/K$ defined by $\varphi:P\ra \kg$  is
locally constant if and only if ${\rm pr}_{z_\kg(\varphi_p)}(d_p\varphi)=0$  for all
$p\in P$.  But, for a vertical tangent vector $v\in T_p(P)$ one obviously has
 $d_p\varphi(v)\in[\kg,\varphi_x]$, so the condition ${\rm
pr}_{z_\kg(\varphi_p)}(d_p\varphi)=0$ is in fact equivalent  with the condition ${\rm
pr}_{z_{\kg_x}(\varphi_x)}(d_{A,x}\varphi)=0$.
\qed

\begin{lm}\label{GIT} Let $K$ be a maximal compact subgroup of the reductive group
$G$. The natural map
$$\qmod{i\kg}{ K}\map  \qqmod{\g}{  G}
$$
is  injective.
\end{lm}

Here we denoted by $\qqmod{\g}{G}$ the GIT quotient of $\g$ with
respect to the
adjoint action of $ G$.  In other words,
$$ \qqmod{\g}{  G} ={\rm Spec}\{\C[\g]^{  G}\}\ .
$$
\pf The moment map for the $  K$-action on $\g$ is
$$m:\gamma\mapsto -\frac{i}{2}[\gamma,\gamma^*]\ ,
$$
where $(\cdot)^*$ is the conjugation with respect to the real structure
$\g=[i\kg]\otimes\C$.

It suffices to  note that $\kg$ is contained in $m^{-1}(0)$.
\qed

\begin{co}\label{Invariants} Let $\iota_1,\dots,\iota_k\in\C[\g]$ be
homogeneous
generators of the invariant algebra $ \C[\g]^{  G}$. Then  two
elements $u$, $v\in
i\kg$ are in the same $  K$-orbit if and only if
$\iota_i(u)=\iota_i(v)$ for $1\leq i\leq
k$.

In particular, a section $\sigma$ in the adjoint bundle $P\times_{\ad}
i\kg$ of a principal $K$-bundle is constant (respectively constant
almost everywhere) if and only if the  complex valued functions
$\iota_i(\sigma)$ are
constant  (respectively constant almost everywhere) for $1\leq i\leq k$.
\end{co}

The second statement in Lemma \ref{ConjClassMap} holds for general possibly
non-connected compact Lie groups. The converse statement is also true. More  precisely,
one has

\begin{pr}\label {ConstConjClassMap}
Let $K$ be a compact Lie group, $P$ a principal $K$-bundle, $A$ a connection on $P$
and   $\varphi\in
\Gamma(X,\ad(P))$ regarded as a $K$-equivariant map $P\rightarrow \kg$.

  The map
$B\textmap{[\varphi]}
\kg/K$ defined by $\varphi$ is locally constant if and only if the  projection
${\rm pr}_{z_{\kg_x}(\varphi_x)}(d_{A,x}\varphi)$
vanishes for every
$x\in B$.
\end{pr}
\pf

We suppose that the base $B$ is connected.

1.  {\it The implication:}
$$[\varphi]\ \hbox{is constant}\ \Rightarrow [{\rm
pr}_{z_{\kg_x}(\varphi_x)}(d_{A,x}\varphi)=0\ \forall x\in B]\ .$$

In the   case when $K$ is connected, one just takes the open chamber $c$ which contains
 the conjugacy class of $[\varphi]$ and applies  Lemma
\ref{ConjClassMap}.

The general case can be  reduced to the connected case   by taking a finite covering of  the
base $B$ on which $P$
  admits a
$K_e$-reduction $P_e\subset P$.  It suffices to note that, in general, a continuous map
 $U\map\kg/{K_e}$ on a connected manifold $U$ is  constant  if and only if the
induced map  $U\map\kg/{K}$ is   constant.\\

2. {\it The implication:}
$$[{\rm
pr}_{z_{\kg_x}(\varphi_x)}(d_{A,x}\varphi)=0\ \forall x\in B]\Rightarrow   [\varphi]
\ \hbox{is constant}\ .$$

Let $j_i:\g^{d_i}\ra\C$ be the $\ad_{G}$-invariant symmetric
multilinear map
which corresponds to $\iota_i$ (see Corollary \ref{Invariants}).

One has
$$d [j_i(\varphi,\dots,\varphi)]=
d_i j_i(d_A\varphi,\varphi,\dots,\varphi)\ .
$$

On the other hand, taking into account that $j_i$ is $\ad_{
G}$-invariant, it follows easily that, for any $v\in
\g$, the linear
functional  $j_i(\cdot,v,\dots,v)$ vanishes on
$[v,\kg]$. The point is that, for  $v\in \kg$, one has an orthogonal
decomposition
$$\kg=[v,\kg]\oplus z_\kg(v)\ .
$$

Therefore, one
gets
$$d [j_i(\varphi,\dots,\varphi)]=d_i j_i(d_A\varphi,\varphi,\dots,
\varphi) =d_i j_i({\rm
pr}_{z_\kg(v)}(d_A(\varphi)),\varphi,\dots,
\varphi)=0 \ .
$$
  The claim follows now from
Corollary \ref{Invariants}.
\qed

In section \ref{ContMeth} we  have used the following technical result:

\begin{pr}\label{Differential} Let $K$ be a maximal compact subgroup of a complex
reductive group
$G$ and $s, \dot s\in i\kg$. Using the decomposition (\ref{deco}), write
$$\dot s= \lambda+[  k,s]\ ,
$$
with $\lambda\in iz_\kg(s)\ ,\ k\in z_\kg(s)^\bot$.
Then
\begin{enumerate}
\item   $$ (d_s\exp)(\dot s) e^{-s}=   \lambda+( k-\ad_{e^s}k) \ .$$
\item  $$\langle \dot s, p_{i\kg}\left\{ (d_s\exp)(\dot
s) e^{-s}\right\}\rangle\geq \|\dot s\|^2 \ .$$
\item
$$\left\langle \left[ p_{\kg}\left\{ (d_s\exp)(\dot
s) e^{-s}\right\},s\right], p_{i\kg}\left\{ (d_s\exp)(\dot
s) e^{-s}\right\}\right\rangle\leq 0
$$
with equality if and only if $\dot s\in i z_\kg(s)$.
\end{enumerate}
\end{pr}
\pf

Put
$$\dot s= \lambda+[  k,s]=\lambda+\frac{d}{dt}(\ad_{e^{t  k}}(s))\ ,\
\lambda\in iz_\kg(s)\ ,\ k\in z_\kg(s)^\bot\ .
$$
One has
$$(d_s\exp)(\dot s)=(d_s\exp)(\lambda)+(d_s \exp)([  k,s])  =
\lambda e^s+\frac{d}{dt}|_{t=0}
\ \Ad_{e^{tk}}\ e^s=$$
$$[ \lambda+( k-\ad_{e^s}  k)]e^s \ ,
$$
$$
$$
so
$$\sigma:=(d_s\exp)(\dot s)e^{-s}=   \lambda+( k-\ad_{e^s}
k)\ .
$$
Let  $R_s$ denote  the  set of
eigenvalues of the endomorphism $[s,\cdot]$ on $\g$. Decompose
$$  k=\sum\limits_{\rho\in R_s\setminus\{0\}}  k_\rho\ ,\ k_\rho\in \g_\rho\ ,
$$
with $ \bar  k_{\rho}=-k_{-\rho}$ (since $k\in \kg$). Then
\begin{equation}\label{RootDec}
[s,k]=\sum\limits_{\rho\in R_s\setminus\{0\}} \rho\ k_\rho\ ,\
\ad_{e^s}  k=\sum\limits_{\rho\in R_s\setminus\{0\}} e^\rho k_\rho\ ,  $$
$$\dot s=\lambda-\sum_{\rho\in R_s\setminus\{0\}} \rho k_\rho\ ,\
\sigma=\lambda+\sum_{\rho\in R_s\setminus\{0\}}(1-e^\rho)k_\rho\ .
\end{equation}
We get
$$ p_{i\kg}(\ad_{e^s}  k) = \frac{1}{2}(\ad_{e^s}  k+\overline{\ad_{e^s}k })=
\frac{1}{2}\sum_{\rho\in R_s\setminus\{0\}} (e^\rho-e^{-\rho})k_\rho\ ,$$
$$p_\kg(\ad_{e^s}  k)=
\frac{1}{2}(\ad_{e^s}  k -\overline{\ad_{e^s}k })=
\frac{1}{2}\sum_{\rho\in R_s\setminus\{0\}} (e^\rho+e^{-\rho})k_\rho\ ,
$$
$$p_\kg(k-\ad_{e^s}  k)=\frac{1}{2}\sum_{\rho\in
R_s\setminus\{0\}}(2-e^\rho-e^{-\rho})k_\rho\ ,
$$
$$[p_\kg(\sigma),s]= \frac{1}{2}\sum_{\rho\in
R_s\setminus\{0\}}
\rho(e^\rho+e^{-\rho}-2)k_\rho=\frac{1}{2}\sum_{\rho\in R_s\setminus\{0\}}
\rho e^{-\rho} (e^\rho-1)^2k_\rho\ ,
$$
$$p_{i\kg}(\sigma)=\lambda+\frac{1}{2}\sum_{\rho\in
R_s\setminus\{0\}}(-e^\rho+e^{-\rho})k_\rho\ .
$$
Therefore,
$$\langle \dot s, p_{i\kg}(\sigma) \rangle=\|\lambda\|^2+
\frac{1}{2}\sum_{\rho\in
R_s\setminus\{0\}}\rho(e^\rho-e^{-\rho})\|k_\rho\|^2\geq
$$
$$\geq \|\lambda\|^2+
 \sum_{\rho\in
R_s\setminus\{0\}}\rho^2 \|k_\rho\|^2=\|\dot s\|^2\ ,
$$
because the inequality $|e^x-e^{-x}|\geq 2|x|$  holds on $\R$.

For the third  inequality we compute
\begin{equation}\label{LongFormula}
\langle [p_\kg(\sigma),s], p_{i\kg}(\sigma)\rangle =$$
$$
=-\sum_{\rho\in
R_s\setminus\{0\}}\frac{1}{4}\rho(e^{\rho}-e^{-\rho})e^{-\rho}(e^\rho-1)^2\|k_\rho\|^2\leq
0\ .
\end{equation}
Equality  occurs s if and only if $\sum_\rho \rho\  k_\rho=0$.
\qed

\subsection{Local maximal torus reductions of a
$K$-bundle}

We quote from [LT] the following Lemma (see [LT], section 7.4)
\begin{lm}\label{filtration} Let $X$ be a manifold, and
$$X = G_r \supset
G_{r-1} \supset \ldots \supset G_1$$
 a filtration by closed
subsets. Define $F_k:= {G_k}\setminus{G_{k-1}}\ ,\ k =
2,\ldots,r\ .$ Then
$$
    W := \bigcup\limits_{k=2}^r {\inw{F}}_k
$$
is open and dense in $X$.

\end{lm}

\begin{pr}\label{LocDiag} For any $f\in \Gamma(B,\ad(P))$ there is
an open dense subset $W \subset X$ such that for  every
$x \in W$ there exists an open face $c_x\subset C$,  an open
neighborhood $U_x$ of
$x$ in $B$ and a
$T$-reduction
$\Pi_x\subset P|_{U_x}$ of the restriction $P|_{U_x}$ such that
$f|_{U_x}$  is defined by  a  smooth map
$$\lambda
\in  {\cal C}^\infty(U_x,c_x)\subset {\cal
C}^\infty(U_x,\tg)=\Gamma(U_x,\ad(\Pi_x))\ .$$
\end{pr}
\pf

Regard $f$ as a $K$-equivariant map $P\ra \kg$, and let $[f]:B\ra
\kg/K=C$ be the
induced $\kg/K$-valued  map on $B$. Consider the skeleton filtration
$$C=C_d\supset
C_{d-1}\supset\dots C_1\supset \dots \supset C_0$$
of $C$, where $C_i$ is the union of all open faces  of dimension less
or equal $i$.  We get
a filtration by closed sets
$$B=B_d\supset
B_{d-1}\supset\dots B_1\supset \dots \supset B_0
$$
of $B$, where $B_i:=[f]^{-1}(C_i)$.

Consider now the locally closed subsets
$$F_i:=B_{i}\setminus
B_{i-1}=[f]^{-1}(C_{i}\setminus
C_{i-1})=[f]^{-1}\left[\coprod_{\matrix{\scriptstyle
c\subset C\
\hbox{\scriptsize open face}\vspace{-1.5mm}\cr\scriptstyle\dim(c)=i}}
c\right]\ .
$$

By Lemma \ref{filtration}, it follows that $W:=\union \cringle{F_i}$ is
dense in $B$.
For any $x\in W$ let $i_x\in\{0,\dots,d\}$ such that $x\in
\cringle{F}_{i_x}$ and
let $U_x$ be an arbitrary contractible open neighborhood of $x$ in
$\cringle{F}_{i_x}$. Since $U_x$ is connected, $[f](U_x)$ is
contained in a single
open face $c_x$ of $C$. Moreover, the fibre bundle
$\left[f|_{U_x}^{-1}(c_x)\right]/T$ over the contractible manifold
$U_x$ is trivial,
hence  one can find a section
$\sigma_x$ in this bundle.  Regarding $\sigma_x$ as a section in
$\left[P|_{U_x}\right]/T$, we get a $T$-reduction of $P|_{U_x}$ on which $f$ is
given by a smooth $c_x$-valued function.
\qed
\begin{re}\label{GenCase}  Proposition \ref{LocDiag} is also true for
non-connected
Lie groups
$K$. The proof reduces to the connected case using a finite cover of $B$, on
which $P$ can be reduced to the connected component $K_e\subset K$.
\end{re}

\subsection{Decomposing a connection with respect to a maximal torus reduction.
Applications} \label{UsefulFormulae}

Let $P$ be a principal $K$-bundle over a manifold $B$ and let
$\iota:\Pi\hookrightarrow P$ be a
$T$-reduction of $P$.  Consider the root set $R\subset\Hom(\tg,i\R)$
and the root
decomposition
$$\g=\bigoplus\limits_{\rho\in R} \g_\rho
$$
of $\g:=\kg^\C$.

Since this decomposition is $T$-invariant, and we have fixed a
$T$-reduction of $P$, one gets a decomposition of the complexified
adjoint bundle
$\ad^\C(P)$:
\begin{equation}\label{RootDecAd}\ad^\C(P)=\bigoplus_{\rho\in
R}\Pi\times_{\ad}\g_\rho.
\end{equation}
Denote by $\tg^\bot$ the natural complement of the Lie algebra $\tg$
in $\kg$, i. e.
$$\tg^\bot=\left[\bigoplus_{\rho\ne 0}\kg_\rho^\C\right]\cap\kg\ .
$$

Let $A$ be a connection on $P$, and $\omega_A\in A^1(P,\kg)$ its connection
form.  Let $A_\Pi$ be the induced connection on $\Pi$, and
$\iota_*({A_\Pi})$ the push-forward connection to $P$. Write
$$\omega_A=\omega_{\iota_*(A_\Pi})+a\ ,
$$
where  $a\in
A^1(\Pi\times_\ad\tg^\bot)$ is the second fundamental form of the
subbundle $\Pi$
with respect to the connection $A$.

Suppose now that the basis $B$ is a complex manifold. Then one can write
$a=a^{10}+a^{01}$, where each term decomposes as
$$a^{10}=\sum_{\rho\ne 0}a^{10}_\rho\ ,\ a^{10}_\rho\in
A^{10}(\Pi\times_{\ad}\g_\rho)\ ,\  a^{01}=\sum_{\rho\ne 0}a^{01}_\rho\ , \
a^{01}_\rho\in A^{01}(\Pi\times_{\ad}\g_\rho).
$$
With respect to the usual real structure $\g=[i\kg]^\C$ of
$\g$, one has
$$\overline{a^{10}_\rho}=-a^{01}_{-\rho}\ .
$$

Let $\varphi\in A^0(\ad(P))$ be a smooth section. Suppose now that,
with respect to
the reduction $\Pi$, the section
$\varphi $ is induced by a smooth map
$\lambda\in{\cal C}^\infty(B,c)$ where $c$ is a fixed open face of
the Weyl chamber
$C$. We know, by Proposition \ref{LocDiag}, that  locally   we can
always come to
this situation.

Via the decomposition (\ref{RootDecAd}), and taking into account that
$[a,\lambda]=-[\lambda,a]$, we  have the following formulae for
$\partial_A \varphi$, $\bar\partial_A \varphi$:
\begin{equation}\label{ds}
\partial_A \varphi=  \partial\lambda-\sum_{\rho\ne 0}
\rho(\lambda) a^{10}_\rho\ ,\ \bar \partial_A \varphi=
\bar\partial\lambda-\sum_{\rho\ne 0}
\rho(\lambda) a^{01}_\rho=   \bar\partial\lambda+\sum_{\rho\ne 0} \rho(\lambda)
\overline{ a^{10}_{\rho}}\ . \\
\end{equation}

We will also need a formula for $e^{-s}\partial_A (e^s)$, where $s\in A^0(i\ad(P))$.
Regard $s$ (and $e^s$) as an  $\ad$-equivariant (respectively $\Ad$-equivariant)
$\g$-valued (respectively $G$-valued) function on
$P$.  Let be $x\in X$,  $p\in \Pi_x$ and
$v\in T_x(X)$ and $w\in T_p(\Pi)$ the $A_\Pi$ - horizontal lift of $v$.

The $A$-horizontal component of  $w$ is

$$w^A=w-[\omega_A(w)]^\#_p=w-a(w)^\#_p\ .
$$
One has
\begin{equation}
e^{-s(p)}d
(e^s)(w^A)=e^{-s(p)}\left[d(e^s)(w)-\frac{d}{dt}|_{t=0}
\left(e^{s(p\exp(ta(w))}\right)
\right]=
$$
$$
=e^{-s(p)}\left[d(e^s)(w)-\frac{d}{dt}|_{t=0}
\left(e^{\ad(\exp(ta(w)))^{-1}s(p)}\right)
\right]=
$$
$$=e^{-s(p)}\left[d(e^s)(w)-\frac{d}{dt}|_{t=0}\left(\Ad(\exp(ta(w)))^{-1}
(e^{s(p)})\right)
\right]=
$$
$$
=e^{-s(p)}\left[d(e^s)(w)- [-a(w)e^{s(p)}+e^{s(p)}a(w)]
\right]\ .
\end{equation}
Denote now by $\lambda$ the restriction  $\lambda:=s|_\Pi$.
We get
$$e^{-s(p)}d (e^s)(w^A)  =d\lambda(w) + \sum_{\rho\ne
0}(e^{-\rho(\lambda)}-1)a_\rho(w)\ ,
$$
hence
\begin{equation}\label{d(exp(s))}
e^{-s}d_A(e^s)= d\lambda+\sum_{\rho\ne
0}(e^{-\rho(\lambda)}-1)a_\rho\ ,
\end{equation}
and similarly
\begin{equation}\label{newd(exp(s))}
d_A(e^s)e^{-s}= d\lambda+\sum_{\rho\ne
0}(1-e^{\rho(\lambda)})a_\rho\ .
\end{equation}

In section \ref{Estimates}  we need the following important technical results.

Let $E$ be a Hermitian space and $f:\R\ra\R$  a real function. Using the spectral
decomposition of Hermitian matrices, one can extend $f$ to a map (denoted by the
same symbol to save on notations) $\Herm(E)\ra\Herm(A)$ by putting
$$f(\sum_i\lambda_i{\rm pr}_{F_i}):=\sum f(\lambda_i){\rm pr}_{F_i}\
$$
for every orthogonal direct sum decomposition of $E$.

It is known (see for instance \cite{Bh}) that
\begin{pr} The extension $\Herm(E)\ra\Herm(E)$ of a    continuous function
$f:\R\ra\R$ is continuous.
\end{pr}

In general, this extension inherits the regularity properties  of $f$.

\begin{re}\label{monotony}
\begin{enumerate}
\item If $|f|\leq |g|$ then
$$\|f(h)(v)\|\leq \|g(h)(v)\|\ \forall v\in E,\ \forall h\in\Herm(E)\  .
$$
\item More generally, if on a subset $A\subset\R$  we have $|f|_A|\leq |g|_A|$, then
$$ \|f(h)(v)\|\leq \|g(h)(v)\|\   \forall v\in E,\ \forall h\in\Herm(E)\  \hbox{with }
{\rm Spec}(h)\subset A\  .
$$
\end{enumerate}
\end{re}
Indeed, writing $h=\sum_i\lambda_i{\rm pr}_{F_i}$ as above, one gets
$$\|f(h)(v)\|=\sqrt{\sum_i f(\lambda_i)^2\|{\rm pr}_{F_i}(v)\|^2}\leq \sqrt{\sum_i
g(\lambda_i)^2\|{\rm pr}_{F_i}(v)\|^2}= \|g(h)(v)\|\ .
$$
\\

Denote by $\eta:\R\ra\R_{>0}$ the \ub{real}
\ub{analytic} function defined by
$$\eta(t)=\left\{
\begin{array}{ccc}
\sqrt{\frac{1-e^{-t}}{t}}&\rm if& t\ne 0\\
1&\rm if&t=0\ .
\end{array}
\right.
$$

\begin{pr}\label{identity}
Let ${\cal Q}$ be a holomorphic  principal  $K^\C$-bundle on a complex
manifold $X$, $P\subset {\cal Q}$ a $K$-reduction, denote by $A^0$ the
corresponding Chern connection, and let $d_0=\bar\partial+\partial_0$ be the
associated covariant derivative. Let $s\in A^0(P\times_{\ad}i\kg)$.
Then   one has
\begin{equation}
\left(i\Lambda_g\bar\partial(
e^{-s}\partial_0(e^s))),s \right)=\frac{1}{2}
P(|s|^2)+|\eta([s,\cdot])(\partial_0(s))|^2\ .
\end{equation}
\end{pr}
\pf   Let $T$ be a maximal torus of $K$ and $R\subset\Hom(i\tg,\R)$ be the root
set. By Proposition
\ref{LocDiag} it suffices to check the formula in the case   when $P$ has a
$T$-reduction $\Pi\subset P$ and
$s$ is given by a smooth function $\lambda:X\ra i\tg$.
By formulae (\ref{ds}) and
(\ref{d(exp(s))}) one has
\begin{equation}\label{partial0-s}\partial_0 s =
\partial\lambda-\sum_{\rho\in
R\setminus\{0\}}
\rho(\lambda) a^{10}_\rho\ ,\
    e^{-s}\partial_0(e^s)  = \partial \lambda+\sum_{\rho\in
R\setminus\{0\}}
(e^{-\rho(\lambda)}-1)  a^{10}_\rho\ .
\end{equation}
Therefore, using the Hermitian inner product:
\begin{equation}\label{compute}
\left(i\Lambda_g\bar\partial(
e^{-s}\partial_0(e^s))),s \right)= i\Lo\dbar\left(
{e^{-s}}\circ\partial_0(e^s),s\right) +
 i\Lo ({e^{-s}}\circ\partial_0(e^s),\partial_0s)=$$
$$=i\Lo\dbar\left(\partial
\lambda,\lambda\right) + i\Lo(\partial\lambda,\partial\lambda)
-\sum\limits_{\rho\in
R\setminus\{0\}} \rho(\lambda) (e^{-\rho(\lambda)} -
1)i\Lo (a^{10}_\rho,a^{10}_\rho)=
$$
$$(P(\lambda),\lambda)-\sum\limits_{\rho\in
R\setminus\{0\}} \rho(\lambda) (e^{-\rho(\lambda)} -
1)i\Lo (a^{10}_\rho,a^{10}_\rho)\ .
\end{equation}

  On the other hand,
\begin{equation}\label{ThirdIneq}
    P(\vert s{\vert^2}) =
i\Lo\dbar\partial(|\lambda|^2) =
2  \left(P(\lambda),\lambda) -
\vert\partial{\lambda}{\vert^2}\right) \leq
2(P(\lambda),\lambda)\ .
\end{equation}
Now it suffices to note that
$$
|\partial\lambda|^2-\sum\limits_{\rho\in
R\setminus\{0\}} \rho(\lambda)
(e^{-\rho(\lambda)} - 1)i\Lambda_g
(a^{10}_\rho,a^{10}_\rho)=|\eta([s,\cdot])(\partial_0(s))|^2\ .
$$
\qed

\subsection{Analytic results}\label{AnalyticResults}

\begin{lm}\label{analysis}
Let $E$, $F$ be  vector bundles with inner products on a measurable space $X$ with
finite measure. Let
$(u_n)_n$ be a sequence of $L^2$-sections on $E$ weakly convergent in $L^2$ to a
section $u$, and let $(v_n)_n$ be a sequence of $L^2$-sections on $F$ strongly
convergent  in $L^2$ to a section $v$. Suppose that $v_n$ converges almost
everywhere to $v$ and that $(|v_n|)_n$ is uniformly bounded by a positive constant
$M$.  Then it holds:
\begin{enumerate}
\item For every $L^2$-section $\psi$ in $(E\otimes F)^\vee$,
the sequence $(v_n\psi)_n$ of sections in  $E^\vee$   converges strongly in $L^2$ to
$v\psi$.
\item The sequence $(u_n\otimes v_n)_n$ converges weakly in $L^2$  to $(u\otimes
v)$.
\end{enumerate}

\end{lm}

\pf 1. This follows
from the Lebesgue dominated convergence theorem, since the sequence of
non-negative functions $(|v_n\psi- v\psi|^2)_n$ converges almost everywhere to 0 and is
bounded almost everywhere by the integrable function
$4M^2[\psi|^2$.

2. We have:
$$ \langle  u_n\otimes v_n ,\psi\rangle_{L^2}-\langle u\otimes
v,\psi\rangle_{L^2} = \langle  u_n, v_n\psi\rangle_{L^2}-\langle u,
v\psi\rangle_{L^2} =
$$
$$= \langle  u_n, v_n\psi-v\psi\rangle_{L^2} +\langle
u_n-u, v \psi\rangle_{L^2}\ . $$
The second term converges to 0  since, by hypothesis, $u_n\ra u$ weakly. The
first term tends to 0 because
$$|\langle  u_n, v_n\psi-v\psi\rangle_{L^2}|
\leq\|u_n\|_{L^2}\|v_n\psi-v\psi\|_{L^2}\ ,$$
$(u_n)_n$ is bounded\footnote{Any weakly convergent sequence in a Hilbert
space is bounded.} in $L^2$ and $\|v_n\psi-v\psi\|_{L^2}\ra 0$ by 1.
\qed

\begin{lm}\label{Fubini} Let  $D^{n_i}\subset \R^{n_i}$ be the standard
$n_i$-dimensional disk,
$i=1,\ 2$, and let  $\varphi\in L^p_1(D^{n_1}\times D^{n_2})$, $1\leq p<\infty$.
Then  for almost all
$y\in D^{n_2}$ the following holds:
\begin{enumerate}
\item  $\varphi|_{{D^{n_1}\times\{y\}}}$ belongs to $L^p_1(D^{n_1})$,
$\partial_x\varphi|_{D^{n_1}\times\{y\}}$ belongs to
$L^p(D^{n_1})$.
\item $d (\varphi|_{{D^{n_1}\times\{y\}}})=(\partial_x\varphi)|_{D^{n_1}\times\{y\}}$
in $L^p(D^{n_1})$.
\end{enumerate}
\end{lm}
\pf By   Fubini's theorem $\varphi|_{D^{n_1}\times\{y\}}$ and
$d\varphi|_{D^{n_1}\times\{y\}}$ belong  to $L^p(D^{n_1})$ for almost all $y\in
D^{n_2}$. In particular, $\partial_x\varphi|_{D^{n_1}\times\{y\}}$ belongs to
$L^p(D^{n_1})$ for almost all $y\in D^{n_2}$, too.

Let $\varphi_n$ be a sequence of smooth functions converging to $\varphi$ in the
$L^p_1$ - topology    as
$n\ra\infty$.  Therefore
$$\int\limits_{D^{n_1}\times D^{n_2}} |\varphi_n-\varphi|^p\ra 0\ ,\
\int\limits_{D^{n_1}\times D^{n_2}} |d\varphi_n-d\varphi|^p\ra 0\ .
$$
In particular,
$$\int\limits_{D^{n_1}\times D^{n_2}} |\partial_x\varphi_n-\partial_x\varphi|^p\ra 0\ .
$$

By  Fubini's theorem, it follows that the positive functions
$$y\mapsto \int\limits_{D^{n_1}\times\{y\}}   |\varphi_n-\varphi|^p\ ,\  y\mapsto
\int\limits_{D^{n_1}\times\{y\}}   |\partial_x\varphi_n-\partial_x\varphi|^p$$
converge to 0 with respect to the $L^1$ - norm,  so (taking a subsequence if necessary)
one can assume that they converge to 0 almost everywhere on $D^{n_2}$.  Therefore, for
almost every $y\in D^{n_2}$, one has
$$\varphi_n|_{D^{n_1}\times\{y\}}\textmap{L^p} \varphi|_{D^{n_1}\times\{y\}}\ ,\
(\partial_x\varphi_n)|_{D^{n_1}\times\{y\}}\textmap{L^p}
(\partial_x\varphi)|_{D^{n_1}\times\{y\}}\ .
$$

But, since the $\varphi_n$ are smooth, one has
$$(\partial_x\varphi_n)|_{D^{n_1}\times\{y\}}=
\partial_x[\varphi_n |_{D^{n_1}\times\{y\}}]\ .
$$
One the other hand,
$$  \partial_x[\varphi_n
|_{D^{n_1}\times\{y\}}]\map \partial_x[\varphi |_{D^{n_1}\times\{y\}}]\ \hbox{as
distribution}\ $$
 (because the distribution limit of $\varphi_n
|_{D^{n_1}\times\{y\}}$ is $\varphi |_{D^{n_1}\times\{y\}}$).  Therefore,
$$\varphi_n|_{D^{n_1}\times\{y\}}\textmap{L^p} \varphi|_{D^{n_1}\times\{y\}}\ ,\
\partial_x[\varphi_n |_{D^{n_1}\times\{y\}}]\textmap{L^p}
\partial_x[\varphi
|_{D^{n_1}\times\{y\}}]=$$
$$=(\partial_x\varphi)|_{D^{n_1}\times\{y\}}\ ,
$$
which shows that 1. and 2. hold for such $y$.
\qed

\subsection{Meromorphic  parabolic reductions defined by weakly holomorphic $L^2_1$ -
  sections}\label{MerRed}

 Let $M$ be a complex manifold. We begin with the following definition:
\begin{dt}\label{HolLines} A \ub{measurable}
map
$F:B^n\ra M$ is called    ``holomorphic on almost all lines" if
\begin{enumerate}
\item  $n=1$, $F$ coincides almost everywhere with a holomorphic map  $B\rightarrow
M$.
\item $n>1$, and for any holomorphic parameterization $D\times
\Delta\textmap{\varphi} U\subset B^n$
 with $D\subset\C$, $\Delta\subset\C^{n-1}$ the restriction $F\circ\varphi(\cdot,\zeta):D\ra
M$ is weakly holomorphic in the sense of 1., for almost all $\zeta\in \Delta$.
\end{enumerate}
\end{dt}

We   recall  the following fundamental results (see Shiffman
\cite{Sh1}, Uhlenbeck-Yau \cite{UY1}, \cite{UY2}):
\begin{thry}\label{HolLinesIsMer}
Any  measurable map $F:B^n\ra M$ from the ball $B^n\subset\C^n$   into a
projective  algebraic manifold $M$, which is holomorphic on almost all lines,
coincides almost everywhere with a meromorphic map.
\end{thry}

Theorem \ref{HolLinesIsMer}  is used combined with an important     regularity
result (see Uhlenbeck-Yau \cite{UY1}, \cite{UY2}) which applies to the case $n=1$.

In order to state this result we need   a brief preparation:

Let $Y$ be a compact manifold (possibly with boundary) and $N$ a closed submanifold
of $\R^m$. We put
$$L^2_1(Y,N):=\{\varphi\in L^2_1(Y,\R^m)|\ \varphi(y)\in N\
\hbox{for almost every}
\ y\in Y\}\ ,
$$
and we endow this set with the  topology induced from $L^2_1(Y,\R^m)$. The
 space
$L^2_1(Y,N)$ obviously depends on the  embedding
$Y\hookrightarrow
\R^m$.  When $Y$ is a surface, one has the following important density property (see
Schoen-Uhlenbeck \cite{SU}).
\begin{pr}\label{density} Suppose that $Y$ is a compact surface with possibly empty
${\cal C}^1$-boundary. Then ${\cal C}^\infty(Y,N)$ is dense in  $L^2_1(Y,N)$.
\end{pr}

This result has the following subtle consequence:
\begin{co} Suppose that $Y$ is a compact surface with possibly empty ${\cal
C}^1$-boundary and let $\varphi\in  L^2_1(Y,M)$. Then $d\varphi(T_y(Y))\subset
T_{\varphi(y)}(M)$ for almost all $y\in Y$.
\end{co}
\pf By  the density property  Proposition \ref{density}, one can find a sequence
$\varphi_j\in {\cal C}^\infty(Y,N)$ which converges in the $L^2_1(Y,\R^m)$ - topology
to
$\varphi$. There exists a subsequence $(\varphi_{j_k})_k$ such that
$\varphi_{j_k}\ra\varphi$ and $d\varphi_{j_k}\ra d\varphi$ almost everywhere as
$k\ra\infty$. It suffices to note that $T(N)$ is closed in $N\times\R^m$.
\qed

Now let $N\subset \R^m$ be a compact complex manifold embedded as a real
submanifold in $\R^m$, and let $J_N$ be the corresponding almost complex structure.
The inclusion
$T_x(N)\hookrightarrow\R^m$ of the tangent space at $x\in Y$ induces an embedding
$$T^{10}_x(N)\oplus T^{01}_x(N)=T^\C_x(N)\hookrightarrow\C^m\ .
$$
Therefore, one gets a filtration
$$0\subset T^{10}(N)\subset T^\C(N)\subset N\times\C^m
$$
of the trivial complex rank $m$ bundle over $N$.

Let $p:N\times\C^m\ra T^{10}(N)$ be any bundle projection which induces the standard
projection
$$T^\C(N)\map T^{10}(N)\ ,\ v\mapsto\frac{1}{2}(v-iJ_N(v))
$$
on $T^{\C}(N)$.

\begin{re}\label{IndHerm}
In the special case when the induced Riemannian metric on the
submanifold $N\subset\R^m$ is    Hermitian, one can just take the orthogonal projection
with respect to the standard Hermitian structure on $N\times\C^m$  because, in this case,
the direct sum
$T^{10}(N)\oplus T^{01}(N)$ will be orthogonal with respect to this standard Hermitian
structure.
\end{re}

We can now state
\begin{thry} \label{Reg} \cite{UY1}, \cite{UY2} Let $B\subset \C$ be the standard ball and let
$N$ be a compact K\"ahler manifold embedded as a real submanifold in $\R^N$.   Let $f:B\ra
N$ be a map with the following properties:
\begin{enumerate}
\item   $f\in L^2_1(B,N)$.
\item $f$ is weakly holomorphic, i. e. one of the following two equivalent conditions holds:
\begin{enumerate}
\item For almost all $b\in B$ the tangent map $T_{b}(B)\ra T_{f(b)}(N)$ is $\C$-linear.
\item For almost all $b\in B$ it holds $p_{f(b)}\circ(\bar\partial_bf)=0$, where
$\bar\partial f$ stands for the $\bar\partial$-derivative $T^{01}_b(B)\ra\C^m$ of
$f:B\ra\C^m$.
\end{enumerate}
\end{enumerate}
Then $f$ coincides almost everywhere with a holomorphic map.
\end{thry}

Combining Theorem \ref{HolLinesIsMer}   with the regularity Theorem \ref{Reg}, one gets
the following important
\begin{thry}\label{meromorphy} Let $M \subset\R^m$  be a complex projective algebraic
manifold embedded as a real submanifold of $\R^m$, and let $p:M\times\C^m\ra
T^{10}(M)$ be a bundle projection which agrees with the standard projection
$T^{\C}(M)\ra T^{10}(M)$ on the subbundle $T^{\C}(M)\subset M\times\C^m$.  Then
any map $f\in L^2_1(B^n,M)$ satisfying the equation
$$p\circ\bar\partial f=0
$$
almost everywhere  on $B^n$,  coincides almost everywhere with a meromorphic map.
\end{thry}
\pf This follows easily from the Theorem \ref{HolLinesIsMer}, the regularity Theorem
\ref{Reg} and   Lemma \ref{Fubini}.
\qed
\\ \\

Let $K$ be a maximal compact subgroup of the reductive group $G$.  We will apply these
results to an algebraic manifold of the form
$$M=\qmod{G}{G(s_0)} \ ,
$$
where $G(s_0)$ is the parabolic subgroup  defined by  an element $s_0\in  i\kg$.  We
will need the following simple result.
\begin{pr}\label{HolTangSp} i) The map
$$\ad_K(s_0)=\qmod{K}{Z_K(s_0)}\ra \qmod{G}{G(s_0)}$$
is a diffeomorphism.

ii)  Endow $i\kg$ with an $\ad_K$-invariant inner product and consider the embedding
$M=\ad_K(s_0)\hookrightarrow i\kg$ of the projective algebraic manifold $M$  in the
Euclidean space $i\kg$ defined by the identification
$M=\ad_K(s_0)$   in i).   Then
\begin{enumerate}
\item For any $s\in \ad_K(s_0)=M$, the images of the  morphisms
$$T^{10}_s(M)\ra (i\kg)\otimes\C=\g,\ T^{01}_s(M)\ra (i\kg)\otimes\C=\g$$
induced by this embedding are $\g^\pm_{s}$, where
$$\g^\pm_{s}:=\bigoplus\limits_{\pm\lambda>0}{\rm Eig}([s,\cdot],\lambda)\ .
$$
\item  The induced  Riemannian metric  on $M$ is K\"ahler and homogeneous.

\end{enumerate}
\end{pr}
\pf i)  We just have to  prove that the natural map
$$\qmod{K}{Z_K(s_0)}\map \qmod{G}{G(s_0)}
$$
is   bijective.  For  surjectivity, let $g\in G$. We can decompose $g$ as $g=k u$, where
$u\in \exp(i\kg)$.  Let $B\subset G(s_0)\cap G_e$ be a Borel subgroup contained in
the parabolic subgroup $G(s_0)\cap G_e$ of the connected reductive group $G_e$ (the connected component of $e$ in
$G$). Since $\exp(i\kg)\in G_e$,  one can further decompose $u=l b$, where $l\in K_e$
and
$b\in B$, so one gets $g=(kl)b$, which shows that $g\equiv kl$ mod $G(s_0)$.

For   injectivity, it suffices to show that
\begin{equation}\label{intersection}
K\cap G(s_0)=Z_K(s_0)\ .
\end{equation}
The inclusion
$Z_K(s_0)\subset K\cap G(s_0) $ is obvious. The opposite inclusion follows by choosing an
embedding $G\hookrightarrow GL(n,\C)$ which maps $K$ into $U(n)$.  Let $\Phi$ be
the filtration  of $\C^n$ defined by the eigenspaces of the Hermitian matrix $s_0$. More
precisely, the
$k$-term of this filtration is the direct sum of the eigenspaces corresponding to the first
$k$ eigenvalues of $s_0$.   It is well known that the parabolic subgroup
$GL(n,\C)(s_0)\subset GL(n,\C)$ is just the stabilizer of the filtration $\Phi$. On the other
hand, one has
$G(s_0)=G\cap  GL(n,\C)(s_0)$. Let  $g\in  K\cap G(s_0)$. The matrix defined by $g$ is
unitary and leaves invariant the filtration $\Phi$, so it leaves invariant all the eigenspaces of
$s_0$ in $\C^n$.  Therefore $g$ commutes with $s_0$.

ii)

1. Let $s=\ad_k(s_0)\in \ad_K(s_0)$. The tangent space $T_s(\ad_K(s_0))$ is just the
subspace $\ad_k([\kg,s_0])=[\kg,s]\subset i\kg$. The point which corresponds to $s$ via
our identification
$\ad_K(s_0)=G/G(s_0)$ is  the congruence class $c:=kG(s_0)\in G/G(s_0)$. The tangent
space of
$G/G(s_0)$ at $c$ is
$$T_c(G/G(s_0))= \qmod{(L_k)_*(\g)}{(L_k)_*(\g(s_0))}\ .
$$
Now note that $\g^+_{s_0}$ is a complement of $\g(s_0)$ in $\g$, so one gets a
natural isomorphism $T_c(G/G(s_0))=(L_k)_*(\g^+_{s_0})=(R_k)_*(\g^+_s)$.

The
isomorphism
$T_s(\ad_K(s_0))\simeq T_c(G/G(s_0))$ induced by our  identification is given
by
$$[\kg,s]\map \g^+_s\textmap{(R_k)_*}(R_k)_*(\g^+_s)\ ,
$$
where the first arrow is just the restriction $\rho_s: [\kg,s]\ra \g^+_s$ of the
 projection $i\kg \ra\g^+_s$. The complex structure on the tangent space
$[\kg,s]$ is the pull-back of the natural complex structure of $\g^+_s$ via
$\rho_s$. Therefore,  the subspaces
$T^{10}_s(\ad_K(s_0))$, $T^{01}_s(\ad_K(s_0))$ of the complexified tangent space
$T_s^\C(\ad_K(s_0))=[\g,s]=\g^+_s\oplus\g^-_s$ are just $\g^+_s$, $\g^-_s$.

2. It is easy to see that the Riemannian metric induced by the embedding $\ad_K(s_0)$ is
Hermitian and homogeneous, i. e. invariant under the natural $K$-action. The fact that the
corresponding K\"ahler metric $\omega$ is closed follows from the identity
$$L_X \omega=d\ \iota_X\omega+ \iota_X d\omega
$$
applied to the fundamental vector fields associated with the transitive $K$-action on
$\ad_K(s_0)$.

\qed

Applying  Theorem \ref{meromorphy}, we get
\begin{pr}\label{AlmostHol} Let $\varphi\in L^2_1(B^n,i\kg)$ such that
$\varphi(x)\in\ad_K(s_0)$ for almost every $x\in  B^n$. Suppose that
$${\rm pr}_{\g^+_{\varphi(x)}}(\bar\partial_x \varphi)=0
$$
for almost every $x\in B^n$. Then $\varphi$ coincides almost everywhere with a
meromorphic map $B^n\ra M=\ad_K(s_0)$.

More generally, let ${\cal Q}$ be a holomorphic $G$-bundle on a complex manifold $X$, let
$P\subset {\cal Q}$ be a $K$-reduction, and let $\sigma\in L^2_1(P\times_{\ad} i\kg)$
be an
$L^2_1$ section  which defines an almost everywhere constant map $[\sigma]:X\ra
i\kg/K$ such that
\begin{equation}\label{WHC}
{\rm pr}_{V^+_{\sigma (x)}}(\bar\partial_x \sigma)=0
\end{equation}
for almost every $x\in X$.  Fix a representative $\sigma_0$ in the conjugacy class defined
by $\sigma$.  Then $\sigma$ defines a meromorphic reduction of ${\cal Q}$ to
$G(\sigma_0)$.
\end{pr}

In the regular case, the holomorphic map $B^n\ra G/G(s_0)$ associated with a  smooth
map $\varphi:B^n\ra \ad_K(s_0)$ satisfying the  holomorphy condition
$${\rm pr}_{V^+_{\sigma (x)}}(\bar\partial_x \sigma)=0\ \forall  x\in B^n
$$
is given explicitly by  the formulae
$$x\mapsto   \hbox{the class mod  $G(s_0)$  of an element  $k\in K$  with}\
  \ad_k(s_0)=\varphi(x)\ ,$$
$$x\mapsto   \{g\in G|\  \ad_{g^{-1}}(\varphi(x))\in\g(s_0)\}\ .$$

Similarly, the $G(s_0)$-reduction defined by a regular section  $\sigma\in
A^0(P\times_{\ad} i\kg)$ satisfying our weak holomorphy condition is
$${\cal Q}(\sigma)=\{q\in{\cal Q}|\ \sigma(q)\in \g(\sigma_0)\}\ .
$$

The converse of Proposition \ref{AlmostHol} is also true; more precisely
\begin{pr} Let $\rho\subset {\cal Q}/G(\sigma_0)$ be a
meromorphic $G(\sigma_0)$-reduction   of a
holomorphic principal
$G$-bundle
${\cal Q}$ on a
\it{compact} complex manifold $X$, and let $P\subset{\cal Q}$ be any $K$-reduction.

Then the  section $\sigma$ associated with $\sigma_0$ via the  $K\cap
G(\sigma_0)=Z_K(\sigma_0)$-reduction
$P\cap {\cal Q}^\rho$ of $P|_{X_\rho}$ defines an element   in
$L^2_1(P\times_{\ad} i\kg)$ which   satisfies the weak holomorphy condition
(\ref{WHC}) and whose associated meromorphic $G(\sigma_0)$-reduction is $\rho$.
\end{pr}
\pf  The crucial point here is the fact that $\sigma\in L^2_1(P\times_{\ad} i\kg)$.
Let $A$ be the Chern connection of $P$.  We make use of the formula (\ref{Kob}) which
was obtained in the proof of  Theorem
\ref{SimImp}.
\begin{equation} \int\limits_{X_\rho}  \left\langle  i \Lambda_g F_{
A},\sigma)\right\rangle vol_g= \frac{2\pi}{(n-1)!}\deg(\rho,
 {h(\xi)})- \sum_{\lambda\in {\rm Spec}^+
 {[ \xi ,\cdot] }}\lambda\nr a_\lambda\nr^2_{L^2}\ .
\end{equation}

On the other hand, by   formula  (\ref{d0sigma}), one has on $X_\rho$
\begin{equation} \partial_0 (\sigma|_{X_\rho}) = \sum_{\lambda\in {\rm Spec}^+
 {[\sigma_0,\cdot] }} [a _\lambda,(\sigma|_{X_\rho}) ]=-\sum_{\lambda\in {\rm
Spec}^+
 {[\sigma_0,\cdot] }} \lambda  a _\lambda \ .
\end{equation}

Comparing  these formulae, one gets immediately
$$\nr  \partial_0 (\sigma|_{X_\rho})\nr_{L^2}<\infty\ .
$$
On the other hand, one can  easily check  that the form $\partial_0 (\sigma
|_{X_\rho}))$, regarded as an $L^2$ - form on the whole  manifold $X$, is in fact the
distribution
$\partial_0$ - derivative of the $L^\infty$ - section $\sigma$. The proof uses the same
argument as for the analogous result for  orthogonal projections on subsheaves of
Hermitian holomorphic bundles.
\qed

\vspace{5mm}
{\small Authors' addresses: \vspace{2mm}\\
Martin L\"ubke: \\
Mathematical Institute,
Leiden University,
P.O. Box 9512
NL, 2300 RA Leiden,\\
e-mail: lubke@math.leidenuniv.nl
\\ \\
Andrei Teleman: \\
CMI,   Universit\'e de Provence,  39  Rue F. Joliot-Curie, F-13453
Marseille Cedex 13,   e-mail: teleman@cmi.univ-mrs.fr
}

 \end{document}